\documentclass[10pt]{amsart}

\usepackage[dvips]{graphicx}

\usepackage{color}
\usepackage{amscd}
\usepackage{amssymb}
\usepackage{amsthm}
\usepackage{amsmath}
\usepackage{latexsym}
\usepackage{subfigure}
\usepackage{euscript}

\usepackage{upref}

\usepackage[colorlinks,backref]{hyperref}                   

\newcommand{\down}{\downarrow}

\newcommand{\eps}{\varepsilon}

\newcommand{\E}{\mathcal{E}}

\newcommand{\D}{\mathcal{D}}

\newcommand{\Div}{\mathrm{div}}
\newcommand{\Grad}{\mathrm{\nabla}}

\newcommand{\R}{\mathbf{R}}

\newcommand{\dx}{\, dx}

\newcommand{\ds}{\, ds}
\newcommand{\dy}{\, dy}

\newcommand{\dal}{\, d\alpha}

\newcommand{\dt}{\, dt}
\newcommand{\dalpha}{\, d\alpha}
\newcommand{\dbeta}{\, d\beta}

\newcommand{\sgn}[1]{\mathrm{sign}^{+}\left(#1\right)}

\newcommand{\dsgne}[1]{\Bigl(\mathrm{sign}^+_{\eps}\Bigl)'\left(#1\right)}

\newcommand{\sgnm}[1]{\mathrm{sign}^{-}\left(#1\right)}

\newcommand{\Const}[1]{\mathrm{Const}_{#1}}

\newcommand{\norm}[1]{\left\|#1\right\|}

\newcommand{\abs}[1]{\left| #1 \right|}

\newcommand{\pt}{\ensuremath{\partial_t}}

\newcommand{\ps}{\ensuremath{\partial_s}}

\newcommand{\dto}{\ensuremath{\downarrow}}

\newcommand{\Do}{\ensuremath{Q_T}}
\newcommand{\Om}{\ensuremath{\Omega}}
\newcommand{\pOm}{\ensuremath{\partial\Omega}}
\newcommand{\oOm}{\ensuremath{\overline{\Omega}}}

\newcommand{\N}{{\mathbb N}}

\newcommand\bel{\begin{equation}\label}
\newcommand\ee{\end{equation}}
\newcommand\ba{\begin{array}}
\newcommand\bal{\begin{array}{l}}
\newcommand\ea{\end{array}}
\def\dsp{\displaystyle}

\newtheorem{theo}{Theorem}[section]
\newtheorem{prop}{Proposition}[section]
\newtheorem{lem}{Lemma}[section]
\newtheorem{rem}{Remark}[section]
\newtheorem{defi}{Definition}[section]
\newtheorem{cor}{Corollary}[section]

\def\ptl{{\partial}}
\newcommand{\sign}{{\rm sign}}

\newcommand{\grad}{{\nabla}}
\def\dist{{\rm dist\,}}
\def\diam{{\rm diam\,}}

\newcommand\1{\hbox{\hbox{1}\kern-.3em I}}

\def\dsp{\displaystyle}
\def\div{{\rm div\,}}
\def\Div{{\rm div\, }}
\def\grad{{\,\nabla }}

\def\dist{{\rm dist\,}}

\def\ph{{\varphi}}
\def\al{{\alpha}}
\def\eps{{\varepsilon}}
\def\Om{{\Omega}}
\def\ptl{{\partial}}

\def\ba{\begin{array}}
\def\bal{\begin{array}{l}}
\def\ea{\end{array}}
\def\be{\begin{equation}}
\def\bel{\begin{equation}\label}
\def\ee{\end{equation}}
\def\iint{\int\!\!\int}

\def\sign{{\rm sign\,}}

\def\K{{\scriptstyle K}}
\def\L{{\scriptstyle L}}
\def\dK{{\scriptstyle K^*}}
\def\dL{{\scriptstyle L^*}}
\def\dM{{\scriptstyle M^*}}
\def\ptK{{\scriptscriptstyle K}}
\def\ptL{{\scriptscriptstyle L}}
\def\ptdK{{\scriptscriptstyle K^*}}
\def\ptdL{{\scriptscriptstyle L^*}}
\def\ptdM{{\scriptscriptstyle M^*}}
\def\D{{D}}
\def\ptD{{\scriptstyle D}}

\def\Tau{{\boldsymbol{\mathfrak T}}}
\def\Drond{{\boldsymbol{\mathfrak D}}}
\def\Frond{{\mathcal F}}
\def\Grond{{\mathcal G}}

\def\Erond{{\mathcal E}}
\def\dErond{{{\mathcal E}^*}}
\def\Nrond{{\scriptstyle \mathcal N}}
\def\dNrond{{{\scriptstyle \mathcal N}^*}}
\def\Srond{{\boldsymbol{\mathfrak S}}}
\def\Mrond{{\boldsymbol{\mathfrak M}}}
\def\Prond{{{\mathcal P}}}
\def\dMrond{{\boldsymbol{\mathfrak M}^*}}
\def\ptTau{{\scriptscriptstyle \boldsymbol{\mathfrak T}}}
\def\ptDrond{{\scriptscriptstyle \boldsymbol{\mathfrak D}}}

\def\ptErond{{\scriptscriptstyle \mathcal E}}
\def\ptdErond{{{\scriptscriptstyle \mathcal E}^*}}
\def\ptNrond{{\scriptscriptstyle \mathcal N}}
\def\ptdNrond{{{\scriptscriptstyle \mathcal N}^*}}
\def\ptSrond{{\scriptscriptstyle \boldsymbol{\mathfrak S}}}
\def\ptMrond{{\scriptscriptstyle \boldsymbol{\mathfrak M}}}

\def\ptdMrond{{\scriptscriptstyle \boldsymbol{\mathfrak M}^*}}

\def\uTau{{u^\ptTau}}
\def\ubarTau{{u^{\overline{\ptTau}}}}

\def\uMrond{{u^\ptMrond}}

\def\udMrond{{u^\ptdMrond}}

\def\ptbarTau{{\overline{\ptTau}}}

\def\xK{{x_\ptK}}
\def\xL{{x_\ptL}}
\def\xdK{{x_\ptdK}}
\def\xdL{{x_\ptdL}}
\def\xdM{{x_\ptdM}}
\def\uK{{u_\ptK}}
\def\uL{{u_\ptL}}
\def\udK{{u_\ptdK}}
\def\udL{{u_\ptdL}}

\def\uKn{{u_\ptK^n}}

\def\udKn{{u_\ptdK^n}}

\def\wK{{w_\ptK}}
\def\wL{{w_\ptL}}
\def\wdK{{w_\ptdK}}
\def\wdL{{w_\ptdL}}
\def\wdM{{w_\ptdM}}

\def\edge{{\sigma}}
\def\dedge{{{\sigma}^*}}

\def\ptKL{{\ptK\!|\!\ptL}}
\def\ptdKdL{{\ptdK\!|\!\ptdL}}

\def\pt{\ensuremath{\partial_t}}

\def\D{{\mathcal D}}
\def\I{{\!|\!}}
\def\ptI{{{\scriptscriptstyle\!|\!}}}

\def\edgeS{{\edge_\ptS}}
\def\dedgeS{{\edge^*_\ptS}}
\def\ProjD{{\mathrm{Proj}_\ptD}}
\def\dProjD{{\mathrm{Proj}^*_\ptD}}
\def\size{{\text{size}}}
\def\reg{{\mathrm{reg}}}

\def\dKL{{d_{\ptK\!\ptL}}}
\def\ddKdL{{d_{\ptdK\!\ptdL}}}

\def\ddMdK{{d_{\ptdM\!\ptdK}}}
\def\nuKL{{\nu_\ptKL}}
\def\nudKdL{{\nu_\ptdKdL}}
\def\nudLdM{{\nu_{\ptdL\!\!,\ptdM}}}
\def\nudMdK{{\nu_{\ptdM\!\!,\ptdK}}}
\def\nuK{{\nu_{\ptK}}}
\def\nudK{{\nu_{\ptdK}}}
\def\nuS{{\nu_\ptS}}
\def\nudS{{\nu^*_\ptS}}

\def\DM{{\scriptstyle D}}
\def\SDM{{\scriptstyle S}}
\def\ptD{{\scriptscriptstyle D}}
\def\ptS{{\scriptscriptstyle S}}

\def\Vrond{{\scriptstyle\mathcal V}}
\def\dVrond{{{\scriptstyle\mathcal V}^*}}
\def\ptVrond{{\scriptscriptstyle\mathcal V}}
\def\ptdVrond{{{\scriptscriptstyle\mathcal V}^*}}

\def\ptptl{{\scriptscriptstyle \ptl}}
\def\KIL{{\K\I\L}}
\def\ptKL{{\ptK\!,\ptL}}
\def\ptdKdL{{\ptdK\!\!,\ptdL}}
\def\dKIdL{{\dK\I\dL}}
\def\ptKIL{{\ptK\ptI\ptL}}
\def\ptdKIdL{{\ptdK\ptI\ptdL}}

\def\mKIL{{m_\ptKIL}}
\def\mdKIdL{{m_\ptdKIdL}}
\def\mK{{m_\ptK}}
\def\mdK{{m_\ptdK}}
\def\mD{{m_\ptD}}

\def\mS{{m_\ptS}}

\def\Delt{{{\scriptstyle \Delta} t}}
\def\Del{{\scriptstyle \Delta}}
\def\diam{{\rm diam\,}}

\def\Bleft{{\biggl[\hspace*{-3pt}\biggl[}}
\def\Bright{{\biggr]\hspace*{-3pt}\biggr]}}

\def\Aleft{{\biggl\{\!\!\!\biggl\{}}
\def\Aright{{\biggr\}\!\!\!\biggr\}}}

\def\char{{1\!\mbox{\rm l}}}

\def\QKn{{Q_\ptK^n}}
\def\QdKn{{Q_\ptdK^n}}
\def\QDn{{Q_\ptD^n}}
\def\barDelt{{\,\overline{\!\Delt\!}\,}}

{%

\begin{enumerate}}%
{\end{enumerate}
}
\newenvironment{Definitions}
{%

\begin{enumerate}}%
{\end{enumerate}
}
\newenvironment{DefinitionsP}
{%

\begin{enumerate}}%
{\end{enumerate} }

\begin{document}

\title[FV for doubly nonlinear degenerate equations] {Discrete duality
finite volume schemes \\ for doubly nonlinear degenerate hyperbolic-parabolic equations}

\date{\today}

\author[B. Andreianov]{B. Andreianov}
\address[Boris Andreianov]{\newline
         Laboratoire de Math\'ematiques\newline
         Universit\'e de Franche-Comt\'e\newline
         16 route de Gray\newline
         25 030 Besan{\rm \c{c}}on
         Cedex, France} \email[]{boris.andreianov\@@univ-fcomte.fr}

\author[M. Bendahmane]{M. Bendahmane}
\address[Mostafa Bendahmane]{\newline
Departamento de Ingenier\'{\i}a Matem\'atica\newline Facultad
de Ciencias F\'{\i}sicas y Matem\'aticas\newline Universidad de
Concepci\'on\newline Casilla 160-C Concepci\'on, Chile} \email[]{
mostafab\@@ing-mat.udec.cl}

\author[K. H. Karlsen]{K. H. Karlsen}
\address[Kenneth H. Karlsen]{\newline
                Centre of Mathematics for Applications \newline
                 University of Oslo\newline
                 P.O. Box 1053, Blindern\newline
                 N--0316 Oslo, Norway\newline
                 and\newline
                Center for Biomedical Computing,\newline
                Simula Research Laboratory\newline
                P.O. Box 134\newline
                N--1325 Lysaker, Norway}
\email[]{kennethk@math.uio.no}
\urladdr{http://folk.uio.no/kennethk}



\thanks{The work of M. Bendahmane was supported by the FONDECYT project 1070682.
The work of K. H. Karlsen was supported by the Research Council of Norway through
an Outstanding Young Investigators Award. A part of this work was
done while B. Andreianov enjoyed the hospitality of the Centre of
Mathematics for Applications (CMA) at the University of Oslo, Norway.
This article was written as part of the the international research program
on Nonlinear Partial Differential Equations at the Centre for
Advanced Study at the Norwegian Academy of Science
and Letters in Oslo during the academic year 2008--09.}

\begin{abstract}
We consider a class of doubly nonlinear degenerate
hyperbolic-parabolic equations with homogeneous Dirichlet boundary
conditions, for which we first establish the existence and
uniqueness of entropy solutions. We then turn to the construction
and analysis of discrete duality finite volume schemes (in the
spirit of Domelevo and Omn\`es \cite{DomOmnes}) for these problems
in two and three spatial dimensions. We derive a series of
discrete duality formulas and entropy dissipation inequalities for
the schemes. We establish the existence of solutions to the
discrete problems, and prove that sequences of approximate
solutions generated by the discrete duality finite volume schemes
converge strongly to the entropy solution of the continuous
problem. The proof revolves around some basic a priori estimates,
the discrete duality features, Minty-Browder type arguments, and
``hyperbolic'' $L^\infty$ weak-$\star$ compactness arguments
(i.e., propagation of compactness along the lines of Tartar,
DiPerna, \ldots). Our results cover the case of non-Lipschitz
nonlinearities.
\end{abstract}

\maketitle

\tableofcontents

\section{Introduction}\label{sec:intro}
In this paper we consider degenerate hyperbolic-parabolic problems
of the form
\begin{equation}\label{eq:prob1}
    \begin{cases}
        \ptl_t u +\div {\mathfrak f}(u)
        -\div {\mathfrak a}(\grad A(u))=\EuScript S , & \text{in $Q:=(0,T)\times\Om$},
        \\ u|_{t=0}=u_{0}, & \text{in $\Omega$}
        \\ u=0, & \text{on $\Sigma=(0,T)\times\ptl\Om$},
    \end{cases}
\end{equation}
where $u:(t,x)\in Q\to\R$ is the unknown function, $T>0$ is a
fixed time, $\Omega\subset\R^d$ is a bounded domain with polygonal
boundary $\pOm$ and outward unit normal $n$. We consider the cases $d=2$ and $d=3$. The
initial data $u_0:\Om\to\R$ are assumed to be a bounded measurable
function, i.e.,
$$
u_0\in L^\infty(\Om),
$$
while the source $\EuScript S:Q\to \R$ is assumed to be a
measurable function for which ${\EuScript S}(t,\cdot)\in
L^\infty(\Om)$ for a.e.~$t\in (0,T)$ and $\int_0^T \|{\EuScript
S}(t,\cdot)\|_{L^\infty(\Om)}\dt<\infty$; we abusively denote it
by
\begin{equation}\label{eq:source-terms-space}
{\EuScript S}\in L^1(0,T;L^\infty(\Om)).
\end{equation}

The function  ${\mathfrak a}:\R^N\to \R^N$ is taken under the form
$$
{\mathfrak a}(\xi)=k(\xi)\xi,
$$
where $k$ is a scalar function. The function ${\mathfrak a}$ is assumed to be
continuous and strictly monotone. We assume that
there exist $p\in (1,+\infty)$ and $C>0$ such that
$$
\frac{1}{C}|\xi|^{p-2}\leq k(\xi)\leq C |\xi|^{p-2}, \qquad \forall \xi\in \R^d\setminus\{0\}.
$$
In particular, the associated operator $w\mapsto -\div\,
\Bigl(k(|\grad w|)\grad w\Bigr)$ is a Leray-Lions operator acting
from $W^{1,p}_0(\Om)$ to $W^{-1,p'}(\Om)$ with $p'=\frac{p}{p-1}$.
A prototype example is the $p$-laplacian, which corresponds to
$k(\xi)=|\xi|^{p-2}$.

We assume that the diffusion function $A(\cdot)$ satisfies
\begin{equation*}
    \text{$A(\cdot)$ is continuous and nondecreasing, normalized by $A(0)=0$},
\end{equation*}
while the convective flux function ${\mathfrak f}(\cdot)$ satisfies
$$
\text{${\mathfrak f}=({\mathfrak f}_1,\dots,{\mathfrak f}_d):Q\times\R\to\R^d$
is continuous and normalized by $\mathfrak{f}(0)=0$.}
$$
We emphasize that the fluxes ${\mathfrak f},A$ are
not necessarily locally Lipschitz continuous.

Problems more general than \eqref{eq:prob1}, for which our results
can be extended, will be discussed in Section \ref{sec:Extensions}.

The class (\ref{eq:prob1}) of nonlinear partial differential
equations includes several important particular cases. The
hyperbolic conservation law
$$
\pt u + \Div {\mathfrak f}(u) =0
$$
is a special case of \eqref{eq:prob1}. The celebrated theory of
$L^\infty$ entropy solutions for scalar conservation laws in
$\R^d$ was developed by Kruzhkov \cite{Kruzkov}, while the $BV$
theory was set up by Vol'pert \cite{Volpert}. The extensions for
the Dirichlet problem in bounded domains are due to Bardos,
LeRoux, N\'ed\'elec \cite{Bardos:BC} (for the $BV$ setting) and
Otto \cite{Otto:BC} (for the $L^\infty$ setting). Note that the
boundary condition is only verified in some generalized sense (see
\cite{Bardos:BC,Otto:BC,Malek,Carrillo,RouvreGagneux,Mascia_etal:2000,Vovelle:02,EGHMichel,MichelVovelle,AmmarCarWitt}).

Many other well-known partial differential equations (usually
possessing more regular solutions) are also special cases of
\eqref{eq:prob1}. Let us mention the heat and porous medium equations
$$
\pt u = \Delta u, \qquad
\pt u = \Delta u^m, \quad m>1,
$$
and more generally degenerate
convection-diffusion equations of the type
\begin{equation}
   \label{eq:ConvDiff}
   \pt u + \div {\mathfrak f}(u)  = \Delta A(u).
\end{equation}
Degenerate parabolic equations like \eqref{eq:ConvDiff} occur in
theories of flow in porous media (see discussion and references \cite{EspKar}) and
sedimentation-consolidation processes \cite{Burger:bok}.

As other famous representatives of the class
of equations that is considered herein, we
mention the $p$-Laplace equation
$$
\pt u = \Div \Bigl(\abs{\Grad u}^{p-2}
\Grad u\Bigr),\qquad p>1,
$$
which arises in the theory of non-Newtonian filtration.
Also well known is the more general polytropic filtration equation
$$
\pt u = \Div \left(\abs{\Grad \left(\abs{u}^{m-1}u\right)}^{p-2}
\Grad \left(\abs{u}^{m-1}u\right)\right), \qquad m,p>1.
$$

A related class of equations consists of the so-called elliptic-parabolic equations
$$
\pt b(v) = \Div\, {\mathfrak a}(v,\Grad v),
$$
where $b:\R\to\R$ is continuous nondecreasing, and
$\mathfrak{a}(r,\xi):\R\times\R^N\to \R^N$ gives rise to a
Leray-Lions operator. We refer to
\cite{AltLuckhaus,BenilanWittbold,Otto:L1_Contr,CarrilloWittbold:99,AmmarWittbold}
and the references cited therein for more information on
elliptic-parabolic equations.

A chief goal of this paper is to propose and analyze a specific
class of finite volume schemes for the problem \eqref{eq:prob1}.
Note that finite volume schemes are well suited for approximation
of equations in divergence form, such as \eqref{eq:prob1}.
Discretization of the aforementioned hyperbolic, porous medium,
convection-diffusion, and elliptic-parabolic equations by finite
volume methods is quite standard by now and often used in
engeneering practice. We refer to
\cite{EyGaHe:book,ClaireChainais,AfifAmaziane,EvjeKarlsen:DD,EvjeKarlsen:SJNA,KR:Rough_Diff,Herm,Ohlberger:FVM,Vovelle:02,EGHMichel,MichelVovelle,AGW,ABH-M2AN,ABH``double'',Droniou} 
and references therein for different convergence results and
numerical experiments. For related works on linear elliptic
problems, see
\cite{Aavatsmark-et-al,Aavatsmark,Herm,DomOmnes,AngotetCo,Herm-3D,EGH:IMAJNA,EGH:AnisoCRAS,EH-SUCCES,EGH:Sushi-IMAJNA}
and the discussion in Section~\ref{sec:Extensions}. 
Alternative numerical approaches have also been 
investigated; here we only mention finite element schemes (see
\cite{Chow,BarrettLiu} and references therein), kinetic schemes
(see \cite{A-D:Nat:Tang,BouGuaNat,GuarMilTer} and references
therein) and operator splitting schemes (see \cite{EspKar}).

Having said that, we are not aware of any papers that construct
convergent numerical schemes for mixed type equations of the
generality considered herein. Indeed, they combine a number of 
difficulties such as nonlinear convection, doubly nonlinear diffusion, strong 
degeneracy, and shocks, which in turn necessitates the use 
of a suitable framework of discontinuous
entropy solutions. Furthermore, in the absence of 
the Lipschitz continuity assumption on the 
convective flux $\mathfrak f(\cdot)$, the CFL condition
does not make sense; therefore we have to discretize the
convective term with a time-implicit scheme.

We begin by providing the entropy solution framework for \eqref{eq:prob1}; this
is the topic of Section \ref{sec:def} and Appendix A. Due to the nonlinearity of
${\mathfrak f}(\cdot)$ and the possible degeneracy of $A(\cdot)$,
the problem \eqref{eq:prob1} will in general possess shock wave
solutions, a feature that can reflect the physical phenomenon of
breaking of waves.  This is well known in the context of
conservation laws. Also the boundary condition cannot be
prescribed pointwise on the whole boundary $\Sigma$ when $A$ is
not strictly increasing.  Due to this loss of regularity, it is necessary to work
with weak solutions; moreover, to single out a physically relevant
and unique weak solution, we need to
impose additional ``entropy inequalities'', in the spirit of
Kruzhkov \cite{Kruzkov}. Early results on hyperbolic-parabolic
equations were obtained by Volpert, Hudjaev \cite{VolHud}; see
also \cite{WuYin,Yin:Double,WZYL:2001}, \cite{Simondon},
\cite{Carrillo:94}, \cite{BenilanToureIII} and references cited therein, and
\cite{Tassa}, \cite{CockGripen}, \cite{ChenDiBen}.
$L^1$ entropy techniques for degenerate convection-diffusion
equations like \eqref{eq:ConvDiff}, which take into account
both hyperbolic and parabolic features, were
developed by Carrillo \cite{Carrillo} for the homogeneous
Dirichlet problem in bounded domains. Since then, many authors
extended the Carrillo results in various directions (see e.g.
\cite{CarrilloWittbold:99,Igbida-Urbano,Mascia_etal:2000,
RouvreGagneux,BurgerEvjeKarlsen:1D_IBVP,EKR:Hyp2000,KO:Unique,KR:Rough_Unique,MichelVovelle,EGHMichel,AmmarCarWitt,
AndrIgbida}). Some additional techniques are required for anisotropic
diffusion problems, where a kinetic approach (see Chen, Perthame
\cite{ChenPerthame}) and an accurate entropic
approach (see Bendahmane, Karlsen \cite{BenKar:Renorm,BenKar:RenormII}) were developed in the few last
years; see also Souganidis, Perthame \cite{SougPerth:03} and
and Chen, Karlsen \cite{ChenKarlsen}. In this paper, we use a
variant of the Carrillo entropic approach.
Following the Tartar-DiPerna idea of measure-valued solutions
and using the techniques of Eymard, Gallou\"et, Herbin
\cite{EyGaHe:book}, we introduce a notion of entropy process
solution for \eqref{eq:prob1}, and establish the related
identification and uniqueness results. More exactly, we work with
entropy double-process solutions arising in the particular context
of discrete duality finite volume schemes. In Section \ref{sec:def}, we
show the existence result for \eqref{eq:prob1} and state
uniqueness; an adaptation of the standard uniqueness (and, more
generally, $L^1$ contraction and comparison principle)
proof is given in Appendix A.

In Section~\ref{sec:FV} we construct discrete duality finite
volume (DDFV) schemes for \eqref{eq:prob1} in two and three
spatial dimensions (some other schemes are briefly discussed in
Section~\ref{sec:Extensions}). We adapt the approximations used by
Eymard, Gallou\"et, Herbin \cite{EyGaHe:book} (see also
\cite{ClaireChainais,Vovelle:02,EGHMichel,MichelVovelle}) for the
nonlinear convection term, and those used by Hermeline
\cite{Herm,Herm-3D}, Domelevo, Omn\`es \cite{DomOmnes} and
Andreianov, Boyer, Hubert \cite{ABH``double''} for the doubly
nonlinear diffusion term. In 3D, we propose new DDFV schemes that
possess convenient discrete duality properties.

Our 3D scheme is a very particular case of the schemes introduced and
studied numerically by Hermeline in \cite{Herm-3D}. In passing, we
mention that different kinds of 3D discrete duality schemes were
constructed in \cite{Pierre,CoudierePierre} and in
\cite{CoudiereHubertPierre}. Appendix B (see also
\cite{ABK-FVCA5,AndrBend}) is devoted to an elementary
reconstruction lemma which underlies our DDFV schemes in 3D. In
contrast to \cite{DomOmnes,ABH``double''}, we are led to penalize
our DDFV schemes to ensure that the two approximations of $A(u)$
actually converge to the same limit (see Section
\ref{sec:Penalization}). The DDFV schemes constructed in Section
\ref{sec:FV} possess several convenient discrete calculus formulas
that we collect in Section \ref{sec:DiscreteDuality}. Related
consistency estimates and properties of the associated spaces of
discrete functions are given in Section \ref{sec:DiscreteSpaces}.
The (few) available {\it a priori} estimates for the discrete
solutions are collected in In Section \ref{sec:DiscrSol}. In the
same section, the existence of discrete solutions is shown.
Furthermore, we establish that, up to an error term in the
equation depending on the discretization parameter, discrete
solutions can be considered as entropy solutions of
\eqref{eq:prob1}. In Section \ref{sec:convergence} we prove that
discrete solutions converge, as the discretization parameter tends
to zero, to an entropy double-process solution that turns out to
be the (unique) entropy solution of \eqref{eq:prob1}.  It should
emphasized that we obtain strong convergence of both convective
and diffusive fluxes, in spite  of the double nonlinearity of the
problem \eqref{eq:prob1}. Section \ref{sec:Extensions} contains
references to some known finite volume schemes for nonlinear
diffusion-convection equations, and discusses the extension of our
results to different generalizations of problem \eqref{eq:prob1}.

\section{Notions of solution and well-posedness}
\label{sec:def}
As it was explained in the introduction, we need the notion of
weak solution for \eqref{eq:prob1} with additional ``entropy''
conditions. In order to use entropy conditions in the interior of
$Q$ and, moreover, take into account the homogeneous Dirichlet
boundary condition on $\Sigma$, following Carrillo \cite{Carrillo}
we will work with the so-called ``semi-Kruzhkov'' entropy-entropy
flux pairs $(\eta^\pm_c,{\mathfrak q}^\pm_c)$ for each $c \in \R$; they are defined as
$$
\eta^{+}_c (z)=(z-c)^{+}, \qquad \eta^{-}_c (z)=(z-c)^{-},
$$
$$
{\mathfrak q}^{+}_c(z)=\sgn{z-c} ({\mathfrak f}(z)-{\mathfrak
f}(c)), \qquad {\mathfrak q}^{-}_c(z)=\sgnm{z-c} ({\mathfrak
f}(z)-{\mathfrak f}(c)).
$$
By convention, we assign $(\eta^\pm_c)'(c)$ to be zero. Here
$(z-c)^\pm$ denote the nonnegative quantities satisfying
$z-c=(z-c)^+ -(z-c)^-$; moreover, we use the notation
\begin{align*}
    \sgn{z\!-\!c} & =(\eta^+_c)'(z)=
    \begin{cases}
        1,& z>c \\
        0, &  z\leq c,
    \end{cases}\\
    \sgnm{z\!-\!c} & =(\eta^-_c)'(z)=
    \begin{cases}
        0,& z\geq c, \\
        -1,& z< c.
    \end{cases}
\end{align*}
At certain points, we will also need smooth regularizations of the
semi-Kruzhkov entropy-entropy flux pairs; it is sufficient to
consider regular ``boundary'' entropy pairs
$(\eta_{c,\eps}^\pm,\mathfrak{q}_{c,\eps}^\pm)$ (cf. Otto
\cite{Otto:BC} and the book \cite{Malek}), which are
$W^{2,\infty}$ pairs with the same support as
$(\eta_c^\pm,\mathfrak{q}_c^\pm)$, converging pointwise to
$(\eta_c^\pm,\mathfrak{q}_c^\pm)$ as $\eps\to 0$. Specifically,
the functions
$$
\sign^+_\eps(z)= \frac{1}{\eps}\min\{z^+,\eps\}, \quad
\sign^-_\eps(z)=\frac{1}{\eps}\max\{-z^-,-\eps\}
$$
will be used to approximate
$\sign^\pm(\cdot)=(\eta_0^\pm)'(\cdot)$.

In view of the monotonicity of $A:\R\to\R$, the following
definition is meaningful.

\begin{defi}\label{AThetaDef}
For any locally bounded piecewise continuous function $\theta:\R\to \R$, we
define (using, e.g., the Stieltjes integral) the function $A_\theta:\R\to\R$ by
\begin{equation}\label{AthetaDef}
    A_\theta(z)= \int_0^z \theta(s)\,dA(s).
\end{equation}
\end{defi}

The ensuing lemma shows that there exists a continuous function $\widetilde A_\theta$
such that $A_\theta(z)=\widetilde A_\theta(A(z))$.
We prove this lemma under rather strong assumptions, but
they are still sufficient for our needs.

\begin{lem}\label{FunctionalDepLemma}
(i) Let $\theta,A_\theta$ be a couple of functions as
introduced in Definition \ref{AThetaDef}.
Then there exists a continuous function
$\widetilde A_\theta:A(\R)\to \R$ such that
$$
A_\theta(z)=\widetilde A_\theta(A(z)), \qquad \forall z\in \R.
$$
Moreover, $\widetilde A_\theta$ is Lipschitz continuous.

(ii) Assume additionally that $\theta \in W^{1,\infty}(\R)$, and let
$(A^\rho)_\rho$ be a sequence of nondecreasing continuous
surjective functions converging to $A$ pointwise on $\R$ as
$\rho\to 0$. Define $\widetilde A^\rho_\theta,A^\rho_\theta$  by
(i) and \eqref{AthetaDef} with $A^\rho$ replacing $A$.
Then $\widetilde A^\rho_\theta$ converges to $\widetilde A_\theta$
uniformly on compact subsets of $A(\R)$.
\end{lem}

\begin{proof}
(i) For $b\in A(\R)$, we can define $\widetilde A_\theta$ by
$\widetilde A_\theta(b)=A_\theta(z)$ for some
$z\in A^{-1}(b)$. If $A(z)=A(\hat z)$, then the measure $dA(s)$ vanishes between
$z$ and $\hat z$; thus
$$
A_\theta(z)-A_\theta(\hat z)=\int_{\hat z}^{z} \theta(s)\,dA(s)=0,
$$
and $\widetilde A_\theta$ is well-defined. For all $b,\hat b\in A(\R)$,
$$
\widetilde A_\theta(b)-\widetilde A_\theta(\hat b)
=A_\theta(z)-A_\theta(\hat z)=\int_{\hat z}^z \theta(s)\,dA(s),
\quad z\in A^{-1}(b), \hat z \in A^{-1}(\hat b).
$$
Consequently,
$$
\abs{\widetilde A_\theta(b)-\widetilde A_\theta(\hat b)}
\leq \norm{\theta}_{L^\infty} \abs{A(z)-A(\hat z)}
=\norm{\theta}_{L^\infty}\abs{b-\hat b}.
$$

(ii) Since the functions $\widetilde A^\rho_\theta$ are monotone, by
the Dini theorem it is sufficient to prove the pointwise convergence.
By the same argument, the convergence of $A^\rho$ to
$A$ is actually uniform on compact subsets of $\R$. Take $b\in A(\R)$ and
$z\in A^{-1}(b)$. Set $b^\rho=A^\rho(z)$; we have
$b^\rho\to b$ as $\rho\to 0$. Using (i) and the integration-by-parts formula
for the Stieltjes integral, we get
\begin{align*}
    \abs{\widetilde A^\rho_\theta(b)-\widetilde A_\theta(b)}
    & \leq \abs{\widetilde A^\rho_\theta(b)-\widetilde A^\rho_\theta(b^\rho)}
    +\abs{\widetilde A^\rho_\theta(b^\rho)-\widetilde A_\theta(b)} \\
    & \leq \norm{\theta}_{L^\infty}\abs{b-b^\rho}+\abs{\int_0^z \theta(s)\,d(A^\rho(s)-A(s))} \\
    & \leq 2\norm{\theta}_{L^\infty}\abs{b-b^\rho}+\abs{\int_0^z(A^\rho(s)-A(s))\theta'(s)\,ds}.
\end{align*}
The right-hand side converges to zero as $\rho\to 0$. Thus the claim follows.
\end{proof}

We have now came to the definition of an entropy solution. Here
and in the sequel,  $\R^{\pm}$ denote $\{k\in\R\,|\, \pm k \geq
0\}$, respectively.

\begin{defi}[entropy solution]\label{def:Entropy-Solution}
An entropy solution of the initial-boundary value problem \eqref{eq:prob1}
is a measurable function $u:\Do\to \R$ satisfying
\begin{Definitions}\label{def:entropy2}
    \item $u\in L^\infty(Q)$ and $w=A(u)\in L^p(0,T;W^{1,p}_0(\Om))$;

    \item for all $\psi\in \D([0,T)\times \Om)$,
    \begin{equation*}
        \begin{split}
            &\int_{Q} \Biggl(u\pt \psi+ {\mathfrak f}(u) \cdot \grad\psi
            -k(\Grad w) \Grad  w\cdot\grad \psi \Biggl) \dx\dt
            \\ & \qquad\quad +\int_\Omega u_0\psi(0,\cdot)\dx
            +\int_Q  {\EuScript S}\psi \dx\dt =0;
        \end{split}
    \end{equation*}
    \item for all pairs $(c,\psi)\in \R^\pm\times\D([0,T)\times \overline{\Om})$, $\psi \geq 0$,
    and also for all pairs $(c,\psi)\in \R\times\D([0,T)\times \Omega)$, $\psi \geq 0$,
    \begin{equation*}
        \begin{split}
            &\int_{Q} \Biggl(\eta_c^\pm(u)\pt \psi
            +{\mathfrak q}_c^\pm(u) \cdot \grad\psi
            -k(\Grad w) \Grad \widetilde A_{(\eta_c^\pm)'}(w) \cdot
            \grad \psi \Biggl) \dx\dt
            \\ &  \qquad\qquad
            +\int_\Omega \eta_c^\pm(u_0)\psi(0,\cdot)\dx
            + \int_Q (\eta_c^\pm)'(u)\, {\EuScript S}\psi \dx\dt  \ge 0.
        \end{split}
    \end{equation*}
    \label{def:entropy2-KHK}
\end{Definitions}
\end{defi}

For the convergence proof we need the notion of entropy
double-process solutions; we adapt
this notion from \cite{EyGaHe:book,ClaireChainais,GaHu,EGHMichel}, where
entropy process solutions have been introduced
for hyperbolic problems and degenerate parabolic problems with linear diffusion.
This definition is based upon the so-called
``nonlinear $L^\infty$ weak-$\star$ convergence'' property, which is well-known
in the equivalent framework of measure-valued solutions
developed earlier by Tartar and DiPerna:
\begin{equation}\label{nonlinweakstar}
    \left |
    \begin{array}{l}
        \text{each sequence $(u_\rho)$ bounded in
        $L^\infty(Q)$ admits a subsequence such that}\\
        \text{$\forall F\in C(\R)$, $\quad F(u_\rho(\cdot,\cdot))
        \to {\dsp\int_0^1} F(\mu(\cdot,\cdot,\alpha))\;d\alpha$
        in $L^\infty(Q)$ weak-$\star$},
    \end{array}\right.
\end{equation}
where the function $\mu\in L^\infty(Q\times(0,1))$ is referred to as
the ``process function''; it is related to the distribution function of the Young measure.
As usual, in \eqref{nonlinweakstar} and elsewhere we do
not bother to (re)label sequences.

We remark that the reason for introducing in the definition below two different
process functions $\mu,\mu^*$, both corresponding to the single unknown
function $u$, is that it permits us to handle the double approximation of $u$ by
pairs $\uMrond,\udMrond$ in the framework
of DDFV schemes (see Section \ref{sec:FV}).

\begin{defi}[entropy double-process solution]\label{def:Entropy-Solutionp}
A triplet $(\mu,\mu^*,w)$ of measurable functions, with
$\mu,\mu^*:Q\times (0,1) \to \R$ and $w:Q\to\R$, is called an
entropy double-process solution of the initial-boundary
value problem \eqref{eq:prob1} if the following conditions are met:
\begin{DefinitionsP}
    \item $\mu,\mu^* \in L^\infty(Q\times(0,1))$,
    $w\in L^p(0,T;W^{1,p}_0(\Om))$, and
    $$
    A(\mu(t,x,\alpha))\equiv w(t,x)\equiv A(\mu^*(t,x,\alpha)),
    $$
    for a.e.~$(t,x,\alpha)\in Q\times (0,1)$.
    \label{def:entropy1p}

    \item For all $\psi\in \D([0,T)\times \Om)$,
    \begin{equation*}
        \begin{split}
            &\int_0^1 \int_{Q} \Biggl(\frac{1}{d}\Bigl(\mu+(d \!-\! 1)\mu^*\Bigr)  \pt \psi
            + \frac{1}{d}\Bigl({\mathfrak f}(\mu)
            +(d \!-\! 1){\mathfrak f}(\mu^*)\Bigr) \cdot \grad\psi\Biggl)\dx\dt \dal
            \\ & \qquad - \int_Q k(\Grad w) \Grad  w \cdot \grad \psi \dx\dt
            +\int_\Omega u_0\psi(0,\cdot)\dx + \int_Q  {\EuScript S}\psi \dx\dt =0.
        \end{split}
    \end{equation*}\label{def:entropyWeakp}
    \item For all pairs $(c,\psi)\in \R^\pm\times\D([0,T)\times \overline{\Om})$, $\psi \geq 0$,
    and also for all pairs $(c,\psi)\in \R\times\D([0,T)\times \Omega)$, $\psi \geq 0$,
    \begin{equation*}
        \begin{split}
            &\int_0^1\!\!\!\int_{Q} \Biggl(\frac{1}{d}\Bigl(\eta^\pm_c(\mu)
            +(d \!-\! 1)\eta^\pm_c(\mu^*)\Bigr)\pt \psi
            +\frac{1}{d}\Bigl({\mathfrak q}_c^\pm(\mu)+(d \!-\! 1){\mathfrak q}_c^\pm(\mu^*)\Bigr) \cdot
            \grad\psi\Biggl) \dx\dt \dal
            \\ & \qquad -\int_{Q} k(\grad w) \Grad \widetilde A_{(\eta_c^\pm)'}(w) \cdot\grad \psi \dx\dt
            +\int_\Omega \eta^\pm_c(u_0)\psi(0,\cdot)\dx
            \\ & \qquad\qquad
            + \int_0^1\!\!\!\int_Q \frac{1}{d}\Bigl((\eta^\pm_c)'(\mu)
            +(d \!-\! 1)(\eta^\pm_c)'(\mu^*)\Bigr)\,{\EuScript S}\psi\dx\dt\dal \ge 0.
        \end{split}
    \end{equation*}
    \label{def:entropy2p}
\end{DefinitionsP}
\end{defi}

\begin{rem}\normalfont \label{RemTwoFormsOfDiffusion}
Since $\grad w=0$ a.e.~on $\{(t,x)\in Q)\,|\, w(t,x)=A(c)\}$ for
any $c\in \R$, the term $k(\grad w) \Grad \widetilde
A_{(\eta_c^\pm)'}(w)$ in the above definitions can be rewritten as
\begin{equation}\label{eq:equivalent-KHK}
    \text{$(\eta_c^\pm)'(z)\mathfrak{a}(\grad w)$ for any $z\in A^{-1}(w)$, and
    also as $\sign^\pm(w-A(c)) \mathfrak{a}(\grad w)$.}
\end{equation}
The form used in \ref{def:entropy2-KHK} and \ref{def:entropy2p} is convenient for
expressing the approximate entropy inequalities at the discrete
level; the equivalent form \eqref{eq:equivalent-KHK} is used in the uniqueness proof.
Both forms are exploited in the existence proof below.
\end{rem}

\begin{rem}\normalfont \label{rem:EPS-ES}
Let $u$ be an entropy solution of \eqref{eq:prob1}. Then the
triplet $(\mu,\mu^*,w)$ defined by
\begin{align*}
    &\text{$\mu(t,x,\al)=\mu^*(t,x,\al)=u(t,x)$ for a.e.~$(t,x,\al)\in Q\times (0,1)$},\\
    &\text{$w(t,x)=A(u(t,x))$ for a.e.~$(t,x)\in Q$}.
\end{align*}
is an entropy double-process solution of \eqref{eq:prob1}.

Conversely, if $(\mu,\mu^*,w)$ is an entropy double-process solution of
\eqref{eq:prob1} for which $\mu(t,x,\al)=\mu^*(t,x,\al)=u(t,x)$
a.e.~on $Q\times (0,1)$ for some function $u:Q\to \R$, then this
$u$ is an entropy solution of \eqref{eq:prob1}.
\end{rem}

Note that in Definition~\ref{def:Entropy-Solution}, we have only considered
$\al$-independent data $u_0,f$. In this case, the notion of
entropy double-process solution is just a technical tool that permits to bypass the lack
of strong compactness of sequences of approximate solutions. As a first illustration of this, we
pass to the limit in vanishing viscosity approximations (without $BV$ estimates)
to prove the existence of an entropy double-process
solution such that $\mu\equiv \mu^*$.

\begin{theo}\label{th:EPS-ESexists}
Under the assumptions stated in Section \ref{sec:intro}, there exists an entropy double-process solution to
the initial-boundary value problem \eqref{eq:prob1} for which $\mu\equiv \mu^*$.
\end{theo}
Notice that the above result holds for any Lipschitz domain $\Om$
in any space dimension. In passing, we also mention that the
existence result of Theorem~\ref{th:EPS-ESexists} has recently
been generalized by Ouaro and the authors \cite{ABKO} to the
case of a triply nonlinear degenerate diffusion equation.

\begin{proof}
The proof is divided into several steps.

\textbf{(i)} We approximate problem \eqref{eq:prob1} by
regular problems $\eqref{eq:prob1}_\rho$ with ${\mathfrak f},A$
replaced by ${\mathfrak f}_\rho, A_\rho$ such that ${\mathfrak f}_\rho,A_\rho,[A_\rho]^{-1}$
are Lipschitz continuous on $\R$ and ${\mathfrak f}_\rho,A_\rho$
converge to ${\mathfrak f},A$,  respectively, uniformly
on compacts sets as $\rho\to 0$.

Using  classical techniques (cf.~Alt, Luckhaus \cite{AltLuckhaus} and
Lions \cite{Lions:Book69_Fr}), we can show
that there exists a weak solution $u_\rho\in L^p(0,T;W^{1,p}_0(\Om))$ to
problem $\eqref{eq:prob1}_\rho$ in the following sense:
\begin{equation}\label{eq:WeakFormRho}
    \begin{cases}
        \pt u_\rho+\div {\mathfrak f}_\rho(u_\rho)=
        \div\mathfrak{a}(\grad A_\rho(u_\rho))+{\EuScript S}\\
        \text{\qquad in \; $\!L^{p'}(0,T;W^{-1,p'}(\Om))\!+\!L^1(Q)$}, \qquad u_\rho|_{t=0}=u_0.
    \end{cases}
\end{equation}
Moreover, since ${\mathfrak f}_\rho\circ A_\rho^{-1}$ is Lipschitz
continuous, the $L^1$ contraction property and comparison
principle for weak solutions can be verified. It can be obtained
either by the technique of Otto \cite{Otto:L1_Contr} (doubling the
time variable) or using the theory of integral solutions and
nonlinear semigroup methods, consult for example 
\cite{CarrilloWittbold:99}. Besides, $u_\rho$ verifies the entropy
formulation of Definition \ref{def:entropy2} with fluxes
${\mathfrak f}_\rho,A_\rho$, where $\eta_c^\pm$ can be replaced by
regular ``boundary'' entropies $\eta_{c,\eps}^\pm$, whenever we
prefer to do so.

\textbf{(ii)} We claim that the following quantities are uniformly bounded in $\rho$:\\
$\bullet$ $\norm{u_\rho}_{L^\infty(\Om)}$ and $\norm{A_\rho(u_\rho)}_{L^p(0,T;W^{1,p}_0(\Om))}$;\\
$\bullet$ space translates of $A_\rho(u_\rho)$ in $L^1(Q)$ (consequence of previous estimate);\\
$\bullet$ time translates of $A_\rho(u_\rho)$ in $L^1(Q)$.

Indeed, for the first point consider the function
$$
M(t)=\norm{u_0}_{L^\infty(\Om)} +\int_0^t \norm{{\EuScript
S}(\tau,\cdot)}_{L^\infty(\Om)}\,d\tau,
$$
which is a solution of $(\ref{eq:prob1})_\rho$ with $x$-constant
data $\norm{u_0}_{L^\infty(\Om)}$, $ \norm{{\EuScript
S}(t,\cdot)}_{L^\infty(\Om)}$. The comparison principle mentioned
in (i) ensures that a.e.~on $Q$,
$$
-M(T)\leq -M(t)\leq u_\rho(t,x) \leq M(t)\leq M(T).
$$
Next, we employ $A_\rho(u_\rho)$ as a test function in
\eqref{eq:WeakFormRho}. The product between $\pt u_\rho$ and
$A_\rho(u_\rho)$ is handled using the usual chain rule argument
(see, e.g., \cite{AltLuckhaus,Otto:L1_Contr,CarrilloWittbold:99}),
where the relevant duality is between the space
$E:=L^{p}(0,T;W^{1,p}_0(\Om))\cap L^\infty(Q)$ and the space
$L^{p'}(0,T;W^{-1,p'}(\Om))+L^1(Q)\subset E^*$. Here we are also
exploiting the $L^\infty$ bound on ${\mathfrak f}_\rho(u_\rho)$ in
a straightforward fashion to treat the term ${\mathfrak
f}_\rho(u_\rho)\cdot\grad A_\rho(u_\rho)$; but notice that using
the Green-Gauss trick \eqref{Green-Gauss} below, we can supply a
finer analysis of this term.

For the third bullet point, we first use \eqref{eq:WeakFormRho} to get,
for a.e. $t,t+\Del\in (0,T)$,
$$
\int_\Om (u_\rho(t+\Del)-u_\rho(t))\,\xi = \int_t^{t+\Del}\!\!
\int_\Om \left[\;\Bigl(-{\mathfrak f}_\rho(u_\rho) +{\mathfrak
a}(\grad A_\rho(u_\rho))\Bigl)  \cdot \grad \xi\; +\;{\EuScript
S}\,\xi\;\right]
$$
for all $\xi\in W^{1,p}_0(\Om)\cap L^\infty(\Om)$. Taking
$\xi=A_\rho(u_\rho(t+\Del))-A_\rho(u_\rho(t))$ and integrating in
$t$, using the two previously obtained estimates, we deduce that
\begin{equation}\label{Transl-base}
    \iint_Q \abs{u_\rho(t+\Del)-u_\rho(t)}
    \abs{A_\rho(u_\rho(t+\Del))-A_\rho(u_\rho(t))}
    \leq \Const{} \abs{\Del}.
\end{equation}
Now, let $\pi$ be a (common for all $\rho$) concave modulus of
continuity for $A_\rho$ on $[-M(T),M(T)]$, $\Pi$ be its inverse,
and set $\tilde\Pi(r)=r\,\Pi(r)$. Let $\tilde\pi$ be the inverse of
$\tilde\Pi$. Note that $\tilde \pi$ is concave, continuous, and
$\tilde\pi(0)=0$. Set $v(t,x)=u_\rho(t+\Del,x)$ and $y(t,x)=u_\rho(t,x)$.
We have
\begin{equation*}
    \begin{split}
        \int_Q \abs{A_\rho(v)-A_\rho(y)}
        & =\int_Q\tilde\pi\Biggl(\tilde\Pi(\abs{A_\rho(v)-A_\rho(y)}) \Biggr)
        \\ & \leq |Q|\,\tilde\pi\Biggl(\frac{1}{|Q|}
        \int_Q\tilde\Pi( \abs{A_\rho(v)-A_\rho(y)})\Biggr).
    \end{split}
\end{equation*}
Since $|A_\rho(v)-A_\rho(y)|\leq \pi(|v-y|)$, we have
$\Pi(|A_\rho(v)-A_\rho(y)|)\leq |v-y|$ and
\begin{equation*}
    \begin{split}
        \tilde\Pi(|A_\rho(v)-A_\rho(y)|)
        & =\Pi(|A_\rho(v)-A_\rho(y)|)|A_\rho(v)-A_\rho(y)|
        \\ & \leq |v-y|\,|A_\rho(v)-A_\rho(y)|.
    \end{split}
\end{equation*}
Therefore, estimate \eqref{Transl-base} implies
\begin{equation}\label{TranslEst-truc3}
    \begin{split}
        &\int_Q \abs{A_\rho(u_\rho(t\!+\!\Del,x))-A_\rho(u_\rho(t,x))}
        \\ & \qquad \leq |Q|\,\tilde\pi\Biggl(\frac{1}{|Q|}\int_Q|v-y|\abs{A_\rho(v)-A_\rho(y)}\Biggr)
        \\ & \qquad
        = |Q|\,\tilde\pi \Biggl(\frac{1}{|Q|} J(\Del) \Biggr)
        \leq C\tilde\pi(C\Del)=:\omega_A(\Del),
    \end{split}
\end{equation}
where $\omega_A\in C(\R^+,\R^+)$, $\omega_A(0)=0$.

\textbf{(iii)} Thanks to the estimates in \textbf{(ii)} and
standard compactness results, there
exists a (not labelled) sequence $\rho \to 0$ such that\\
$\bullet$ $w_\rho=A_\rho(u_\rho)$  converges strongly in
$L^1(Q)$ and pointwise a.e.~on $Q$;\\
$\bullet$  $\grad w_\rho$ converges weakly in $L^p(Q)$;\\
$\bullet$ ${\mathfrak a}(\grad w_\rho)$ converges weakly in $L^{p'}(Q)$ to some limit $\chi$;\\
$\bullet$ $u_\rho$ converges to $\mu:Q\times(0,1)\in \R$ in 
the sense of \eqref{nonlinweakstar}.

Let us introduce the function
\begin{equation}\label{eq:weak-star-limit-u}
	u(t,x)=\int_0^1 \mu(t,x,\al)\dal,
	\qquad \text{for a.e.~$(t,x)\in Q$}.
\end{equation}
Thanks to the convergence of $A_\rho$ to $A$, we can identify the limit of
$w_\rho(\cdot,\cdot)$ with $\int_0^1A(\mu(\cdot,\cdot,\al))\,d\al$. Moreover, since
$w_\rho$ is converging strongly, $A(\mu(\cdot,\cdot,\al))$ is actually
independent of $\alpha\in (0,1)$ and equals $A(u(\cdot,\cdot))$.
Using distributional derivatives, we also identify the
limit of $\grad w_\rho$ with $\grad A(u)$.

\textbf{(iv)} We have now come to the main step of the proof,
namely to improve the weak convergence of $\grad A_\rho(u_\rho)$
to strong convergence, and to identify the weak limit of
${\mathfrak a}(\grad A_\rho(u_\rho))$ with ${\mathfrak a}(\grad
A(u))$, where  $u$ is defined in \eqref{eq:weak-star-limit-u}; of
course, the chief difficulty comes from the lack of strong
convergence of $u_\rho$.

We begin by specifying the test
function in \eqref{eq:WeakFormRho} as $w_\rho\, \zeta$, yielding
\begin{equation}\label{weakform-of-approxeqn-tmp}
    \underbrace{\int_0^T \langle\pt u_\rho,w_\rho\, \zeta\rangle}_{I_{1,\rho}}
    -\underbrace{\int_Q \!{\mathfrak f}_\rho(u_\rho)\cdot\grad w_\rho\, \zeta}_{I_{2,\rho}}
    + \int_Q {\mathfrak a}(\grad w_\rho) \cdot \grad w_\rho\, \zeta
    - \underbrace{\int_Q {\EuScript S} w_\rho\, \zeta}_{I_{3,\rho}}=0,
\end{equation}
where $w_\rho=A_\rho(u_\rho)$ and $\zeta\in \mathcal D([0,T))$ is
nonincreasing with $\zeta(0)=1$. Next, we pass to the limit into
the weak formulation \eqref{eq:WeakFormRho}, obtaining
\begin{equation}\label{eq:LimitEqTheor}
    \begin{cases}
        \pt u+\div \int_0^1  {\mathfrak f}(\mu)\,d\al= \div \chi\; +\; {\EuScript S}\\
        \text{\qquad in \; $\!L^{p'}(0,T;W^{-1,p'}(\Om))\!+\!L^1(Q)$}, \qquad u_\rho|_{t=0}=u_0.
    \end{cases}
\end{equation}
In \eqref{eq:LimitEqTheor}, we take $w\,\zeta$ as test function, where $w=A(u)$, $u$ is
defined in \eqref{eq:weak-star-limit-u}, and $\zeta$ is as specified above. The result is
\begin{equation}\label{weakform-of-limiteqn-tmp}
    \underbrace{\int_0^T \langle\pt u,w\, \zeta\rangle}_{I_{1}}
    - \underbrace{\int_Q \int_0^1 {\mathfrak f}(\mu)\cdot\grad w\, \zeta}_{I_{2}}
    + \int_Q \chi\cdot \grad A(u)\, \zeta
    - \underbrace{\int_Q {\EuScript S} w\, \zeta}_{I_{3}}=0.
\end{equation}

In order to later use the Minty-Browder trick, we shall 
combine \eqref{weakform-of-limiteqn-tmp} and
the ``$\rho\to 0$" limit of \eqref{weakform-of-approxeqn-tmp} to conclude
the validity of the following inequality:
\begin{equation}\label{eq:MintyIneqTheor}
    \int_Q \chi\cdot \grad A(u)
    \geq \liminf_{\rho\to 0} \int_Q
    {\mathfrak a}(\grad w_\rho) \cdot \grad w_\rho.
\end{equation}

A crucial role is played by the following calculation, which
reveals that the lack of strong convergence of ${\mathfrak f}_\rho(u_\rho)$
is not an obstacle. Indeed, a componentwise application
of Lemma \ref{FunctionalDepLemma} (i) yields the existence
of a Lipschitz continuous vector-valued
function $\widetilde A_{\mathfrak{f}}$ such that
\begin{equation}\label{eq:Af-KHK}
    \int_0^{z} {\mathfrak f}(s)\,dA(s)
    =\widetilde A_{\mathfrak{f}}(A(z)).
\end{equation}
Hence, by the chain rule and the Green-Gauss formula, we can calculate as follows:
\begin{align*}
	\int_Q \int_0^1 {\mathfrak f}(\mu)\cdot\grad A(u)
	& =\int_Q \int_0^1 {\mathfrak f}(\mu)\cdot \grad A(\mu)
	= \int_0^1 \int_0^T\int_\Om
	\div \widetilde A_{\mathfrak{f}}(A(\mu))
	\\ & =\int_0^T\int_{\pOm} \widetilde A_{\mathfrak{f}}(A(u))\cdot n=0,
\end{align*}
because for a.e.~$\alpha\in (0,1)$,
$$
A(\mu(\cdot,\cdot,\al))= A(u(\cdot,\cdot))\in L^p(0,T;W^{1,p}_0(\Om)).
$$

By similar (simpler) arguments and
$u_\rho\in L^p(0,T;W^{1,p}_0(\Om))$, we also have
\begin{equation}\label{Green-Gauss}
    \begin{split}
        \int_Q {\mathfrak f}_\rho(u_\rho)\cdot\grad A_\rho(u_\rho)
        &=\int_0^T\int_\Om \div
        \left(\int_0^{u_\rho} {\mathfrak f}_\rho(s)\,dA_\rho(s) \right)=0.
    \end{split}
\end{equation}
Consequently, we can make $I_2$
and $I_{2,\rho}$ (for each $\rho>0$) vanish.

Next, let us prove that $I_1\le  \lim_{\rho\to 0} I_{1,\rho}$.
As above, the duality products $\langle \pt u_\rho,A_\rho(u_\rho)\rangle$,
$\langle \pt u,A(u)\rangle$ are treated via the chain rule argument (cf.~\cite{AltLuckhaus}).
Set $B(z)=\int_0^z A(s)\ds$, $B_\rho(z)=\int_0^z A_\rho(s)\ds$, and note that
these functions are convex. Also, $B_\rho\to B$ uniformly on compact subsets of $\R$.
With the help of Jensen's inequality,
\begin{align*}
	I_1 & =\int_0^T \langle \pt u,A(u)\,\zeta\rangle
	=-\int_Q B(u) \zeta' - \int_\Om B(u_0)
	\\ & = \int_Q B\left(\int_0^1 \mu(t,x,\al)\dal\right) (-\zeta') -\int_\Om B(u_0)
	\\ & \leq \int_Q \int_0^1 B( \mu(t,x,\al))\dal\, (-\zeta') -\int_\Om  B(u_0)
	\\ & = \lim_{\rho\to 0} \left( -\int_Q B_\rho( u_\rho)\,\zeta' -\int_\Om B_\rho(u_0) \right)
	\\ & = \lim_{\rho\to 0} \int_0^T \langle \pt u_\rho,A_\rho(u_\rho)\,\zeta\rangle
	=  \lim_{\rho\to 0} I_{1,\rho}.
\end{align*}

Finally, it is clear that $I_{3,\rho} \to I_{3}$ as $\rho\to 0$.
Letting $\zeta$ tend to $\char_{[0,T)}$, the desired inequality \eqref{eq:MintyIneqTheor}
follows from subtracting the ``$\rho\to 0$" limit of \eqref{weakform-of-approxeqn-tmp} from
\eqref{weakform-of-limiteqn-tmp} and the above calculations.

Starting off from \eqref{eq:MintyIneqTheor}, we can use the
Minty-Browder trick (see, for example,
\cite{Lions:Book69_Fr,Brow,AltLuckhaus,BocMuratPuel} and the proof
of Theorem~\ref{ThFVConverge} in Section \ref{sec:convergence}) to deduce that
\begin{equation}\label{eq:MintyConvergence}
    {\mathfrak a}(\grad w_\rho)-{\mathfrak a}(\grad A(u))
    \to 0 \quad \text{weakly in $L^{p'}(Q)$ as $\rho\to 0$}.
\end{equation}
Thus $\chi={\mathfrak a}(\grad A(u))$.
Simultaneously, from the strict monotonicity of ${\mathfrak a}(\cdot)$
we deduce that, firstly, the convergence in \eqref{eq:MintyConvergence}
also takes place a.e.~in $Q$; secondly, that \eqref{eq:MintyIneqTheor}
actually holds with an equality sign. Next, we consider the
functions $g_\rho:={\mathfrak a}(\grad w_\rho) \cdot\grad w_\rho\ge 0$ and
$g:={\mathfrak a}(\grad A(u))\cdot\grad A(u)\ge 0$, and observe that
$$
\text{$g_\rho\to g$ a.e.~in $Q$}, \qquad \int_Q g_\rho \to \int_Q
g \quad \text{as $\rho\to 0$}.
$$
Hence, we deduce that a subsequence of $(g_\rho)_\rho$ converges to $g$ strongly in $L^1(Q)$,
cf.~\cite{Brow}, \cite[Lemma 5]{BocMuratPuel}, \cite[Lemma 8.4]{Droniou}.
Due to the coercivity of $\mathfrak{a}(\cdot)$, $\Bigl(|\grad w_\rho|^p\Bigr)_\rho$
is equi-integrable, so the Vitali theorem yields the strong $L^p$ convergence
of $\grad w_\rho$, along a subsequence if necessary, to a limit already
identified as $\grad w$, $w=A(u)$.

\textbf{(v)} By \eqref{eq:LimitEqTheor}, we readily conclude that
$(\mu,\mu,w)$ verifies \ref{def:entropyWeakp}. Now we can pass to
the limit in the entropy inequalities corresponding to
(\ref{eq:prob1})$_\rho$ and deduce \ref{def:entropy2p}.

Let us first show that $\grad \widetilde A^\rho_{(\eta_{c,\eps}^\pm)'}(w_\rho)$
converges weakly to $\grad \widetilde A_{(\eta_{c,\eps}^\pm)'}(w)$
in $L^p(Q)$. By Lemma \ref{FunctionalDepLemma} (i), $A^\rho_{(\eta_{c,\eps}^\pm)'}(\cdot)$
are uniformly Lipshitz continuous functions. Thus $\grad \widetilde A^\rho_{(\eta_{c,\eps}^\pm)'}(w_\rho)$
are uniformly bounded and weakly compact in $L^p(Q)$.
Moreover, $\widetilde A^\rho_{(\eta_{c,\eps}^\pm)'}(w_\rho)$
converges to $\widetilde A_{(\eta_{c,\eps}^\pm)'}(w)$
by Lemma \ref{FunctionalDepLemma} (ii) and because of the
pointwise convergence of $w_\rho$ to $w$. Using
the distributional convergence, we eventually work out our claim.

Now note that if $p> 2$, then $k$ is continuous. By the last
result of \textbf{(iv)},  we can assume without loss of generality
that $k(\grad w_\rho)$ converges to $k(\grad w)$ a.e.~in $Q$.
Moreover, $(k(\grad w_\rho))$ is bounded in $L^\frac{p}{p-2}(Q)$,
since $(\grad w_\rho)$ is bounded in $L^p(Q)$. Applying the Egorov
theorem and H\"older's inequality with exponents $p',p$ in the
product $\Bigl(k(\grad w_\rho)\grad \psi\Bigr)\cdot\grad
\widetilde A_{(\eta_{c,\eps}^\pm)'}(w_\rho)$, we deduce that
\begin{equation}\label{PassageInTildeA}
    \lim_{\rho\to 0} \int_Q k(\grad w_\rho)
    \grad \widetilde A_{(\eta_{c,\eps}^\pm)'}(w_\rho)\cdot\grad \psi
    = \int_Q k(\grad w)\grad \widetilde A_{(\eta_{c,\eps}^\pm)'}(w)
    \cdot \grad \psi.
\end{equation}
If $p\leq 2$, we fix a small $\delta>0$ and truncate $k(\cdot)$ in
the $\delta$-neigbourhood of the origin (if $k(\cdot)$ is replaced
by $k_\delta(\cdot)=\min\{k(\cdot),\min_{|\xi|\leq \delta}
k(\xi)\}$, the argument used for $p>2$ applies), and we analyze
separately the set $\{(t,x)\,|\,|\grad w^\rho(t,x)|<\delta \}$. On
this set,
$$
k(\grad w_\rho)\, \Bigl|\grad \widetilde
A_{(\eta_{c,\eps}^\pm)'}(w_\rho)\Bigr| \leq C\,
\|(\eta^\pm_{c,\eps})'\|_\infty |\grad w_\rho|^{p-1}\leq \Const
\,\delta^{p-1},
$$
uniformly in $\rho$. To conclude that \eqref{PassageInTildeA} still holds, we first pass
to the limit as $\rho\to 0$ for a fixed $\delta>0$, and then send $\delta\to 0$.

Let us take regular ``boundary'' entropy pairs $\eta_{c,\eps}^\pm$
such that $(\eta_{c,\eps}^\pm)'$ approximate $(\eta_c^\pm)'$
(extended by zero at the point $c$, by our convention), pointwise
a.e.~in $\R$ as $\eps\to 0$. We use \eqref{PassageInTildeA} to
pass to the limit in the entropy inequality corresponding to
(\ref{eq:prob1})$_\rho$. We pass to the limit in the remaining
terms in this entropy inequality using the continuity of
$\eta^\pm_{c,\eps}, {\mathfrak
q}^\pm_{c,\eps},(\eta^\pm_{c,\eps})'$ and the nonlinear $L^\infty$
weak-$\star$ convergence property \eqref{nonlinweakstar}. Finally,
we pass to the limit as $\eps\to 0$, rewriting $\grad \widetilde
A_{(\eta_{c,\eps}^\pm)'}(w)$ as $(\eta^\pm_{c,\eps})'(u) \grad w$
(consult Remark \ref{RemTwoFormsOfDiffusion}) and using the
Lebesgue dominated convergence theorem and the pointwise
convergences of $\eta^\pm_{c,\eps}, {\mathfrak q}^\pm_{c,\eps},
(\eta^\pm_{c,\eps})'$. The passage to the limit in the weak
formulation is similar.

\textbf{(vi)} We  conclude that $(\mu,\mu^*,A(u))$ is an entropy
double-process solution of \eqref{eq:prob1} such that $\mu^*=\mu$.
\end{proof}

Given Theorem \ref{th:EPS-ESexists},  the uniqueness of an
entropy double-process solution can be established using
Kruzhkov's method, along the lines of Carrillo \cite{Carrillo}.

\begin{theo}\label{th:EPSunique}
Suppose the assumptions stated in Section \ref{sec:intro} hold.
Let $(\mu,\mu^*,w)$ be an entropy double-process solution
of the initial-boundary value problem \eqref{eq:prob1}.
Then it is unique. Moreover, there exists a function
$u \in L^\infty(Q)$ such that
$$
\text{$\mu(t,x,\al)=u(t,x)=\mu^*(t,x,\al)$ for a.e.~$(t,x,\alpha)\in Q\times(0,1)$.}
$$
\end{theo}

We refer to Appendix A for a sketch of the proof.

Theorems \ref{th:EPS-ESexists} and \ref{th:EPSunique} as well as
the arguments of Appendix A imply

\begin{cor}[well-posedness]\label{CorExistUniqEntropySol}
Under the assumptions stated in Section \ref{sec:intro}, there
exists a unique entropy solution of the initial-boundary value
problem \eqref{eq:prob1}. Let $u$ and $v$ be two entropy solutions
of \eqref{eq:prob1} with initial data $u|_{t=0} =u_0 \in
L^{\infty}(\Om)$ and $v|_{t=0} = v_0 \in L^{\infty}(\Om)$ and
source terms $\EuScript S$ and $\EuScript T$ of the kind
\eqref{eq:source-terms-space}, respectively. 
For a.e.~$t \in (0,T)$, we have
$$
\int_{\Om} \left(u(t,x)-v(t,x)  \right)^+ \dx 
\leq \int_{\Om}\left(u_0-v_0 \right)^+ \dx
+ \int_0^t\!\!\int_\Om (\EuScript S-\EuScript T)^+.
$$
Consequently, if $u_0 \leq v_0$ a.e.~in $\Om$ and $\EuScript S
\leq \EuScript T$ a.e. on $Q$, then $u \leq v$ a.e.~in $Q$.
Finally, if $u_0=v_0$ a.e.~in $\Om$ and $\EuScript S = \EuScript
T$ a.e. on $Q$, then $u=v$ a.e.~in $Q$.
\end{cor}

The upcoming sections are concerned with the construction of
finite volume schemes for which the corresponding
discrete solutions converge to the unique entropy solution of \eqref{eq:prob1}
as the discretization parameter (mesh size) tends to zero.
The convergence proof will attempt to mimic the
proof of Theorem \ref{th:EPS-ESexists}.

\section{Discrete duality finite volume (DDFV) schemes}\label{sec:FV}

Let $\Om$ be a polygonal (respectively, polyhedral) open bounded
subset of $\R^d$, $d=2$ (respectively, $d=3$). In what follows, we
introduce most of the notation related to DDFV
schemes; each piece of new notation is given in italic script.

\subsection{Construction of ``double'' conformal meshes}~

 $\bullet$ A {\it partition} of $\Om$ is a finite set of disjoint
open polygonal (respectively, polyhedral) subsets of $\Om$ such
that $\Om$ is contained in their union, up to a set of zero
$d$-dimensional measure.

Following Hermeline \cite{Herm}, Domelevo, Omn\`es \cite{DomOmnes}
and Andreianov, Boyer, Hubert \cite{ABH``double''}, we consider a
DDFV mesh which is a triple
$\Tau=\Bigl(\overline{\Mrond},\overline{\dMrond},\Srond \Bigr)$
described below.

$\bullet$ We let $\Mrond$ be a partition of $\Om$ into triangles
(respectively, tetrahedra); a more general case is discussed in
Section~\ref{sec:Extensions}\footnote{In particular, in the two
dimensional case we can partition $\Omega$ into polygons that
admit a circumscribed circle. In the three dimensional case , we
can partition $\Omega$ in polyhedra that have triangular faces and
admit a circumscribed ball.}.
We assume that the mesh satisfies the
Delaunay condition (see, e.g., \cite{EyGaHe:book});
for simplicity of the representation, the reader may assume that
each triangle (respectively, tetrahedron) contains the centre if
its circumscribed circle (respectively, ball). We assume in addition
\begin{equation}\label{eq:AssumptMesh}
    \left|
    \begin{array}{l}
        \text{if $d=3$, each face of each tetrahedron of $\Mrond$}\\
        \text{contains the centre of its circumscribed circle.}
    \end{array}\right.
\end{equation}
Although the definition of the scheme does not require condition
\eqref{eq:AssumptMesh} (see Remark~\ref{rem:Delonnay} below), we
do need this condition in order to deduce the discrete entropy
inequalities and to prove that the scheme converges.

Each {\it control volume} $\K\in\Mrond$ is supplied with a {\it centre}
$\xK$ that we choose to be the centre of the circle (respectively, ball)
circumscribed around $\K$. We call $\ptl\Mrond$ the set of all
edges (respectively, faces) of control volumes that are included in
$\ptl\Om$. These edges (respectively, faces) are considered as {\it
boundary control volumes}; for $\K\in\ptl\Mrond$, we choose the
middle of $\K$ (respectively, the centre of the circle circumscribed
around $\K$) for the centre $\xK$. We denote by
$\overline{\Mrond}$ the union $\Mrond\cup\ptl\Mrond$. We call {\it
vertex} (of $\Mrond$) any vertex of any control volume
$\K\in\overline{\Mrond}$.

\begin{figure}[bth]
\begin{center}
\scalebox{0.9}{\input{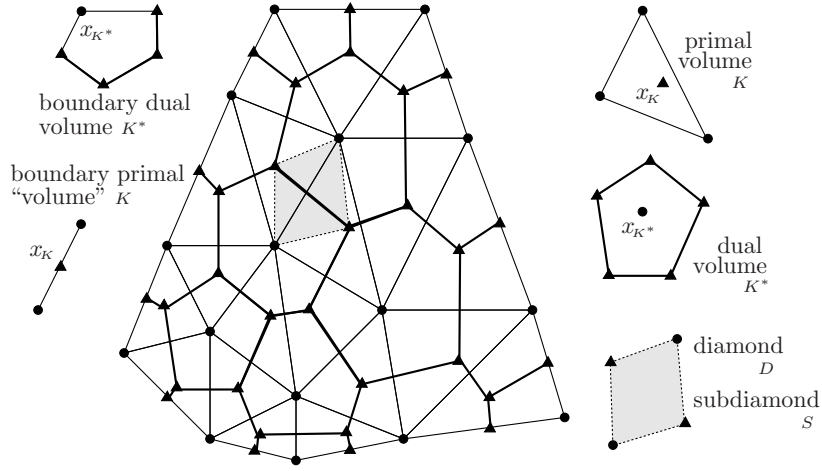}}
\end{center}
\caption{2D primal and dual meshes; diamond (subdiamond).}
\label{Fig-2D}
\end{figure}

$\bullet$ (see Figure~\ref{Fig-2D}) We take $\overline{\dMrond}$
as the partition of $\Om$ into {\it dual control volumes} $\dK$,
supplied with {\it dual centres} $\xdK$, such that $\xdK$ is a
vertex of $\Mrond$ and $\dK$ is the subset of points of $\Om$ that
are closer\footnote{in order to avoid pathological situations
which could appear in non-convex domains, e.g., in domains with
cracks, here the distance between two points $x,y$ of
$\overline{\Om}$ is understood as the length of the shortest path
which connects $x$ with $y$ and which lies {\it within}
$\overline{\Om}$ } to $\xdK$ than to any other vertex of $\Mrond$.
In other words, $\overline{\dMrond}$ is the Vorono\"{\i} mesh
constructed from the vertices of $\Mrond$. If $\xdK\in\Om$, we say
that $\dK$ is a {\it dual control volume} and write
$\dK\in\dMrond$; and if $\xdK\in\ptl\Om$, we say that $\dK$ is a
{\it boundary dual control volume} and write $\dK\in\ptl\dMrond$.
Thus $\overline{\dMrond}=\dMrond\cup \ptl\dMrond$. We call {\it
dual vertex} (of $\dMrond$) any vertex of any dual control volume
$\dK\in\overline{\dMrond}$. Note that by the choice of $\xK$, the
set of centres coincides with the set of dual vertices, and the
set of vertices coincides with the set of dual centres. In other
words, $\overline{\Mrond}$ and $\overline{\dMrond}$ are finite
volume meshes that are dual each one to the other.

$\bullet$ We call {\it neighbours} of $\K$, all control volumes
$\L\in\overline{\Mrond}$  such that $\K$ and $\L$ have a common
edge (respectively, common face). The set of all neighbours of $\K$ is
denoted by $\Nrond(\K)$. Note that if $\L\in \Nrond(\K)$, then
$\K\in\Nrond(\L)$; in this case we simply say that $\K$ and $\L$
are (a couple of) neighbours.

\begin{figure}[b]
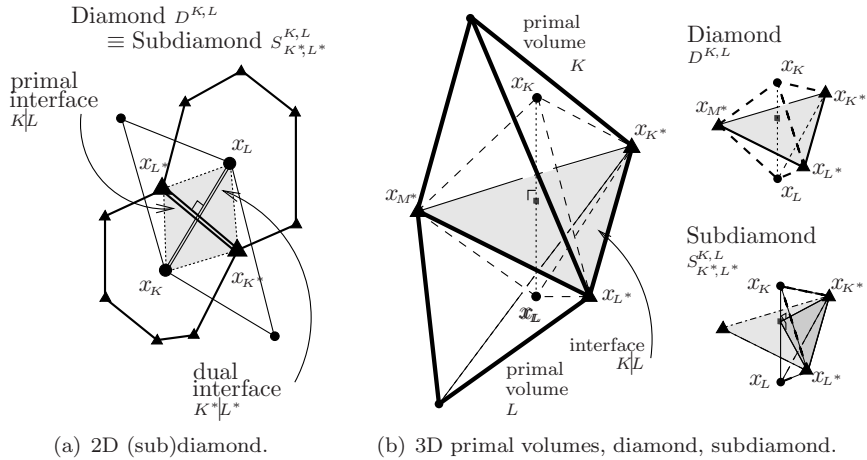

\centering
\subfigure[2D (sub)diamond.] 
{ \label{Fig-2D-diamond}
 \scalebox{0.9}
 {\input{ABK-2D-subdiamond.pstex_t}}
} \hspace{5mm}
\subfigure[3D primal volumes, diamond, subdiamond.] 
{
    \label{Fig-3D}
 \scalebox{0.9}
 {\input{ABK-3D-neigh-volumes.pstex_t}}}
\caption{Diamonds and subdiamonds.}
\label{fig:2D/3D-diamonds} 
\end{figure}

$\bullet$ (see Figures~\ref{Fig-2D} and \ref{Fig-3D}) If $\K$ and
$\L$ are neighbours, we denote by $\KIL$ the {\it interface}
$\ptl\K\cap\ptl\L$ between $\K$ and $\L$. The set of all
interfaces is denoted by $\Erond$.

$\bullet$ In the same way, we denote by $\dNrond(\dK)$ the set of
{\it (dual) neighbours} of a dual control volume $\dK$, and by
$\dKIdL$, the {\it (dual) interface} $\ptl\dK\cap\ptl\dL$ between
dual neighbours $\dK$ and $\dL$. The set of all dual interfaces is
denoted by $\dErond$.

 $\bullet$ (see Figure~\ref{fig:2D/3D-diamonds}) The meshes $\overline{\Mrond}$ and
$\overline{\dMrond}$ induce partitions of $\Om$ into diamonds and subdiamonds. 
Let us describe them separetely for $d=2$ and $d=3$.

For $d=2$ (see Figure~\ref{Fig-2D-diamond}), if
$\K,\L\in\overline{\Mrond}$ are neighbours, then there exists a
unique couple of dual neighbours $\{\dK,\dL\}$ such that the
interface $\KIL$ is the segment with summits $\xdK$ and $\xdL$.
Then the quadrilateral $\DM^{\ptK,\ptL}_{\ptdK\!\!,\ptdL}$ which
is either the union (if $\xK,\xL$ lie on different sides from
$\KIL$) or the difference (if $\xK,\xL$ lie on the same side from
$\KIL$) of the triangles $\xK\xdK\xdL$, $\xL\xdK\xdL$ is called a
{\it diamond}; it is also unambiguously denoted by $\DM^{\ptK,\ptL}$.

For $d=2$, every diamond is also called a {\it subdiamond}; the subdiamond which coincides with a diamond
$\DM^{\ptK,\ptL}$ is denoted by $\SDM^{\ptK,\ptL}_{\ptdK\!\!,\ptdL}$.

For $d=3$ (see Figure~\ref{Fig-3D}), if $\K,\L\in\overline{\Mrond}$ 
are neighbours, then there exists a
unique triple of dual neighbours $\{\dK,\dL,\dM\}$ (which are
neighbours pairwise) such that the interface $\KIL$ is the
triangle with summits $\xdK$, $\xdL$ and $\xdM$. Then the
polyhedron $\DM^{\ptK,\ptL}_{\ptdK\!\!,\ptdL\!\!,\ptdM}$ which is
either the union (if $\xK,\xL$ lie on different sides from $\KIL$)
or the difference (if $\xK,\xL$ lie on the same side from $\KIL$)
of the pyramids $\xK\xdK\xdL\xdM$, $\xL\xdK\xdL\xdM$ is called a
{\it diamond}; it is also unambiguously denoted by
$\DM^{\ptK,\ptL}$. Each diamond is split into three subdiamonds;
e.g., the {\it subdiamond}  $\SDM^{\ptK,\ptL}_{\ptdK\!\!,\ptdL}$
is the convex hull of $\xK,\xdK,\xL,\xdL$.

We denote by  $\Drond$,$\Srond$ the sets of all
diamonds and the set of all subdiamonds, respectively. Generic
elements of $\Drond$,$\Srond$ are denoted by $\DM$,$\SDM$, respectively.

\begin{rem}\normalfont \label{rem:Delonnay} If we drop condition
\eqref{eq:AssumptMesh}, the orthogonal projection of $\xK$ (which
coincides with the projection of $\xL$) on $\KIL$ may not be
contained within $\KIL$. To cope with this problem, one could
consider subdiamonds of signed volume, not necessarily contained
within the corresponding diamonds. Up to a permutation of the
subscripts $\ptdK,\ptdL,\ptdM$, we have instead of the decomposition
$\DM^{\ptK,\ptL}_{\ptdK\!\!,\ptdL\!\!,\ptdM}=\SDM^{\ptK,\ptL}_{\ptdK\!\!,\ptdL}\cup
\SDM^{\ptK,\ptL}_{\ptdL\!\!,\ptdM} \cup
\SDM^{\ptK,\ptL}_{\ptdM\!\!,\ptdK}$, the decomposition
$\DM^{\ptK,\ptL}_{\ptdK\!\!,\ptdL\!\!,\ptdM}=\Bigl(\SDM^{\ptK,\ptL}_{\ptdK\!\!,\ptdL}\cup
\SDM^{\ptK,\ptL}_{\ptdL\!\!,\ptdM} \Bigr)\setminus
\SDM^{\ptK,\ptL}_{\ptdM\!\!,\ptdK}$; in this case the volume of
$\SDM^{\ptK,\ptL}_{\ptdM\!\!,\ptdK}$ will be taken with the sign
``minus''. Under this convention, Lemma~\ref{geomlem} below holds
true, so that formulas \eqref{DiscrGradD3},
\eqref{DiscrDiv}-\eqref{DiscrDivdM} below still yield consistent
discrete gradient and discrete divergence operators which enjoy
the discrete duality property \cite{AndrBend}. But the
discrete entropy dissipation inequalities of
Proposition~\ref{ParabEntropyDissipation} would fail, which
undermines the subsequent convergence analysis.
\end{rem}

$\bullet$ For all bounded  set $E\subset \R^d$, set $\diam(E)=\sup_{x,\hat x\in E} |x-\hat x|$.

$\bullet$ We denote by $m_E$ the measure of an object $E$ in its
natural dimension (i.e., the  $d$-dimensional measure, if $E$ is a
control volume, a dual control volume, a subdiamond or a diamond;
and the $(d \!-\! 1)$-dimensional measure, if $E$ is an interface
or a part of an interface). According to
Remark~\ref{rem:Delonnay}, for the definition of the scheme we
could drop \eqref{eq:AssumptMesh}, in which case for a subdiamond
$\SDM^{\ptK,\ptL}_{\ptdK\!,\ptdL}$ such that
$\SDM^{\ptK,\ptL}_{\ptdK\!,\ptdL}\cap \DM^{\ptK,\ptL}=\O$ its
volume is taken with the sign ``minus''.

\subsection{Mesh parameters and regularity of meshes}
~

$\bullet$ We define the {\it size} of the mesh by
$\size(\Tau)=\max\nolimits_{E\in \overline{\ptMrond}
\cup\overline{\ptdMrond} \cup \ptDrond} \diam(E)$.

$\bullet$  Following \cite{ABH``double''}, we call the maximum
among
$$
\max_{\ptdK} \;\text{card}(\dNrond(\dK)),\quad
\max_{\ptK} \frac{(\diam(\K))^d}{\mK}, \quad
\max_{\ptdK} \frac{(\diam(\dK))^d}{\mdK},
$$
$$
\max_{\ptK\cap \ptD\neq \O} \Biggl(\frac{\diam(\K)}{\diam(\DM)}+
\frac{\diam(\DM)}{\diam(\K)}\Biggr), \quad \max_{\ptdK\cap
\ptD\neq \O} \Biggl(\frac{\diam(\dK)}{\diam(\DM)}+
\frac{\diam(\DM)}{\diam(\dK)}\Biggr),
$$
(where the maximums are taken over all $\K\in\overline{\Mrond}$,
$\dK\in\overline{\dMrond}$, $\DM\in \Drond$) the {\it regularity constant} 
of the mesh and we denote it by $\reg(\Tau)$. 
Roughly speaking, this constant controls the ratio of dimensions
of neighbouring control volumes, diamonds and dual control
volumes, as well as the proportions of each volume.

In all the discrete estimates and convergence results stated below, we require the family of
meshes $(\Tau_h)_h$ to have regularity constants $\reg(\Tau_h)$ 
that are uniformly bounded in $h$. In the sequel, whenever there 
is a dependency of various constants on $\reg(\Tau)$, we tacitly assume
that this dependency is increasing.

\subsection{Discrete gradient and divergence operators}
Diamonds permit to define the discrete gradient operator, while subdiamonds
permit to define the discrete divergence operator
(see  \eqref{DiscrGrad}, \eqref{DiscrGradD2}, \eqref{DiscrGradD3} and
\eqref{DiscrDiv}, \eqref{DiscrDivM}, \eqref{DiscrDivdM} below, respectively).
Both are needed to discretize the second order ``diffusion'' operator in equation
\eqref{eq:prob1}. But first we need to introduce some more notation.

\begin{figure}[bth]
\begin{center}
\scalebox{0.9} {\input{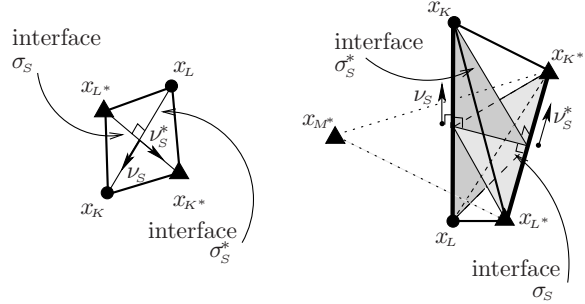}}
\end{center}
\caption{Notation in a subdiamond (2D and 3D).}
\label{Fig-in-a-subdiamond}
\end{figure}

$\bullet$ (see Figure~\ref{Fig-in-a-subdiamond}) For a subdiamond
$\SDM=\SDM^{\ptK,\ptL}_{\ptdK\!,\ptdL}$, we denote by
$\edge=\edgeS$, $\dedge=\dedgeS$ the (parts of the) interfaces
$\SDM\cap \KIL$ and $\SDM\cap \dKIdL$, respectively, and by
$\nuS$, $\nudS$, unit normal vectors to $\edgeS$ and $\dedgeS$,
respectively (their orientation is chosen arbitrarily).

$\bullet$ For a diamond $\DM=\DM^{\ptK,\ptL}$, we denote by
$\ProjD$, $\dProjD$ the operators of orthogonal projection of
$\R^d$ on the subspaces $<\overrightarrow{\xK\xL}>$ and on
$<\overrightarrow{\xK\xL}>^\bot$, respectively. One should note that we have
$<\nuS>=<\overrightarrow{\xK\xL}>$ and $<\nudS>\subset
<\overrightarrow{\xK\xL}>^\bot$ for all $\SDM \in \Srond$ such
that $\SDM\subset \DM$.

$\bullet$ For a couple of neighbours $\K,\L\in\overline{\Mrond}$, denote by $\dKL$, $d_{\ptK,\ptKIL}$, 
and $\nuKL$ the distance between $\xK$ and $\xL$ , the distance from $\xK$ to
$\KIL$, and the unit normal vector to $\K\I\L$ pointing from $\K$
to $\L$, respectively. More generally, if $\K\in {\Mrond}$, then
$\nuK$ denotes the exterior unit normal vector to $\ptl\K$. In the
same way, for neighbours $\dK,\dL\in\overline{\dMrond}$ we define
$\ddKdL$, $d_{\ptdK,\ptdKIdL}$, and $\nudKdL$; for
$\dK\in\overline{\dMrond}$, we define $\nudK$.

\begin{rem}\normalfont \label{rem:conformal}
Note that by construction both meshes
$\overline{\Mrond},\overline{\dMrond}$ are {\it conformal}
(orthogonal) in the sense if \cite{EyGaHe:book}); combined with
the Delaunay condition, this means that
$\nuKL\cdot\overrightarrow{\xK\xL}=\dKL$, $\nudKdL\cdot
\overrightarrow{\xdK\xdL}=\ddKdL$ for all neighbours $\K,\L$ and
$\dK,\dL$, respectively.

The conformity property is particularly important for our $L^1$
framework imposed by the possible degeneracy of the diffusion term
and the presence of the hyperbolic convective term. On the other
hand, if this term is dropped, non-conformal double meshes can be
considered for $d=2$ (see \cite{Herm,DomOmnes,ABH``double''}) and
$d=3$ (see
\cite{Pierre,CoudierePierre,Herm-3D,ABK-FVCA5,AndrBend,CoudiereHubertPierre})
within the variational framework.
\end{rem}

 $\bullet$
A {\it discrete function on ${\Om}$} is a set
$w^\ptTau=\Bigl(\uMrond,\udMrond\Bigr)$ consisting of two sets of
real values $w^\ptMrond=(w_\ptK)_{\ptK\in\ptMrond}$ and
$w^\ptdMrond=(w_\ptdK)_{\ptdK\in\ptdMrond}$. The set of all such
functions is denoted by $\R^\ptTau$.

A {\it discrete function on $\overline{\Om}$} is a set
$$w^\ptbarTau=\Bigl(w^\ptMrond,w^\ptdMrond,w^{\ptptl\ptMrond},w^{\ptptl\ptdMrond}
\Bigr)\equiv\Bigl(w^\ptTau,w^{\ptptl\ptMrond},w^{\ptptl\ptdMrond}
\Bigr)$$ consisting of four sets of real values
$$w^\ptMrond\!\!=(w_\ptK)_{\ptK\in\ptMrond},
\; w^\ptdMrond\!\!\!=(w_\ptdK)_{\ptdK\in\ptdMrond}, \;
w^{\ptptl\ptMrond}\!\!\!=(w_\ptK)_{\ptK\in\ptptl \ptMrond}, \;
w^{\ptptl\ptdMrond}\!\!\!\!=(w_\ptdK)_{\ptdK\in\ptptl
\ptdMrond}.$$ The set of all such functions is denoted by
$\R^\ptbarTau$. In case all the components of $w^{\ptptl\ptMrond}$
and of $w^{\ptptl\ptdMrond}$ are zero, we write $w^\ptbarTau\in
\R^\ptbarTau_0$.

$\bullet$ A {\it discrete field on $\Om$} is a set
$\Frond^\ptTau=\Bigl( {\mathcal F}_\ptD  \Bigr)_{\ptD\in\ptDrond}$
of vectors of $\R^d$. The set of all such functions is denoted by $(\R^d)^\ptDrond$.

$\bullet$ On the set $\R^\ptbarTau$ of discrete functions
$w^\ptbarTau$ on $\overline{\Om}$, we define the {\it discrete gradient} operator
$\grad^\ptTau[\cdot]$ by
\begin{equation}\label{DiscrGrad}
    \grad^\ptTau: w^\ptbarTau\in \R^\ptbarTau
    \mapsto \grad^\ptTau
    w^\ptbarTau=\Bigl( \grad_\ptD w^\ptbarTau \Bigr)_{\ptD\in\ptDrond}
    \in (\R^d)^\ptDrond
\end{equation}
where $\grad^\ptTau w^\ptbarTau$ is the discrete field on $\Om$ with values

for $d=2$:
\begin{equation}\label{DiscrGradD2}
    \grad_\ptD w^\ptbarTau=\frac{\wL-\wK}{\dKL}\,\nuKL
    +\frac{\wdL-\wdK}{\ddKdL}\,\nudKdL \qquad
    \text{for $\DM=\DM^{\ptK,\ptL}=\SDM^{\ptK,\ptL}_{\ptdK\!\!,\ptdL}$};
\end{equation}

for $d=3$:
\begin{equation}\label{DiscrGradD3}
    \left| \begin{array}{l}
        \dsp \grad_\ptD w^\ptbarTau
        =\frac{\wL-\wK}{\dKL}\nuKL +\frac{2}{m_\ptD}
        \Biggl( m_{S^{\ptK,\ptL}_{\ptdK\!,\ptdL}}\frac{\wdL-\wdK}{\ddKdL}
        \,\nudKdL\\[8pt] \dsp\hspace*{15mm}
        + m_{S^{\ptK,\ptL}_{\ptdL\!,\ptdM}}\frac{\wdM-\wdL}{\ddKdL}\,\nudLdM
        + m_{S^{\ptK,\ptL}_{\ptdM\!,\ptdK}}\frac{\wdK-\wdM}{\ddMdK}\,\nudMdK \Biggr)
        \qquad\\[12pt]
        \dsp  \text{for $\DM=\DM^{\ptK,\ptL}_{\ptdK\!\!,\ptdL\!\!,\ptdM}
        =\SDM^{\ptK,\ptL}_{\ptdK\!\!,\ptdL}\cup
        \SDM^{\ptK,\ptL}_{\ptdL\!\!,\ptdM} \cup
        \SDM^{\ptK,\ptL}_{\ptdM\!\!,\ptdK} $}.
    \end{array}\right. \hspace*{-15mm}
\end{equation}

\begin{rem}\normalfont \label{ExplicDiscrGrad}
Formulas \eqref{DiscrGradD2} and
\eqref{DiscrGradD3} have the following common meaning.
The vector $\grad_\ptD w^\ptbarTau$ is the unique element of $\R^d$ such that
$\ProjD(\grad_\ptD w^\ptbarTau)=\frac{\wdL-\wdK}{\ddKdL}
\,\nudKdL$. Further, for $d=2$, $\dProjD(\grad_\ptD w^\ptbarTau)$
is the gradient of the (unique) affine function on the interface
$\KIL$ (which is a segment with summits $\xdK,\xdL$) that takes
the values $\wdK$,$\wdL$ at the points $\xdK$ and $\xdL$,
respectively. Similarly, for $d=3$, $\dProjD(\grad_\ptD
w^\ptbarTau)$ is the gradient of the (unique) affine function on
the interface $\KIL$ (which is a triangle with summits
$\xdK,\xdL,\xdM$) that takes the values $\wdK$,$\wdL$,$\wdM$ at
the points $\xdK$,$\xdL$,$\xdM$, respectively.

Thus, the primal mesh $\overline{\Mrond}$ serves to reconstruct
one component of the gradient, which is the one in the direction
$\overrightarrow{\xK\xL}$. The dual mesh $\overline{\dMrond}$
serves to reconstruct the $(d \!-\! 1)$ other components which are
the components in the $(d \!-\! 1)$-dimensional hyperplane
containing $\KIL$ and is orthogonal to
$\overrightarrow{\xK\xL}$.
\end{rem}

The first and second assertions of Remark~\ref{ExplicDiscrGrad}
are evident. Note that formula \eqref{DiscrGradD2} easily
generalizes to quite arbitrary non conformal double meshes (see
\cite[Lemma 2.4]{ABH``double''}). The third assertion is a direct
consequence of the 2D reconstruction result of
Lemma~\ref{lem:ReconstrPropTriangles} given and proved in Appendix
B (see also \cite{ABK-FVCA5,AndrBend}).

\begin{rem}\normalfont \label{rem:consistAffines}
The discrete gradient is exact on affine functions. 
More precisely, let $\DM$ be a diamond
($\DM=\DM^{\ptK,\ptL}_{\ptdK\!,\ptdL}$, if $d=2$;
$\DM=\DM^{\ptK,\ptL}_{\ptdK\!,\ptdL\!,\ptdM}$, if $d=3$). Let
$w(x):=w_0+r\cdot x$, $w_0,r\in \R^d$, be an affine function. 
If $w^\ptbarTau$ is a discrete function with values
$$
\begin{array}{l}
w_\ptK=w(x_\ptK),\, w_\ptL=w(x_\ptL);\\[3pt] w_\ptdK=w(x_\ptdK),\,
w_\ptdL=w(x_\ptdL) \; \text{(and $w_\ptdM=w(x_\ptdM)$ if $d=3$)},
\end{array}
$$
then $\grad_\ptD w^\ptbarTau=r\equiv \grad w$. 
This property follows by a straightforward comparison of the
formulas \eqref{DiscrGradD2} and \eqref{DiscrGradD3} for the discrete
gradient with the reconstruction formulas of the next lemma.
\end{rem}

\begin{lem}\label{geomlem}
Consider $\DM=\DM^{\ptK,\ptL}\in\Drond$.  With the notation above,
for all $r\in \R^d$ one has the following reconstruction
properties:\\[5pt]
for $d=2$, $\;\,r\;=\;(r\cdot\nuKL)\,\nuKL + (r\cdot\nudKdL)\,\nudKdL;$\\[15pt] 
for $d=3$,
\begin{equation*}
\begin{split}
r & = (r\cdot\nuKL)\,\nuKL
+\frac{2}{m_\ptD}\Biggl( m_{S^{\ptK,\ptL}_{\ptdK\!,\ptdL}}
(r\cdot\nudKdL)\,\nudKdL
+ m_{S^{\ptK,\ptL}_{\ptdL\!,\ptdM}}(r\cdot\nudLdM)\,\nudLdM
\\ & \qquad\qquad\qquad\qquad\qquad\qquad \qquad\qquad\qquad
+ m_{S^{\ptK,\ptL}_{\ptdM\!,\ptdK}}(r\cdot\nudMdK)\,\nudMdK
\Biggr).
\end{split}
\end{equation*}
\end{lem}

\begin{proof}
For $d=2$, the claim is a straightforward consequence of the
conformity of the meshes (see Remark~\ref{rem:conformal}); $\nuKL$, $\nudKdL$ form
an orthonormal basis of $\R^2$. When $d=3$,
the claim follows from the orthogonality of $\nuKL$ to $\KIL$ and
from the 2D reconstruction property of
Lemma~\ref{lem:ReconstrPropTriangles} (cf.~Appendix B) applied
in the plane containing $\KIL$.
\end{proof}

\begin{rem}\normalfont \label{rem:Gen4D}
The fourth assertion of Remark~\ref{ExplicDiscrGrad} indicates
possible generalizations to the multi-dimensional case.
Unfortunately, it can be shown that if $d\geq 4$, the direct
generalization of the reconstitution formula of
Lemma~\ref{geomlem} holds only for meshes $\overline{\Mrond}$ with
very special geometries, such as the uniform simplicial meshes
(see Remark~\ref{rem:ReconstrGeneralizationMultiD}, which has to
be combined with an induction argument on the dimension $d$ in
order to link the weighted projections on the edges appearing in
Lemma~\ref{geomlem} with the weighted projections on the faces
appearing in Lemma~\ref{lem:ReconstrPropTriangles}).
\end{rem}

$\bullet$ For $\SDM\in\Srond$  such that $\SDM\subset \DM$ with
$\DM\in\Drond$, we assign $\grad_\ptS \uTau=\grad_\ptD \uTau$.
More generally, if $\Frond^\ptTau$ is a discrete field on $\Om$,
we assign $\Frond_\ptS=\Frond_\ptD$ for $\SDM\subset \DM$.

For $\K\in\Mrond$, we denote by $\Vrond(\K)$ the set of all
subdiamonds $\SDM\in \Srond$ such that $\K\cap \SDM \neq \emptyset$. In
the same way, for $\dK\in\overline{\dMrond}$ we define the set
$\dVrond(\dK)$ of subdiamonds intersecting $\dK$.

$\bullet$ On the set $(\R^d)^\ptDrond$ of discrete fields
$\Frond^\ptTau$, we define the {\it discrete divergence}
operator $\div^\ptTau[\cdot]$ by
\begin{equation}\label{DiscrDiv}
\div^\ptTau:\Frond^\ptTau\in (\R^d)^\ptDrond \mapsto
v^\ptTau=\div^\ptTau [\Frond^\ptTau] \in \R^\ptTau,
\end{equation}
 where the discrete function
$v^\ptTau=\Bigl( v^\ptMrond,v^\ptdMrond \Bigr)$ on $\Om$ is given
by
\begin{equation}\label{DiscrDivM}
v^\ptMrond=(v_\ptK)_{\ptK\in\ptMrond} \; \text{~with} \quad
v_\ptK=\frac{1}{\mK}\sum\nolimits_{\ptS\in\ptVrond(\ptK)}\!\!\!\!
m_\edgeS \Frond_\ptS \cdot \nuK, \; \text{where} \quad
\nuK=\nuK|_\ptS;
\end{equation}
\begin{equation}\label{DiscrDivdM}
v^\ptdMrond=(v_\ptdK)_{\ptdK\in\ptdMrond}, \;  \quad
v_\ptdK=\frac{1}{\mdK}\sum\nolimits_{\ptS\in\ptdVrond(\ptdK)}\!\!\!\!
m_\dedgeS \Frond_\ptS \cdot \nudK,  \; \quad\nudK=\nudK|_\ptS.
\end{equation}
In \eqref{DiscrDivM},\eqref{DiscrDivdM} for $\SDM$ given,
$\nuK=\nuK|_\ptS$ denotes the restriction on $\sigma_\ptS$ of the
unit normal vector $\nuK$ to $\ptl K$ exterior to $\K$; therefore
it means the one of the vectors $\nuS,-\nuS$ that is exterior to
$\K$ (see Figure~\ref{Fig-in-a-subdiamond}). Similarly,
$\nudK=\nudK|_\ptS$ is the one of the vectors $\nudS,-\nudS$ that
is exterior to $\dK$.

In fact, formulas \eqref{DiscrDivM}, \eqref{DiscrDivdM} can be conveniently
expressed in terms of  vector products involving the discrete
field $\Frond_\ptS$ and specific geometric objects depicted in
Figure~\ref{Fig-in-a-subdiamond} (see \cite{ABK-FVCA5,AndrBend}).

\subsection{Penalization operator}\label{sec:Penalization}
On the set $\R^\ptbarTau$ of discrete functions $w^\ptbarTau$ on
$\overline{\Om}$, we define the operator $ \Prond^\ptTau [\cdot]$
of {\it double mesh penalization} by
\begin{equation*}
\Prond^\ptTau:w^\ptbarTau\in \R^\ptbarTau  \mapsto
v^\ptTau=\Prond^\ptTau [w^\ptbarTau] \in \R^\ptTau,
\end{equation*}
 where the discrete function
$v^\ptTau=\Bigl( v^\ptMrond,v^\ptdMrond \Bigr)$ on $\Om$ is given
by
\begin{equation}\label{DiscrPenalM}
v^\ptMrond=(v_\ptK)_{\ptK\in\ptMrond} \quad \text{with} \quad
v_\ptK=(d-1)\,\frac{1}{\size(\Tau)} \,\frac{1}{\mK}
\sum_{\ptdK\in \overline{\ptdMrond}}
m_{\ptK\cap\ptdK}(w_\ptK-w_\ptdK) ;
\end{equation}
\begin{equation}\label{DiscrPenaldM}
v^\ptdMrond=(v_\ptdK)_{\ptdK\in\ptdMrond} \quad \text{with} \quad
v_\ptdK=\frac{1}{\size(\Tau)} \,\frac{1}{\mdK}
\sum_{\ptK\in\overline{\ptMrond}}
m_{\ptK\cap\ptdK}(w_\ptdK-w_\ptK).
\end{equation}
The penalization is needed in order to
ensure (without using the strong convergence of $\grad^\ptTau
w^\ptbarTau$, cf.~the proof of \cite[Theorem 5.1]{ABH``double''}),
that the two components of a discrete ``double'' function
$w^\ptbarTau$ converge to the same limit.
\begin{rem}\normalfont \label{rem:penalChoice}
The choice of penalization operator we propose here is just the
simplest possibility. In \eqref{DiscrPenalM},\eqref{DiscrPenaldM},
the difference $(w_\ptdK-w_\ptK)$ could be replaced by
$|w_\ptdK-w_\ptK|^{p-2}(w_\ptdK-w_\ptK)$, which seems more natural
with respect to the assumptions on $\mathfrak{a}$; the power of
$\size(\Tau)$ in the denominator can be chosen arbitrarily. The
convergence of the scheme would remain true. The question of
optimal choice of the penalization operator is beyond the scope of
this paper.
\end{rem}

\subsection{Discrete convection operator}
Let ${\mathfrak f}:\R\to\R^d$ be continuous. Denote by
$\omega_M(\cdot)$ a modulus of continuity of ${\mathfrak f}$ on
$[-M,M]$, i.e., a continuous concave function on $[0,M]$ with
$\omega_M(0)=0$ and
\begin{equation*}
\max_{a,b\in [-M,M],|a-b|\leq r} \|{\mathfrak f}(a)-{\mathfrak
f}(b)\|\leq \omega_M(r).
\end{equation*}
Note that we can always choose $\omega_M$ strictly increasing,
upon replacing $\omega_M$ by $\omega_M+Id$ if needed.

Following Eymard, Gallou\"et, Herbin \cite{EyGaHe:book}, we now
define discrete convection fluxes, separately for each of the
meshes $\Mrond$, $\dMrond$. This will allow to discretize the
convective part of equation \eqref{eq:prob1}.

$\bullet$ Let $\KIL\in\Erond$. To approximate
$f(u)\cdot\nuKL$ by means of the two values $\uK,\uL$ that are available in
the neighbourhood of the interface $\KIL$, let us use some function $g_\ptKL$
of the couple $(\uK,\uL)\in\R^2$.
More exactly, take a collection of {\it numerical convection flux
functions} $(g_{\ptK,\ptL})_{\ptKIL\in\ptErond}$, $g_\ptKL\in
C(\R^2,\R)$, with the following properties:
\begin{equation}\label{Hypfluxes}
\begin{cases}
\text{(a) $g_\ptKL(\cdot,b)$ is nondecreasing for all $
b\in\R$,}\\
\hspace*{15pt} \text{and $g_\ptKL(a,\cdot)$ is nonincreasing
 for all $ a\in \R$};\\
\text{(b) $g_\ptKL(a,a)={\mathfrak f}(a)\cdot\nuKL$
for all $a\in \R$};\\
\text{(c) $g_\ptKL(a,b)=-g_{\ptK,\ptL}(b,a)$ $\;\forall
a,b\in\R$, for all neighbours $\K,\L\in\overline{\Mrond}$};
\\ \text{(d)
$g_\ptKL$ has the
same modulus of continuity as ${\mathfrak f}$, i.e.,}\\
\quad\text{ there exists $C$ independent of $\KIL$ such that  $\forall a,b,c,d\in [-M,M]$,}\\
 \dsp\quad\qquad  |g_\ptKL(a,b)-g_\ptKL(c,d)|\leq C\,
\bigl(\,\omega_M(|a-c|)+\omega_M(|b-d|)\,\bigr).\end{cases}
\end{equation}
These assumptions (see \cite{EyGaHe:book}) are by now standard.
Note that the assumption \eqref{Hypfluxes}(d) usually states that
${\mathfrak f},g_\ptKL$ are Lipschitz continuous, with the same
Lipschitz constant; here, we adapt it to the case of general
continuous function ${\mathfrak f}$.

Note that \eqref{Hypfluxes} (b) and (c) are compatible. Also note that the consistency requirement
\eqref{Hypfluxes}(b) together with the Green-Gauss formula imply
\begin{equation}\label{Hypfluxes-prop}
\sum\nolimits_{\ptL\in\ptNrond(\ptK)} \!\!\!\mKIL
g_\ptKL(a,a)={\mathfrak f}(a)\cdot \int_{\ptl\K}\!\!\nuK= 0 \quad
\text{for all $a\in\R$, for all $\K\in\Mrond$}.
\end{equation}
Practical examples of numerical convective flux functions can be
found in \cite{EyGaHe:book}.
These include the Godunov, Lax-Friedrichs, Engquist-Osher and
Rusanov fluxes as particular cases.

$\bullet$  Numerical convective flux functions $g_\ptdKdL$,
$\dKIdL \in \dErond$, are defined similarly.

 $\bullet$ On the set $\R^\ptbarTau$ of discrete functions $\ubarTau$ on $\overline{\Om}$, we define the operator
 $ (\div_c {\mathfrak f})^\ptTau [\cdot]$ of {\it discrete convection} by
\begin{equation*}
(\div_c {\mathfrak f})^\ptTau:\ubarTau\in \R^\ptbarTau  \mapsto
v^\ptTau=(\div_c {\mathfrak f})^\ptTau [\ubarTau] \in \R^\ptTau,
\end{equation*}
 where the discrete function
$v^\ptTau=\Bigl( v^\ptMrond,v^\ptdMrond \Bigr)$ on $\Om$ is given by
\begin{equation*}
v^\ptMrond=(v_\ptK)_{\ptK\in\ptMrond} \quad \text{with} \quad
v_\ptK=\frac{1}{\mK}\sum\nolimits_{\ptL\in\ptNrond(\ptK)} \mKIL
g_\ptKL(\uK,\uL);
\end{equation*}
\begin{equation*}
v^\ptdMrond=(v_\ptdK)_{\ptdK\in\ptdMrond}, \quad  \quad
v_\ptdK=\frac{1}{\mdK}\sum\nolimits_{\ptdL\in\ptdNrond(\ptdK)}
\mdKIdL g_\ptdKdL(\udK,\udL).
\end{equation*}

\subsection{Projection operators and test functions}\label{sec:TestFunctions}~

$\bullet$ On $ L^1(\Om)$, we define the {\it mesh projection operator}
$\mathbb{P}^\ptTau[\cdot]$ on the space of discrete functions on $\Om$ by
\begin{equation*}
\mathbb{P}^\ptTau :{\EuScript S}\in L^1(\Om) \mapsto {\EuScript
S}^\ptTau=\mathbb{P}^\ptTau[{\EuScript S}] \in \R^{\ptTau},
\end{equation*}
where the discrete function ${\EuScript S}^\ptTau=
\Bigl({\EuScript S}^{\ptMrond},{\EuScript S}^{\ptdMrond}\Bigr)$ on $\Om$ is given by
\begin{equation}\label{DiscrProjM}
\begin{array}{l}
\dsp {\EuScript S}^{\ptMrond}=({\EuScript
S}_\ptK)_{\ptK\in\ptMrond} \; \text{~with~} \;
{\EuScript S}_\ptK=\frac{1}{\mK}\int_\ptK {\EuScript S}(x)\,dx;\\[10pt]
\dsp {\EuScript S}^\ptdMrond=({\EuScript S}^n_\ptdK)_{\ptdK\in\ptdMrond} \; \text{~with~}
\; {\EuScript S}_\ptdK=\frac{1}{\mdK}\int_\ptdK {\EuScript
S}(x)\,dx.
\end{array}
\end{equation}

$\bullet$ For a sufficiently regular function $\psi$ on
$\overline{\Om}$, we will often employ the notations
$\psi^\ptTau=\mathbb{P}^\ptTau[\psi]$ and
$(\grad\psi)^\ptTau=\mathbb{P}^\ptTau[\grad\psi]$
($\grad \psi$ being $\R^d$-valued, the projection
is taken component per component).
Further, for $\KIL\in\Erond$ and
$\dKIdL\in\dErond$, we introduce
\begin{equation}\label{PsiSigma}
\dsp \psi_\ptKIL=\frac{1}{\mKIL}\int_\ptKIL \psi, \quad \dsp
\psi_\ptdKIdL=\frac{1}{\mdKIdL}\int_\ptdKIdL \psi.
\end{equation}

For $\L\in\ptl\Mrond$, there exists $\KIL\subset\ptl\Om$
that coincides with $\L$; in this case we assign
$\psi_\ptL=\psi_\ptKIL$. If $\psi|_{\ptl \Om}=0$, we have
$\psi_\ptL=0$ for all $\L\in\ptl\Mrond$.
 For $\dL\in\ptl\dMrond$, we assign
$\psi_\ptdL=\frac{1}{m_\ptdL}\int_\ptdL \psi$. If $\psi$ has a
compact support in $\Om$ and $\size(\Tau)$ is small enough, we
have $\psi_\ptdL=0$ for all $\dL\in\ptl\dMrond$.

Combining the above notation, we write
$\psi^\ptbarTau=\Bigl(\mathbb{P}[\psi],(\psi_\ptK)_{\ptK\in\ptl\ptMrond},
(\psi_\ptdK)_{\ptdK\in \ptl\ptdMrond} \Bigr)$ for the projection of a sufficiently regular
function $\psi$ on the space $\R^\ptbarTau$, and denote the
corresponding projection operator by $\mathbb{P}^\ptbarTau$.

\smallskip
\subsection{Dependency on $t$ and further notation}~

$\bullet$ Let $\Tau$ be a DDFV mesh as described above. Let $\Delt>0$
be the time discretization step. Set
$h=\max\{\size(\Tau),\Delt\}$. By convention, we will use $h$
as the parameter for a sequence of finite volume
schemes; our interest lies in studying convergence of
corresponding discrete solutions as $h\downarrow 0$.

Denote by $N$ the integer part of $T/\Delt$. In the sequel, in our
notation we omit the dependency of $N$, $\Tau$ and $\Delt$ on $h$.

$\bullet$ For a  functional space $X$ on $\Om$, we denote by
$\mathbb{S}^\Delt$ the projection operator
\begin{equation}\label{TimeProject}
\mathbb{S}^\Delt:{\EuScript S}\in L^1(0,T;X) \mapsto ({\EuScript
S}^n)_{n=1,\ldots,N}\in (X)^N, \qquad {\EuScript
S}^n=\frac{1}{\Delt}\int_{(n-1)\Delt}^{n\Delt} {\EuScript
S}(t)\,dt.
\end{equation}

$\bullet$ A {\it discrete function on $Q$} is a set
$u^{\ptTau,\Delt}=(u^{\ptTau,n})_{n=1,\ldots,N}$, where for each $n$,
$u^{\ptTau,n}$ is a discrete function on $\Om$. The set of all
such functions is denoted  $\R^{N\times \ptTau}$.

A {\it discrete function on $\overline{Q}$} is a set
$u^{\ptbarTau,\barDelt }=(u^{\ptbarTau,n})_{n=0,\ldots,N}$, where for
each $n$, $u^{\ptbarTau,n}$ is a discrete function on
$\overline{\Om}$. The set of all such
functions is denoted by $\R^{(N+1)\times \ptbarTau}$.
We also use discrete functions
$u^{\ptbarTau,\Delt}\in \R^{N\times \ptbarTau}$ and $u^{\ptTau,\barDelt }\in\R^{(N+1)\times \ptTau}$.
Each of $u^{\ptTau,\Delt},u^{\ptbarTau,\Delt},u^{\ptTau,\barDelt}$ is
therefore a restriction of $u^{\ptbarTau,\barDelt}$.
The entries of $u^{\ptbarTau,n}$ are denoted by $\uKn$ (respectively,
$\udKn$) for $\K\in\Mrond\cup \ptl\Mrond$ (respectively,
for $\dK\in\dMrond \cup \ptl\dMrond$).

A {\it discrete field on $Q$} is a set $\Frond^{\ptTau,\Delt}
=\Bigl(\Frond^{\ptTau,n} \Bigr)_{n=1,\ldots,N}$
where for each $n$, $\Frond^{\ptTau,n}$ is a discrete field on
$\Om$. The set of all such fields is denoted by $(\R^d)^{N\times
\ptDrond}$.

$\bullet$ Any discrete function can be composed with a mapping
$A:\R\to\R^m$, $m\in\N$; for instance, $A(u^{\ptbarTau,\barDelt})$
stands for $w^{\ptbarTau,\barDelt}$ with values $w_\K^n=A(u_\K^n)$
for $\K\in \overline{\Mrond}$ and $w_\dK^n=A(u_\dK^n)$ for $\dK\in
\overline{\dMrond}$, for $n=0,\ldots,N$. Similarly, any discrete field
can be composed with a mapping $\ph:\R^d\to\R^m$; one has
$\ph(\Frond^\ptTau)=\Bigl(\ph(\Frond_\ptD)\Bigr)_{\ptD\in
\ptDrond}$.

$\bullet$ We say that a discrete function is nonnegative ( respectively,
nonpositive), if all its entries are nonnegative (respectively,
nonpositive); e.g., for $v^\ptTau\in \R^\ptTau$ the notation
$v^\ptTau \geq 0$ means that $v_\ptK\geq 0$ for all $\K\in\Mrond$
and $v_\ptdK\geq 0$ for all $\dK\in\dMrond$.

\subsection{The finite volume scheme}
With the notation introduced above, the finite volume
discretization of problem \eqref{eq:prob1} takes the following compact form:
\[ \text{ find a discrete function $u^{\ptbarTau,\barDelt}$ on
$\overline{Q}$ satisfying for $n=1, \ldots,N$ the equations }
\]
\begin{equation}\label{AbstractScheme}
\left|\begin{array}{l} 
\dsp \frac{u^{\ptTau,n} - u^{\ptTau,(n-1)}}{\Delt} + (\div_c {\mathfrak
f})^\ptTau[u^{\ptbarTau,n}]-\div^\ptTau[{\mathfrak a}(\grad^\ptTau
w^{\ptbarTau,n})] \\[8pt]
\dsp \hspace*{145pt}
+{\Prond}^\ptTau[w^{\ptbarTau,n}]=\mathbb{P}^\ptTau(\mathbb{S}^\Delt[{\EuScript S}])^n,\\
\dsp  w^{\ptbarTau,n}=A(u^{\ptbarTau,n}),
\end{array}\right.
\end{equation}
$$
\text{together with the boundary and initial conditions}
$$
\begin{equation}\label{DiscrBC}
\text{for all $n=1,\ldots,N$,} \qquad
\begin{cases}
    \text{$u_\ptK^n=0$}  \quad \text{ for all $\K\in\ptl\Mrond$}\\
    \text{$u_\ptdK^n=0$}  \quad \text{for all $\dK\in\ptl\dMrond$};
\end{cases}
\end{equation}
\begin{equation}\label{DiscrIC}
\begin{cases}
    \text{$u_\ptK^0=\frac{1}{\mK}\int_\ptK u_0$}
    \quad\quad \,
    \text{for all $\K\in\Mrond$}\\
    \text{$u_\ptdK^0=\frac{1}{\mdK}\int_\ptdK u_0$}  \quad \text{for
    all $\dK\in\dMrond$}.
 \end{cases}
\end{equation}

Let us state \eqref{AbstractScheme} in a more explicit form:\\
\begin{equation*}
\begin{cases}
    \text{for all} \quad n=1,\ldots,N, \\[10pt]
    \quad \mK\frac{u^n_\ptK - u^{(n-1)}_\ptK}{\Delt}\\[10pt]
    \qquad + \;\sum_{\ptL\in\ptNrond(\ptK)}\! \mKIL g_\ptKL(u^n_\ptK,u^n_\ptL)
    - \sum_{\ptS\in\ptVrond(\ptK)}\!
    m_\edgeS {\mathfrak a}(\grad_\ptS A(u^{\ptbarTau,n})) \cdot \nuK \\[15pt]
    \qquad \dsp +\,\frac{d\!-\!1}{\size(\Tau)} \,
    \sum_{\ptdK\in\overline{\ptdMrond}}
    m_{\ptK\cap\ptdK}(A(u_\ptK)-A(u_\ptdK))= \mK {\EuScript S}^n_\ptK,
    \qquad \text{for all $\K\in\ptl\Mrond$},\\[15pt]
    \quad\mdK \frac{u^n_\ptdK - u^{(n-1)}_\ptdK}{\Delt}\\[10pt]
    \qquad +\;
    \sum_{\ptdL\in\ptdNrond(\ptdK)}\! \mdKIdL g_\ptdKdL(u^n_\ptdK,u^n_\ptdL)
    -
    \sum_{\ptS\in\ptdVrond(\ptdK)}\!
    m_\dedgeS {\mathfrak a}(\grad_\ptS A(u^{\ptbarTau,n})) \cdot
    \nudK \\[15pt]
    \qquad \dsp +\,\frac{1}{\size(\Tau)} \,\sum_{\ptK\in\overline{\ptMrond}}
    m_{\ptK\cap\ptdK}(A(u_\ptdK)-A(u_\ptK))= \mdK {\EuScript S}^n_\ptdK,
    \qquad \text{for all $\dK\in\ptl\dMrond$}.
\end{cases}
\end{equation*}
~\\

\noindent Here ${\EuScript S}^n_\ptK,{\EuScript S}^n_\ptL$ are
given by \eqref{DiscrProjM},\eqref{TimeProject};
$g_\ptKL,g_\ptdKdL$ are some numerical convection fuxes satisfying
\eqref{Hypfluxes}; $\nuK,\nudK$ for $\SDM$ given have the same
meaning as in \eqref{DiscrDivM},\eqref{DiscrDivdM}; finally, for
$\SDM$ given such that $\SDM\subset \DM\in \Drond$, $\grad_S
A(u^{\ptbarTau,n})$ is the vector of $\R^d$ constructed from the
values $w_\ptK=A(u^n_\ptK)$, $w_\ptdK=A(u^n_\ptdK)$ by formulas
\eqref{DiscrGradD2} (for $d=2$) or \eqref{DiscrGradD3} (for
$d=3$), i.e., in the  way indicated in Remark
\ref{ExplicDiscrGrad}.

\section{Elements of discrete calculus for DDFV schemes}\label{sec:DiscreteDuality}

In this section, we list convenient formulations of various
summation-by-parts formulas and chain rules needed for
the analysis of the discrete problem \eqref{AbstractScheme}.

\subsection{ Discrete duality formulas for the diffusion terms}~

$\bullet$ Recall that $\R^\ptTau$ is the space of all discrete
functions on $\Om$. For $m\in \N$ and $w^\ptTau,v^\ptTau \in
\Bigl(\R^\ptTau\Bigr)^m$, set
\begin{equation}\label{ScalProdFunct}
\Bleft w^\ptTau,\,v^\ptTau   \Bright =
\frac{1}{d}\sum_{\ptK\in\ptMrond} \mK\;w_\ptK \cdot v_\ptK \;+\;
\frac{d-1}{d} \sum_{\ptdK\in\ptdMrond} \mdK\;w_\ptdK \cdot v_\ptdK
\end{equation}
(here $\,\cdot\,$ denotes the scalar product in $\R^m$); it is
clear that $\Bleft \cdot,\cdot\Bright$ is a scalar product on
$\Bigl(\R^\ptTau\Bigr)^m$. We will use it for $m=1$ or $m=d$.

$\bullet$  Recall that  $(\R^d)^\ptDrond$ is the space of all
discrete fields on $\Om$. For $\Frond^\ptTau,\Grond^\ptTau \in
(\R^d)^\ptDrond$, set
\begin{equation}\label{ScalProdFields}
\Aleft \Frond^\ptTau,\,\Grond^\ptTau   \Aright =
\sum_{\ptD\in\ptDrond} \mD\,\Frond_\ptD \cdot\Grond_\ptD;
\end{equation}
it is clear that $\Aleft \cdot,\cdot\Aright$ is a scalar product
on $(\R^d)^\ptDrond$.

A key property of DDFV schemes (see
\cite{DomOmnes,ABH``double''}) is the following discrete analogue
of the duality between the  $-\div [\cdot]$ and the $\grad[\cdot]$
operators; it is sometimes called the {\it discrete duality}
property for finite volumes.
\begin{prop}\label{DualityDiff}
Let  $v^\ptbarTau\in \R^\ptbarTau_0$ and $\Frond^\ptTau\in
(\R^d)^\ptDrond$. Then
$$
\Bleft - \div^\ptTau [\Frond^\ptTau] \, ,\, v^\ptTau
\Bright=\Aleft \Frond^\ptTau\, ,\, \grad^\ptTau v^\ptbarTau
\Aright.
$$
\end{prop}

\begin{proof}
The proof is straightforward, using the summation-by-parts
procedure. Let us give it for the case $d=2$. Note that for
$\DM=\SDM=\SDM^{\ptK,\ptL}_{\ptdK\!\!,\ptdL}$,
$\mS=\frac{1}{2}\mKIL\dKL=\frac{1}{2}\mdKIdL\ddKdL$. By
\eqref{ScalProdFunct}, by \eqref{DiscrDivM},\eqref{DiscrDivdM},
and finally by \eqref{DiscrGradD2},\eqref{ScalProdFields}, we get
\[
\begin{array}{l}
\dsp \Bleft - \div^\ptTau [\Frond^\ptTau] \, ,\, v^\ptTau
\Bright\\[10pt]
\dsp  = -\frac{1}{2}\sum_{\ptK\in\ptMrond}
\Biggl(\sum_{\ptS\in\ptVrond(\ptK)} m_\edgeS \Frond_\ptS \cdot
\nuK\Biggr)\; v_\ptK \,-\, \frac{1}{2} \sum_{\ptdK\in\ptdMrond}
\Biggl( \sum_{\ptS\in\ptdVrond(\ptdK)} m_\dedgeS \Frond_\ptS
\cdot \nudK\Biggr)\; v_\ptdK\\[10pt]
\dsp = -\frac{1}{2}\sum_{\ptK\in\overline{\ptMrond}}
\Biggl(\sum_{\ptS\in\ptVrond(\ptK)} m_\edgeS \Frond_\ptS \cdot
\nuK\Biggr)\; v_\ptK \,-\, \frac{1}{2}
\sum_{\ptdK\in\overline{\ptdMrond}} \Biggl(
\sum_{\ptS\in\ptdVrond(\ptdK)} m_\dedgeS \Frond_\ptS \cdot
\nudK\Biggr)\; v_\ptdK\\[10pt]
\dsp =
\frac{1}{2}\sum\nolimits_{\ptS\in\ptSrond,\ptS=\ptS^{\ptK,\ptL}_{\ptdK\!\!,\ptdL}}
\Frond_\ptS\,\cdot\,\Bigl( \mKIL(v_\ptL-v_\ptK)\,\nuKL +
\mdKIdL (v_\ptdL-v_\ptdK)\,\nudKdL \Bigr)\\[10pt]
\dsp
=\sum\nolimits_{\ptS\in\ptSrond,\ptS=\ptS^{\ptK,\ptL}_{\ptdK\!\!,\ptdL}}
\mS\;\Frond_\ptS\,\cdot\, \Bigl( \frac{v_\ptL-v_\ptK}{\dKL}\,\nuKL
+
\frac{v_\ptdL-v_\ptdK}{\ddKdL}\,\nudKdL\Bigr)\\[10pt]
\dsp
=\sum\nolimits_{\ptS\in\ptSrond}\mS\;\Frond_\ptS\cdot\grad_\ptS
v^\ptbarTau=\sum\nolimits_{\ptD\in\ptDrond}\mD\;\Frond_\ptD\cdot\grad_\ptD
v^\ptbarTau= \Aleft \Frond^\ptTau\, ,\, \grad^\ptTau v^\ptbarTau
\Aright.\\[-15pt]
\end{array}
\]
\end{proof}

Furthermore, we have the following ``entropy dissipation''
inequalities:
\begin{prop}\label{ParabEntropyDissipation}
 Let $u^\ptbarTau\in \R^\ptbarTau_0$
and $\psi\in {\mathcal D}(\overline{\Om})$, $\psi\geq 0$.
 Let $\theta:\R\to\R$ be a nondecreasing function. Assume that
\begin{equation}\label{CarrilloChoice}
 \text{either $\theta(0)=0$, or $\psi\in{\mathcal D}(\Om)$ and $\size(\Tau)$ is small enough}.
 \end{equation}
Denote $\psi^\ptbarTau={\mathbb{P}}^\ptbarTau[\psi]$. Then
\begin{equation}\label{EntrChainRuleDiff}
\begin{array}{l}
\dsp \Bleft \div^\ptTau \Bigl[ k(\grad^\ptTau
A(u^\ptbarTau))\grad^\ptTau A(u^\ptbarTau)\Bigr] \, ,\,
\theta(u^\ptTau) \psi^\ptTau \Bright\\[6pt]
\hspace*{45mm} \dsp \leq -\Aleft k\Bigl(\grad^\ptTau
A(u^\ptbarTau)\Bigr)\,\grad^\ptTau A_\theta(u^\ptbarTau) \, ,\,
\grad^\ptTau\psi^\ptbarTau \Aright.
\end{array}
\end{equation}
\end{prop}
\begin{rem}\normalfont \label{eq:ConditionsForEntropyIneq}
Note that the conformity of the meshes (see
Remark~\ref{rem:conformal}) is essential for this result, as well
as the particular form of ${\mathfrak a}$ and (for $d=3$)
condition \eqref{eq:AssumptMesh}.
\end{rem}

\begin{proof}
Let us treat  the left-hand side of \eqref{EntrChainRuleDiff} term
by term. It is the sum of generic terms
of the form $T_{\ptK,\ptS},T^*_{\ptdK,\ptS}$; here
$$
\begin{array}{l}
\dsp T_{\ptK,\ptS}=\frac{1}{d}\,\mK \frac{1}{\mK} \mKIL
k(\grad_{\!\ptS}\, A(u^\ptbarTau))\grad_{\!\ptS}\, A(u^\ptbarTau)
\cdot \nuKL \theta(u_\ptK)
\psi_\ptK \\[7pt]
\qquad\qquad \dsp =\frac{1}{d}\,\mKIL k\Bigl(\grad_{\!\ptS}\,
A(u^\ptbarTau)\Bigr)\;\grad_{\!\ptS}\, A(u^\ptbarTau) \cdot
\nuKL\; \theta(u_\ptK) \psi_\ptK
\end{array}
$$ with
$\SDM=\SDM^{\ptK,\ptL}\in\Vrond(\K)$.
 The notation $T^*_{\ptdK,\ptS}$ stands
for analogous terms involving $\dK$ and $\SDM\in \dVrond(\dK)$.
Notice that thanks to assumption \eqref{CarrilloChoice},
$\theta(u^\ptbarTau)\psi^\ptbarTau\in \R^\ptbarTau_0$, so that we
can also add the terms $T_{\ptK,\ptS},T^*_{\ptdK,\ptS}$
corresponding to $\K\in\ptl\Mrond$, $\dK\in\ptl\dMrond$,
respectively. The summation of $T_{\ptK,\ptS},T^*_{\ptdK,\ptS}$
therefore runs on all subdiamonds
$\SDM=\SDM^{\ptK,\ptL}_{\ptdK,\ptdL}\in\Srond$, with the
associated $\K,\L\in\overline{\Mrond}$,
$\dK,\dL\in\overline{\dMrond}$.

The  convexity argument yields
\begin{equation}\label{ConvexityIneq}
(A(z)-A(\hat z))\theta(\hat z)\leq A_\theta(z)-A_\theta(\hat z)
\qquad\text{for all $z,\hat z\in\R$}.
\end{equation}
By \eqref{DiscrGrad}, using the positivity of $\psi_\ptK$ and
applying inequality \eqref{ConvexityIneq}, we get
\begin{equation}\label{ApplConvexityIneq}
\begin{array}{l}
\dsp T_{\ptK,\ptS} = \frac{1}{d}\,\mKIL k\Bigl(\grad_{\!\ptS}\,
A(u^\ptbarTau)\Bigr)
\frac{A(u_\ptL)-A(u_\ptK)}{\dKL}\theta(u_\ptK) \psi_\ptK\\[7pt]
\dsp\qquad\qquad \leq
 \frac{1}{d}\,\mKIL
k\Bigl(\grad_{\!\ptS}\, A(u^\ptbarTau)\Bigr)
\frac{A_\theta(u_\ptL)-A_\theta(u_\ptK)}{\dKL} \psi_\ptK.
\end{array}
\end{equation}
The terms $T_{\ptdK,\ptS}$ are treated in the same way. Now by the
same computation as in the proof of Proposition~\ref{DualityDiff},
one shows that the right-hand sides of \eqref{ApplConvexityIneq}
and of the corresponding inequality for  $T_{\ptdK,\ptS}$ sum up
to yield the right-hand side of \eqref{EntrChainRuleDiff}. This
concludes the proof.
\end{proof}

\subsection{Summation formulas for the penalization terms}
For the penalization operator ${\Prond}^\ptTau$, we have the
following summation formulas.
\begin{lem}\label{SummationPenalize}
Let $w^\ptbarTau\in \R^\ptbarTau$ and $\psi^\ptbarTau\in
\R^\ptbarTau_0$. Then
 \begin{equation}\label{SummPenalWeak}
\Bleft \Prond [w^\ptbarTau], \psi^\ptTau \Bright=\frac{d-1}{d}
\sum_{\ptK\in\overline{\ptMrond}, \ptdK\in\overline{\ptdMrond}}
m_{\ptK\cap\ptdK}\frac{(w_\ptK-w_\ptdK)(\psi_\ptK-\psi_\ptdK)}{\size(\Tau)}.
 \end{equation}
Further, let $A,\theta:\R\to\R$ be nondecreasing. Assume
$u^\ptbarTau\in\R^\ptbarTau$ is such that $A(u^\ptbarTau)$ belongs
to $\R^\ptbarTau_0$. Let $\psi\in {\mathcal D}(\overline{\Om})$,
$\psi\geq 0$; denote $\psi^\ptbarTau=\mathbb{P}^\ptbarTau[\psi]$.
Assume \eqref{CarrilloChoice}. Then
 \begin{equation}\label{SummPenalWeakAvecTheta}
\Bleft \Prond [A(u^\ptbarTau)], \theta(u^\ptTau)\psi^\ptTau
\Bright \geq \frac{d\!-\!1}{d} \!\!\!\!\!\!\!\!\!\!
\sum_{\ptK\in\overline{\ptMrond}, \ptdK\in\overline{\ptdMrond}}
\!\!\!\!\!\!\!\!\!\!
m_{\ptK\cap\ptdK}\;\theta(u_\ptK)\,\frac{(A(u_\ptK)\!-\!A(u_\ptdK))(\psi_\ptK\!-\!\psi_\ptdK)}{\size(\Tau)}.
 \!\!\!\!\!\!\end{equation}
In both formulas
\eqref{SummPenalWeak},\eqref{SummPenalWeakAvecTheta}, the values
$\psi_\ptK,\psi_\ptdK$ for $\K\in
\overline{\Mrond},\dK\in\overline{\dMrond}$ are those of the
corresponding discrete function $\psi^\ptbarTau$.
\end{lem}
The proof is straightforward from the definitions of
$\Bleft\cdot,\cdot\Bright$ and $\Prond^\ptTau$, using the
summation-by-parts procedure.

\subsection{Discrete duality formulas for the evolution terms}~

\begin{lem}\label{DualityEvol}
Let $\theta:\R\to\R$ be a nondecreasing function, and $ \eta=\int \theta(s)\ds$ be its
primitive. Let $\psi\in {\mathcal D}(\overline{Q})$, $\psi\geq 0$.
Denote $\psi^{\ptTau,\Delt}=\mathbb{P}^\ptTau\circ
\mathbb{S}^\Delt[\psi]$. Then for all $u^{\ptTau,\barDelt} \in
\R^{(N+1)\times\ptTau}$  one has
\begin{align*}
    & \sum_{n=1}^N \Delt\Bleft \frac{u^{\ptTau,n} - u^{\ptTau,(n-1)}}{\Delt} \;,\;
    \theta(u^{\ptTau,n})\;\psi^{\ptTau,n} \Bright
    \\ & \qquad \geq  -\sum_{n=1}^{N-1} \Delt \Bleft \eta(u^{\ptTau,n})\;,\;
    \frac{\psi^{\ptTau,(n+1)} - \psi^{\ptTau,n}}{\Delt}\Bright
    \\ & \qquad\qquad\qquad
    + \Bleft \eta(u^{\ptTau,N})\;,\; \psi^{\ptTau,N} \Bright
    -\Bleft \eta(u^{\ptTau,0})\;,\; \psi^{\ptTau,1} \Bright.
\end{align*}
\end{lem}
\begin{proof}
The formula follows by the Abel transformation
combined with the convexity inequality: $ (z-\hat z)\theta( z)\geq
\eta(z)-\eta(\hat z) \quad\text{for all $z,\hat z\in\R$}. $
\end{proof}

\subsection{ Discrete duality formulas for the convection terms}
For the convection terms, we have a more involved ``entropy
dissipation'' duality formula. For later use, we state it in the
double framework, although each of the meshes $\Mrond$,$\dMrond$
is treated separately in the proof.

\begin{prop}\label{DualityConv}
Let $u^\ptbarTau\in \R^\ptbarTau_0$ and
$\psi\in {\mathcal D}(\overline{\Om})$.
Let $\theta:\R\to\R$ be a nondecreasing function.
Assume \eqref{CarrilloChoice}.
Consider the associated entropy-flux pair
$$
\eta=\int\theta(s)\,ds, \qquad \dsp{\mathfrak q}
= \theta\,{\mathfrak f}
-\int{\mathfrak f}(s)\,d\theta(s).
$$
Denote $\psi^\ptTau=\mathbb{P}^\ptTau[\psi]$ and
$(\grad \psi)^\ptTau=\mathbb{P}^\ptTau[\grad\psi]$.
One has
\begin{equation}\label{DualityConvect}
\begin{split}
    & \Bleft  (\div_c {\mathfrak f})^\ptTau [u^\ptbarTau] \, ,\,
    \theta(u^\ptTau)\psi^\ptTau \Bright
    \\ & = -\dsp\Bleft {\mathfrak q}(u^\ptTau)\,,\,(\grad \psi)^\ptTau \Bright
    +I_\theta[u^\ptMrond,\psi]
    +R_\theta[u^\ptMrond,\psi]+I^{*}_\theta[u^\ptdMrond,\psi]
    +R^{*}_\theta[u^\ptdMrond,\psi],
\end{split}
\end{equation}
where
\begin{equation}\label{SignedTerm}
I_\theta
[u^\ptMrond\!,\psi]=\frac{1}{d}\sum_{\ptKIL\in\ptErond}\!\!\!
\mKIL I_\theta^\ptKIL\, \psi_\ptKIL, \;
I^*_\theta[u^\ptdMrond\!\!,\psi]=\frac{d\!-\!1}{d}\!\!\!\!\sum_{\ptdKIdL\in\ptdErond}\!\!\!
\mdKIdL I_\theta^\ptdKIdL \psi_\ptdKIdL
\end{equation}
with
\begin{equation}\label{SignedTermEntries}
\begin{array}{l}
\dsp I^\ptKIL_\theta=\int^{u_\ptL}_{u_\ptK}
\bigl(g_\ptKL(s,s)-g_\ptKL(u_\ptK,u_\ptL)\bigr)\,
d\theta(s),\\[9pt]
\dsp I^\ptdKIdL_\theta=\int^{u_\ptdL}_{u_\ptdK}
\bigl(g_\ptdKdL(s,s)-g_\ptdKdL(u_\ptdK,u_\ptdL)\bigr)\,
d\theta(s).
\end{array}
\end{equation}

Further,  one has $I^\ptKIL_\theta\geq 0$ for all $\KIL\in\Erond$,
and the remainder term $R_\theta$ satisfies
\begin{equation}\label{TermRtheta}
 \Bigl|R_\theta [u^\ptMrond,\psi]\Bigr|\leq
\bigl(\max_{\ptK\in\ptMrond}
|\theta(u_\ptK)|\bigr)\;\sum_{\ptKIL\in\ptErond} \mKIL \,\Bigl(
R_\ptK^\ptKIL+R_\ptL^\ptKIL \Bigr) \,
\Bigl(|\psi_\ptK\!-\!\psi_\ptKIL|+|\psi_\ptL\!-\!\psi_\ptKIL|\Bigr),
\end{equation}
\begin{equation}\label{TermRthetaEntries}
R_\ptK^\ptKIL=|g_\ptKL(u_\ptK,u_\ptK)-g_\ptKL(u_\ptK,u_\ptL)|,
\quad
R_\ptL^\ptKIL=|g_\ptKL(u_\ptL,u_\ptL)-g_\ptKL(u_\ptK,u_\ptL)|.
\end{equation}

Similarly, one has $I^\ptdKIdL_\theta\geq 0$ for all
$\dKIdL\in\dErond$, and the remainder term $R^*_\theta$ satisfies
the analogue of \eqref{TermRtheta}, \eqref{TermRthetaEntries} with
$\K,\L,\Mrond,\Erond$ replaced by $\dK,\dL,\dMrond,\dErond$.
\end{prop}

Note that our notation is consistent: we have
$I_\theta^\ptKIL=I_\theta^{\ptL\ptI\ptK}$,
$R_\ptK^\ptKIL=R_\ptL^{\ptL\ptI\ptK}$ for all neighbours $\K,\L$
(for dual neighbours $\dK,\dL$, similar identities hold).

\begin{proof}
We exploit the ideas of \cite{EyGaHe:book} and \cite{Carrillo}.

Thanks to \eqref{CarrilloChoice} and because $u^\ptbarTau$ is zero
on boundary volumes, we have
\begin{equation}\label{CarrilloBoundaryValues}
  \text{$\eta(u_\ptL)\psi_\ptL=0$ for all $\L\in \ptl\Mrond$};
  \quad \text{$\eta(u_\ptdL)\psi_\ptdL=0$ for all $\dL\in
  \ptl\dMrond$}.
\end{equation}

Separating the contributions of $\Mrond$ and $\dMrond$, we write the
left-hand side of \eqref{DualityConvect} as
$\frac{1}{d}I+\frac{d-1}{d}I^*$, where
\begin{equation}
I:= \sum_{\ptK\in\ptMrond} \mK \Bigl(\frac{1}{\mK}
\sum_{\ptL\in\ptNrond(\ptK)} \mKIL g_\ptKL(u_\ptK,u_\ptL)
\Bigr)\,\theta(u_\ptK)\psi_\ptK.
\end{equation}
Applying \eqref{Hypfluxes}(c) and \eqref{Hypfluxes-prop}, using
\eqref{CarrilloBoundaryValues} in the summation-by-parts procedure, we get
$$
\begin{array}{l}
\dsp I=\sum_{\ptKIL\in\ptErond} \mKIL \Biggl(
\theta(u_\ptL)\Bigr(g_\ptKL(u_\ptL,u_\ptL)-g_\ptKL(u_\ptK,u_\ptL)\Bigr)\,\psi_\ptL
\\[5pt]
\dsp\hspace*{40mm}
-\theta(u_\ptK)\Bigr(g_\ptKL(u_\ptK,u_\ptK)-g_\ptKL(u_\ptK,u_\ptL)\Bigr)\,\psi_\ptK
\Biggr).
\end{array}
$$
Hence, choosing $\psi_\ptKIL$ as defined in \eqref{PsiSigma}, we have
$$
\begin{array}{l}
\dsp I =\sum_{\ptKIL\in\ptErond} \mKIL
\Biggl(\theta(u_\ptL)\Bigr(g_\ptKL(u_\ptL,u_\ptL)-g_\ptKL(u_\ptK,u_\ptL)\Bigr)\\[5pt]
\dsp\hspace*{45mm} -
\theta(u_\ptK)\Bigr(g_\ptKL(u_\ptK,u_\ptK)-g_\ptKL(u_\ptK,u_\ptL)\Bigr)
\Biggr)\,\psi_\ptKIL\\[5pt]
\dsp\qquad +\sum_{\ptKIL\in\ptErond} \mKIL \Biggl(
\theta(u_\ptL)\Bigr(g_\ptKL(u_\ptL,u_\ptL)-g_\ptKL(u_\ptK,u_\ptL)\Bigr)(\psi_\ptL\!-\!\psi_\ptKIL)\\[5pt]
\dsp\hspace*{45mm}
-\theta(u_\ptK)\Bigr(g_\ptKL(u_\ptK,u_\ptK)-g_\ptKL(u_\ptK,u_\ptL)\Bigr)\,
(\psi_\ptK\!-\!\psi_\ptKIL) \Biggr).
\end{array}
$$
Now recall that $g=g_\ptKL$ satisfies \eqref{Hypfluxes}(b). Thus
the following integration-by-parts formula holds true:
\begin{equation*}
\begin{array}{l}
\dsp ({\mathfrak q}(b)\!-\!{\mathfrak
q}(a))\cdot\nuKL=\Bigl(\theta(b){\mathfrak
f}(b)\!-\!\theta(a){\mathfrak f}(a)\! -\! \int_a^b\!
{\mathfrak f}(s)\,d\theta(s)\Bigr)\cdot\nuKL\\[5pt]
\dsp =\theta(b)\bigl(g(b,b)\!-\!g(a,b)\bigr)-\theta(a)
\bigl(g(a,a)\!-\!g(a,b)\bigr)-\!\int_a^b\!
\bigl(g(s,s)\!-\!g(a,b)\bigr)\, d\theta(s).
\end{array}
\end{equation*}
We deduce $I=J+I_\theta+R_\theta$, where
\begin{align*}
    J&=\sum_{\ptKIL\in\ptErond} \mKIL \Bigl({\mathfrak q}(u_\ptK)
    -{\mathfrak q}(u_\ptL)\Bigr)\cdot\nuKL \psi_\ptKIL
    =\sum_{\ptK\in\ptMrond} {\mathfrak q}(u_\ptK)\cdot
    \Bigl(\sum_{\ptL\in\ptNrond(\ptK)} \mKIL \psi_\ptKIL\, \nuKL\Bigr)\\
    & =\sum_{\ptK\in\ptMrond} {\mathfrak q}(u_\ptK)\cdot
    \int_{\ptl \ptK}\psi\, \nuK = \sum_{\ptK\in\ptMrond} \int_\ptK
    \div\Bigl({\mathfrak q}(u_\ptK) \psi\Bigr)=\sum_{\ptK\in\ptMrond}
    \int_\ptK {\mathfrak q}(u_\ptK)\cdot \grad \psi\\
    & =\sum_{\ptK\in\ptMrond} \mK\,{\mathfrak q}(u_\ptK)
    \cdot (\grad \psi)_\ptK,
\end{align*}
and
$$ I_\theta= \sum_{\ptKIL\in\ptErond} \mKIL
\Bigl(\int^{u_\ptL}_{u_\ptK}
\bigl(g_\ptKL(s,s)-g_\ptKL(u_\ptK,u_\ptL)\bigr)\,
d\theta(s)\Bigr)\, \psi_\ptKIL,
$$
$$
\begin{array}{l}
\dsp |R_\theta|\leq \sum_{\ptKIL\in\ptErond} \mKIL \,\Biggl(
|g_\ptKL(u_\ptK,u_\ptK)-g_\ptKL(u_\ptK,u_\ptL)| +
|g_\ptKL(u_\ptL,u_\ptL)-g_\ptKL(u_\ptK,u_\ptL)|\Biggr) \,\\[9pt]
\dsp \hspace*{45mm} \times
\Biggl(|\psi_\ptK\!-\!\psi_\ptKIL|+|\psi_\ptL\!-\!\psi_\ptKIL|\Biggr)\;\times\;\bigl(\max_{\ptK\in\ptMrond}
|\theta(u_\ptK)|\bigr).
\end{array}$$
 In the same way, $I^*=J^*+I^*_\theta+R^*_\theta$ with
analogous estimates. We have the equality
$\frac{1}{d}\,J+\frac{d-1}{d}\,J^*=\dsp\Bleft {\mathfrak
q}(u^\ptTau)\,,\,(\grad \psi)^\ptTau \Bright$.
 With the notation of
\eqref{SignedTerm}-\eqref{TermRthetaEntries}  the result of the
proposition follows.
\end{proof}

\section{Properties of discrete operators and functional spaces}\label{sec:DiscreteSpaces}

In this section we state important embedding and compactness
properties of spaces of discrete functions, as well as the
asymptotic (as $h\to 0$) properties of various discrete operators.

\subsection{Discrete functions and fields as elements
of Lebesgue spaces}\label{EmbedInLebesgue} ~

For any $E\subset\overline{Q}$, denote by $\char_E$ its
characteristic function.

For $n=1,\ldots,N$, set
\begin{align*}
    \QKn&=[(n-1)\Delt,n\Delt[\times\K, \qquad \quad
    \text{ for $\K\in \Mrond$};\\
    \QdKn&=[(n-1)\Delt,n\Delt[\times\dK, \qquad
    \text{ for $\dK\in\dMrond$};\\
    \QDn & =[(n-1)\Delt,n\Delt[\times\DM,
    \qquad \quad \text{ for $\DM\in\Drond$}.
\end{align*}
For a discrete function $v^{\ptTau,\Delt}$ on $Q$, denote by
$v^{\ptMrond,\Delt}$ (respectively, by $v^{\ptdMrond,\Delt}$)
the piecewise constant function
\begin{align*}
    &v^{\ptMrond,\Delt}(t,x)=\!\sum_{n=1}^{N}\sum_{\ptK\in\ptMrond}\!
    u^n_\ptK \char_\QKn(t,x) \\
    & \qquad \Biggl(\text{respectively,~~}
    v^{\ptdMrond,\Delt}(t,x)
    =\!\!\sum_{n=1}^{N}\sum_{\ptdK\in\ptdMrond}\!\!
    u^n_\ptdK \char_\QdKn(t,x) \Biggr).
\end{align*}
Whenever it is convenient, we identify the discrete function
$v^{\ptTau,\Delt}\in\R^{N\times\ptTau}$ with the function on $Q$
given by
$$
v^{\ptTau,\Delt}(t,x)=\frac{1}{d}\;v^{\ptMrond,\Delt}(t,x)+\frac{d\!-\!1}{d}\;v^{\ptdMrond,\Delt}(t,x).
$$
In a similar way, we identify a discrete field $\Frond^{\ptTau,\Delt}\in \R^{N\times\ptDrond}$
on $Q$ with the function
$$
\Frond^{\ptTau,\Delt}(t,x)=\sum_{n=1}^{N}\sum_{\ptD\in\ptDrond}
\Frond_{\ptD}^{n} \char_\QDn(t,x).
$$

Analogous conventions apply to time-independent discrete functions
and discrete fields, in which case we suppress the superscript
$\Delt$ in the notation.

\subsection{Consistency properties of discrete operators}
In the proposition below we show the consistency properties of
the projection and discrete gradient operators in Lebesgue spaces.
Also note the property $(iv)$, which, combined with formula
\eqref{SummPenalWeak}, expresses the fact that the penalization
operator introduced in Section~\ref{sec:Penalization}
vanishes (in an appropriate sense) as $\size(\Tau)\to 0$.

\begin{prop}\label{PropConsistency}
Let $\Tau$ be a double mesh of $\Om$, $\Delt>0$,
$h=\max\{\size(\Tau),\Delt\}$, and $q\in[1,+\infty]$.
Then \\[5pt]
(i)~~ there exists a constant $C$ that only depends on $\Omega$,
$q$ and $\reg(\Tau)$ such that
 $$ \text{$\forall w\in L^q(Q)$,
$\|(\mathbb{P}^\ptTau\circ\mathbb{S}^\Delt
w)^{\ptMrond,\Delt}\|_{L^q}+
\|(\mathbb{P}^\ptTau\circ\mathbb{S}^\Delt
w)^{\ptdMrond,\Delt}\|_{L^q} \leq C\,\|w\|_{L^q},$ }$$
  and
 $$ \text{$\forall w\in L^q(0,T;W^{1,q}_0(\Om)))$,
$\|\grad^\ptTau\mathbb{P}^\ptbarTau\circ\mathbb{S}^\Delt
w\|_{L^q}\leq C\,\|\grad w\|_{L^q}$;}
$$
 (ii)~~
for all $w\in L^q(Q)$, $q<+\infty$, both $(\mathbb{P}^\ptTau\circ
\mathbb{S}^\Delt w)^{\ptMrond,\Delt}$ and $
(\mathbb{P}^\ptTau\circ \mathbb{S}^\Delt w)^{\ptdMrond,\Delt}$
converge to $w$ in $L^q(Q)$ as $h\to 0$;
\\[5pt]
(iii)~~ for all $w\in L^q(0,T;W^{1,q}_0(\Om))$,  $q<+\infty$, the
discrete fields $\grad^\ptTau \mathbb{P}^\ptbarTau\circ
\mathbb{S}^\Delt w$ converge to $\grad w$ in $(L^q(Q))^d$ as $h\to 0$; \\[5pt]
(iv)~~ let $\psi\in {\mathcal D}(\overline{\Om})$, and
$\psi^{\ptbarTau,n}=\mathbb{P}^\ptbarTau(\mathbb{S}[\psi])^n$,
$n=1,\ldots,N$. There exists a constant $C$ that only depends on $Q$ and
$\reg(\Tau)$ such that
$$\dsp\sum_{n=1}^N \Delt\;
\sum_{\ptK\in\overline{\ptMrond}, \ptdK\in\overline{\ptdMrond}}
m_{\ptK\cap\ptdK}\frac{(\psi_\ptK-\psi_\ptdK)^2}{\size(\Tau)} \leq
C\, \|\grad \psi\|_{L^\infty}\times \size(\Tau).
$$
\end{prop}
\begin{proof}
The proof of (i)-(iii) is a straightforward 
generalization of \cite[Lemma 3.3, Proposition 3.4
and Corollary 3.5]{ABH``double''}. We need to take into account
the fact that $\|\mathbb{S}^\Delt w\|_{L^q(X)}\leq \|w\|_{L^q(X)}$
and (for $q\neq+\infty$) $\|\mathbb{S}^\Delt w -w\|_{L^q(X)}
\to 0$ as $\Delt \to 0$  for all $w\in L^q(0,T;X)$, where $X$
stands for $L^q(\Om)$ or for $W^{1,q}_0(\Om)$. 
Remark~\ref{rem:consistAffines} is important for (iii) (thus, the
Delaunay property of $\Mrond$ is used). 
Further, in a standard way similar to \cite[Lemma 3.3]{ABH``double''}
one proves that for all $\K\in\overline{\Mrond}, \dK\in\overline{\dMrond}$ such that
$\K\cap\dK\neq \text{\O}$, one has
$|\psi^n_\ptK-\psi^n_\ptdK|\leq C(\reg(\Tau)) \|\grad
\psi\|_{L^\infty}\times \size(\Tau)$  for all $n=1,\ldots,N$.
Hence the claim (iv) follows.
\end{proof}

\subsection{Discrete embedding and compactness results}
Next we state a version of the Poincar\'e inequality and an
embedding-kind translation estimate on double discrete functions.
\begin{prop}\label{PoincareAndTranslations}
Assume $\Tau$ is a double mesh on $\Om$, $\Delt>0$. Let
$q\in[1,+\infty)$. There exists a constant $C>0$ that only depends
on $\diam(\Omega)$ and $q$ such that
\\[5pt]
 (i)~~~  for all
$w^{\ptbarTau,\Delt}\in \R^{N\times\ptbarTau}_0$ one has
$\|w^{\ptMrond,\Delt}\|_{L^q}+\|w^{\ptdMrond,\Delt}\|_{L^q} \leq C\; 
\|\grad^\ptTau w^{\ptbarTau,\Delt} \|_q;$\\[5pt]
(ii)~~~  for all $w^{\ptbarTau}\in \R^{\ptbarTau}_0$, for all
$\Del\in \R^d$ one has
$$
\|w^{\ptMrond}(\cdot+\Del)-w^{\ptMrond}(\cdot)\|_{L^q}
+\|w^{\ptdMrond}(\cdot+\Del)-w^{\ptdMrond}(\cdot)\|_{L^q} \leq
C\;\|\grad^\ptTau w^{\ptbarTau}\|_{L^q}\times |\Del|^{1/q}.
$$
\end{prop}

\begin{proof}~
The proof follows the lines of \cite[Lemma 1]{AGW} and \cite[Lemma
3.6]{ABH``double''}. Note that if $d=3$, the fact that all
interfaces $\KIL$ are triangles plays an important role in the
proof.
\end{proof}
Here is the asymptotic compactness result for ``discrete
$L^p(0,T;W^{1,p}_0(\Om)$'' spaces.
\begin{prop}\label{DiscrRellich} Let $p\in(1,+\infty)$.
Assume we are given a family $\{w^{\ptbarTau,\Delt}\}_h$ of
discrete functions in $\R^{N\times\bar\ptTau}_0 $ corresponding to
a family of double meshes $\Tau$ such that $\reg(\Tau)$ is
uniformly bounded (recall that we parametrize the meshes by
$h=\max\{\size(\Tau),\Delt\}$).\\[5pt]
 (i)~~~
Assume that there exists a constant $C>0$ such that
$$\|\grad^{\ptTau} w^{\ptbarTau,\Delt}\|_{L^p}\leq C. $$ Then
there exists a (not labelled) sequence of meshes such that as
$h\to 0$
$$
\text{
$w^{\ptTau,\Delt}=\frac{1}{d}\,w^{\ptMrond,\Delt}+\frac{d\!-\!1}{d}\,w^{\ptdMrond,\Delt}$
converge weakly in $L^p(Q)$ to some limit $w$;}$$
furthermore, $w\in L^p(0,T;W^{1,p}_0(\Om))$ and
$$\text{  the discrete fields
$\grad^{\ptTau} w^{\ptbarTau,\Delt}$ converge weakly in
$(L^p(Q))^d$ to $\grad w$ as $h\to 0$}.
$$
(ii)~~~ If, in addition,
$$\dsp \sum_{n=1}^N \Delt\Bleft \Prond^{\ptTau}[w^{\ptbarTau,n}], w^{ \ptTau,n}
\Bright \leq C,$$
 where $\Prond^{\ptTau}$ are the penalization
operators introduced in Section~\ref{sec:Penalization}, then
$$
\text{both $w^{\ptMrond,\Delt}$ and $w^{\ptdMrond,\Delt}$ converge
to $w$ weakly in $L^p(Q)$ as $h\to 0$}.$$
\end{prop}

\begin{rem}\normalfont
Note that upon providing uniform estimates on time translates of
$w^{\ptTau,\Delt}$ in $L^p(Q)$, strong convergence to $w$ in
$L^p(Q)$ holds true (see Section~\ref{sec:convergence}).
\end{rem}

\begin{proof} ~(i)~~~ The proof is very similar to the one of \cite[Lemma 3.8]{ABH``double''}.

First, by Proposition~\ref{PoincareAndTranslations}(i), both
families $\{w^{\ptMrond,\Delt}\}_h,\{w^{\ptdMrond,\Delt}\}_h$ of
components of $w^{\ptTau,\Delt}$ are bounded in $L^p(Q)$.
Therefore we can choose a common sequence such that both
components converge weakly in $L^p(Q)$. Also
$w^{\ptTau,\Delt}=\frac{1}{d}\,w^{\ptMrond,\Delt}+\frac{d\!-\!1}{d}\,w^{\ptdMrond,\Delt}$
converge weakly to some limit that we denote $w$. We can also
assume that the corresponding sequence $\{\grad^\ptTau
w^{\ptbarTau,\Delt}\}_h$ converges weakly in $(L^p(Q))^d$ to some
limit $\chi$. Let us show that $w\in L^p(0,T;W^{1,p}_0(\Om))$ and
$\chi=\grad w$.

Take any field $\Frond\in (L^{p'}(0,T;W^{1,p'}(\Om))^d$.
Denote by $\Frond^{\ptTau,\Delt}$
the discrete field on $Q$ with entries
$$\Frond^n_\ptD=\frac{1}{\Delt\times
\mD}\int_{(n\!-\!1)\Delt}^{n\Delt}\int_\ptD \Frond.$$  Denote by
$(\div \Frond)^{\ptTau,\Delt}$ the discrete function
$\mathbb{P}^\ptTau\circ\mathbb{S}^\Delt[\div \Frond]$ on $Q$,
which has the entries
$$(\div
\Frond)^n_\ptK=\frac{1}{\Delt\mK}\int_{(n\!-\!1)\Delt}^{n\Delt}\int_\ptK
\div\Frond=\frac{1}{\Delt\mK}\int_{(n\!-\!1)\Delt}^{n\Delt}
\sum\nolimits_{\ptS\in\ptVrond(\ptK)} \int_\edgeS \Frond \cdot
\nuK,
$$
$$ (\div
\Frond)^n_\ptdK=\frac{1}{\Delt\mdK}\int_{(n\!-\!1)\Delt}^{n\Delt}\int_\ptdK
\div\Frond=
\frac{1}{\Delt\mdK}\int_{(n\!-\!1)\Delt}^{n\Delt}\sum\nolimits_{\ptS\in\ptdVrond(\ptdK)}
\int_\dedgeS \Frond \cdot \nudK.
$$
By Proposition~\ref{DualityDiff}, by definitions of
$\Aleft\cdot,\cdot\Aright$, $\Bleft\cdot,\cdot\Bright$ and using
the notation introduced in Section~\ref{EmbedInLebesgue}, we have
\begin{align*}
    \dsp 0 & =\sum_{n=1}^N \Delt\Aleft \Frond^{\ptTau,n}\;,\;
    \grad^\ptTau w^{\ptbarTau,n} \Aright + \sum_{n=1}^{N} \Delt\Bleft
    \div^\ptTau [\Frond^{\ptTau,n}]\;,\; w^{\ptTau,n} \Bright \\
    \dsp  & = \sum_{n=1}^N \Delt\Aleft \Frond^{\ptTau,n}\;,\;
    \grad^\ptTau w^{\ptbarTau,n} \Aright
    + \sum_{n=1}^{N} \Delt\Bleft (\div \Frond)^{\ptTau,n}\;,\; w^{\ptTau,n} \Bright \\
    & \qquad \dsp  +\; \sum_{n=1}^{N} \Delt\Bleft \div^\ptTau
    [\Frond^{\ptTau,n}] - (\div \Frond)^{\ptTau,n}\;,\; w^{\ptTau,n}
    \Bright
    \\ & = \int_Q \Frond^{\ptTau,\Delt} \cdot \grad^\ptTau
    w^{\ptbarTau,\Delt}
    \dsp + \int_Q \Bigl(\div \Frond\Bigr)\;\Bigl(
    \frac{1}{d}\,w^{\ptMrond,\Delt}+\frac{d\!-\!1}{d}\,w^{\ptdMrond,\Delt}\Bigr) \\
    & \qquad \dsp  +\; \sum_{n=1}^{N} \Delt\Bleft \div^\ptTau
    [\Frond^{\ptTau,n}] - (\div \Frond)^{\ptTau,n}\;,\; w^{\ptTau,n}\Bright.
\end{align*}

As in Proposition~\ref{PropConsistency}, one shows that
$\|\Frond^{\ptTau,n}-\Frond\|_{L^{p'}}$ tends to zero as $h\to 0$.
Therefore we deduce
\begin{equation}\label{SecondIdent}
0= \int_Q \Frond \cdot \chi + \int_Q \Bigl(\div \Frond\Bigr)\;w
 \;+\; \lim_{h\to0}\sum_{n=1}^{N} \Delt\Bleft \div^\ptTau
[\Frond^{\ptTau,n}] - (\div \Frond)^{\ptTau,n}\;,\;
w^{\ptTau,n}\Bright.
\end{equation}

By definition of $\Bleft\cdot,\cdot\Bright$, we have
$$
\begin{array}{l}
\dsp
 \sum_{n=1}^{N} \Delt\Bleft \div^\ptTau [\Frond^{\ptTau,n}] -
(\div \Frond)^{\ptTau,n}\;,\; w^{\ptTau,n}\Bright\\[10pt]
\dsp =\frac{1}{d} \sum_{n=1}^N \Delt \Biggl(
 \mK\, w^n_\ptK\;
\frac{1}{\Delt\mK}\int_{(n\!-\!1)\Delt}^{n\Delt}
\sum\nolimits_{\ptS\in\ptVrond(\ptK)} \Bigl(\int_\edgeS \Frond -
m_\edgeS\Frond_\ptS^n) \cdot \nuK\Biggr)\\[10pt]
\dsp
 +\frac{d\!-\!1}{d} \sum_{n=1}^N \Delt
\Biggl(  \mdK\,w^n_\ptdK\;
\frac{1}{\Delt\mdK}\int_{(n\!-\!1)\Delt}^{n\Delt}
\sum\nolimits_{\ptS\in\ptdVrond(\ptdK)} \Bigl(\int_\dedgeS \Frond
- m_\dedgeS\Frond_\ptS^n) \cdot \nudK \Biggr).
  \end{array}
$$
Denote by $R+R^*$ the right-hand side above. Summing by parts, we
get
$$
R=\frac{1}{d} \sum_{n=1}^N \sum_{\ptKIL\in\ptErond}
 \dKL
\int_{(n\!-\!1)\Delt}^{n\Delt}\int_\ptKIL \Biggl(\Frond -
\frac{1}{\mD}\int_\ptD \Frond\Biggr) \cdot
\nuKL\;\frac{w^n_\ptK-w^n_\ptL}{\dKL}
$$
where $\DM$ stands for the diamond $\DM^{\ptK\ptL}$ containing the
interface $\KIL$. By the H\"older inequality, we deduce that $|R|$
is controlled by
\begin{align*}
	&\Biggl(\sum_{n=1}^N \sum_{\ptKIL\in\ptErond} \dKL\!\!
	\int_{(n\!-\!1)\Delt}^{n\Delt}\int_\ptKIL \Biggl|\Frond -
	\frac{1}{\mD}\int_\ptD \Frond\Biggr|^{p'}
	\Biggr)^{\frac1{p'}}
	\\ & \qquad \times \Biggl(\sum_{n=1}^N
	\Delt\!\!\sum_{\ptKIL\in\ptErond}
	\mKIL\dKL\;\Bigl|\frac{w^n_\ptK-w^n_\ptL}{\dKL}
	\Bigr|^p\Biggr)^{\frac1p}.
\end{align*}
Using standard estimates similar to \cite[Lemma 3.2]{ABH``double''} and
the definition of $\grad^\ptTau w^{\ptbarTau,\Delt}$, we conclude that
$$
\begin{array}{l}
\dsp |R|\leq C(\reg(\Tau))\times \size(\Tau) \times
\|\Frond\|_{L^{p'}(W^{1,{p'}})}\|\grad^\ptTau
w^{\ptbarTau,\Delt}\|_{L^p}\\[5pt]
\dsp\qquad\qquad \leq C(\reg(\Tau))\times
h\times\|\Frond\|_{L^{p'}(W^{1,{p'}})}\times C \to 0
\end{array}
$$
as $h\to 0$. In the same way, we find $|R^*|\to 0$ as $h\to 0$.

Thus for all $\Frond\in (L^{p'}(0,T;W^{1,p'}(\Om))^d$, the last
term in \eqref{SecondIdent} is zero, so that  $w\in
L^p(0,T;W^{1,p}_0(\Om))$ and $\chi=\grad w$.\\[5pt]
(ii)~~~ If also $\dsp
 \sum_{n=1}^N \Delt\Bleft \Prond^{\ptTau}[w^{\ptbarTau,n}], w^{\ptTau,n}
\Bright \leq C$, then by Lemma~\ref{SummationPenalize} with
$\psi^\ptbarTau=w^\ptbarTau$ we get
$$  \frac{d-1}{d} \sum_{n=1}^N \Delt
\sum_{\ptK\in\overline{\ptMrond}, \ptdK\in\overline{\ptdMrond}}
m_{\ptK\cap\ptdK}(w^n_\ptK-w^n_\ptdK)^2 \leq Ch.
$$
This means that
$\|w^{\ptMrond,\Delt}-w^{\ptdMrond,\Delt}\|_{L^2}\to 0$ as $h\to
0$, which permits to identify the weak limits of both
$w^{\ptMrond,\Delt}$ and $w^{\ptdMrond,\Delt}$ with $w$.
\end{proof}

\section{Properties of discrete solutions}\label{sec:DiscrSol}

\subsection{A priori estimates}

\begin{prop}\label{APrioriEstimates}
Assume we are given a family of double meshes $\Tau$ of $\Om$
and associated time steps $\Delt$ such that
$h=\max\{\size(\Tau),\Delt\} \to 0$.
Assume that $\reg(\Tau)$ is uniformly bounded.

Let $u^{\ptbarTau,\barDelt}$ be a solution to
\eqref{AbstractScheme}, \eqref{DiscrBC}, \eqref{DiscrIC} (recall
that $w^{\ptbarTau,\barDelt}=A(u^{\ptbarTau,\barDelt})$). Then the
following {\it a priori} estimates hold uniformly in $h$:\\[5pt]

(i) $\dsp\max\{\|u^{\ptMrond,\Delt}\|_{L^\infty},\|u^{\ptdMrond,\Delt}\|_{L^\infty}
\}\leq M:=\|u_0\|_{L^\infty}+\int_0^T \|{\EuScript S}(t,\cdot)\|_{L^\infty}\,dt$;\\[5pt]

(ii) there exists $C>0$ such that
$$
\|\grad^\ptTau w^{\ptbarTau,\Delt}\|_{L^p}\leq C
\quad \text{and}\quad
\dsp \sum_{n=1}^N \Delt
\Bleft \Prond^{\ptTau}[w^{\ptbarTau,n}], w^{ \ptTau,n}\Bright
\leq C;
$$

(iii) there exists $C>0$ such that (with the
notation of Proposition \ref{DualityConv})
$$
\sum_{n=1}^N
\Delt\;\Bigl(I_{\mathrm{Id}} [u^{\ptMrond,n},1]
+I^*_{\mathrm{Id}}[u^{\ptdMrond,n},1]\Bigr)\leq C;
$$

(iv) there exists a modulus of continuity
$\omega_A(\cdot)$ such that for all $\Del>0$,
$$
\dsp \int_Q
|w^{\ptMrond,\Delt}(t+\Del,x)-w^{\ptMrond,\Delt}(t,x)| +
|w^{\ptdMrond,\Delt}(t+\Del,x)-w^{\ptdMrond,\Delt}(t,x)|
\leq \omega_A(\Del),
$$
where $w^{\ptMrond,\Delt}$,$w^{\ptdMrond,\Delt}$
are extended by zero on $(N\Delt,+\infty)\times\Om$.
\end{prop}

\begin{proof}
(i) Denote ${\EuScript S}^i=(\mathbb{S}^\Delt[{\EuScript S}])^i$
and ${\EuScript
S}^{\ptTau,i}=\mathbb{P}^\ptTau[(\mathbb{S}^\Delt[{\EuScript
S}])^i]$. For $n=0,\dots,N$, set $c^n=\|u_0\|_{L^\infty}+
\sum_{i=1}^n \Delt \|{\EuScript S}^i\|_{L^\infty}$; note that
$c^n\leq \|u_0\|_{L^\infty}+ \int_0^T \|{\EuScript
S}(t,\cdot)\|_{L^\infty}\,dt=M$ for all $n=1,\dots,N$.

Let us prove by induction that $\|u^{\ptMrond,n}\|_{L^\infty}\leq
c^n$,  $\|u^{\ptdMrond,n}\|_{L^\infty}\leq c^n$. This claim is
clear for $n=0$. Assume it holds true for $n=k-1$. Take the scalar
product $\Bleft\cdot,\cdot\Bright$ of equations
\eqref{AbstractScheme} corresponding to $n=k$ with the discrete
function $\theta(u^{\ptTau,k}):=\sign^+(u^{\ptTau,k}-c^k)$. We get
\begin{equation}\label{AprioriLinfty}
\begin{array}{l}
\dsp \Bleft\frac{u^{\ptTau,k} -
u^{\ptTau,(k-1)}}{\Delt}-{\EuScript
S}^{\ptTau,k}\;,\;\theta(u^{\ptTau,k})\Bright
 + \Bleft(\div_c {\mathfrak
f})^\ptTau[u^{\ptbarTau,k}]\;,\;\theta(u^{\ptTau,k})\Bright\\[10pt]
\dsp\qquad\qquad
 -\Bleft \div^\ptTau[{\mathfrak a}(\grad^\ptTau A(u^{\ptbarTau,k}))]\;,\;
\theta(u^{\ptTau,k})\Bright
+ \Bleft{\Prond}^\ptTau[w^{\ptbarTau,k}]\;,\;
\theta(u^{\ptTau,k})\Bright=0.
\end{array}
\end{equation}
Let us apply to the last three terms above
Proposition~\ref{DualityConv},
Proposition~\ref{ParabEntropyDissipation} and
Lemma~\ref{SummationPenalize} respectively,  with $\psi\equiv 1$.
Note that $\theta(0)=0$, so that \eqref{CarrilloChoice} holds.
 We
conclude that each of the three last terms in
\eqref{AprioriLinfty} is nonnegative. Hence
$$
\begin{array}{l}
\dsp 0\geq\dsp \Bleft\frac{u^{\ptTau,k} -
u^{\ptTau,(k-1)}_\ptK}{\Delt}-{\EuScript S}^{\ptTau,k}\;,\;\theta(u^{\ptTau,k})\Bright\\[10pt]
\dsp \quad =\Bleft\frac{(u^{\ptTau,k}-c^k) -
(u^{\ptTau,(k-1)}-c^{(k-1)})}{\Delt}+\Bigl(\|{\EuScript
S}^k\|_{L^\infty}-{\EuScript S}^{\ptTau,k}\Bigr)\;,\;
\sign^+(u^{\ptTau,k}-c^k)\Bright\\[10pt]
\dsp \quad \geq \Bleft\frac{(u^{\ptTau,k}-c^k) -
(u^{\ptTau,(k-1)}-c^{(k-1)})}{\Delt}\;,\;\sign^+(u^{\ptTau,k}-c^k)\Bright\\[10pt]
\dsp \quad \geq \Bleft
(u^{\ptTau,k}-c^k)^+-(u^{\ptTau,(k-1)}-c^{(k-1)})^+\;,\;1^\ptTau\Bright,
\end{array}
$$
where $1^\ptTau=\mathbb{P}^\ptTau[1]$. By the induction hypothesis
we deduce that $(u^{\ptTau,k}-c^k)^+\leq 0$, which proves our
claim for $n=k$.\\[5pt]

(ii) For $n=1,\dots,N$, take the scalar product
$\Bleft\cdot\;,\;\cdot\Bright$ of equations \eqref{AbstractScheme}
with the discrete function $w^{\ptTau,n}=A(u^{\ptTau,n})$.
Multiply by $\Delt$ and sum up in $n$. We get
\begin{equation}\label{EstimWithW}
\begin{array}{l}
\dsp \sum_{n=1}^N \Delt\;\Bleft \frac{u^{\ptTau,n} -
u^{\ptTau,(n-1)}}{\Delt}\;,\;A(u^{\ptTau,n})\Bright
 + \sum_{n=1}^N \Delt\;\Bleft(\div_c {\mathfrak
f})^\ptTau[u^{\ptbarTau,n}]\;,\;A(u^{\ptTau,n})\Bright\\[10pt]
\dsp  - \sum_{n=1}^N \Delt\;\Bleft \div^\ptTau[{\mathfrak
a}(\grad^\ptTau w^{\ptbarTau,n})]\;,\;w^{\ptTau,n}\Bright
+ \sum_{n=1}^N \Delt\;\Bleft {\Prond}^\ptTau[w^{\ptbarTau,n}]\;,\;
w^{\ptTau,n}\Bright\\[10pt]
\dsp = \sum_{n=1}^N
\Delt\;\Bleft\mathbb{P}^\ptTau(\mathbb{S}^\Delt[{\EuScript
S}])^n\;,\;w^{\ptTau,n}\Bright.
\end{array}
\end{equation}
Note that with $\theta(\cdot)=A(\cdot)$ and $\psi\equiv 1$,
\eqref{CarrilloChoice} holds. Applying Lemma~\ref{DualityEvol},
Proposition~\ref{DualityConv}, Proposition~\ref{DualityDiff} and
Lemma~\ref{SummationPenalize}, respectively, to the terms on the
left-hand side of \eqref{EstimWithW}, we find
\begin{equation}\label{EstimWithW-2}
\begin{array}{l}
\dsp \Bleft \eta(u^{\ptTau,N})\;,\; 1^{\ptTau} \Bright
 + \sum_{n=1}^N
\Delt\; \Bigl( I_{A} [u^{\ptMrond,n},1] +
I^*_{A}[u^{\ptdMrond,n},1]\Bigr)\\[10pt]
\dsp\qquad\quad +
\sum_{n=1}^N \Delt\; \Aleft
\mathfrak{a}(\grad^\ptTau w^{\ptbarTau,n})\;,\;\grad^\ptTau
w^{\ptbarTau,n} \Aright + \sum_{n=1}^N \Delt\;  \Bleft
{\Prond}^\ptTau[w^{\ptbarTau,n}]\;,\; w^{\ptTau,n}\Bright\\[10pt]
\dsp= \sum_{n=1}^N
\Delt\;\Bleft\mathbb{P}^\ptTau(\mathbb{S}^\Delt[{\EuScript
S}])^n\;,\;w^{\ptTau,n}\Bright +\Bleft B_A(u^{\ptTau,0})\;,\;
1^{\ptTau} \Bright,
\end{array}
\end{equation}
where $B_A(z)=\int_0^{z} A(s)\,ds$ and $I_A,I^*_A$ are defined in
Proposition~\ref{DualityConv}. The first two terms in
\eqref{EstimWithW-2} are nonnegative; the next one is lower bounded
by a constant times $\Bigl(\|\grad^\ptTau w^{\ptbarTau,\Delt}\|_{L^p}\Bigr)^p$
due to the coercivity assumption on $\mathfrak{a}$.
By H\"older's inequality, Proposition \ref{PropConsistency}(i) and
Proposition~\ref{PoincareAndTranslations}(i), the first term in
the right-hand side of \eqref{EstimWithW-2} is majorated by
$C(\reg(\Tau))\times \|f\|_{L{p'}}\times\|\grad^\ptTau
w^{\ptbarTau,\Delt}\|_{L^p}$. Finally, the last term in
\eqref{EstimWithW-2} is upper bounded by a constant times
$m_{\Om}\int_{-\|u_0\|_{L^\infty}}^{\|u_0\|_{L^\infty}}A(s)\,ds$.
Hence, (ii) follows.\\[5pt]

(iii) We proceed as in (ii), multiplying equations
\eqref{AbstractScheme} by $u^{\ptTau,n}$ instead of
$A(u^{\ptTau,n})$. As in \eqref{EstimWithW-2} above, taking
$\theta=Id$, $\psi\equiv 1$, applying
Proposition \ref{ParabEntropyDissipation} instead of
Proposition \ref{DualityDiff}, neglecting the nonnegative
terms on the left-hand side, we get
\begin{equation*}
\begin{array}{l}
\dsp \sum_{n=1}^N \Delt\;
\Bigl( I_{\mathrm{Id}} [u^{\ptMrond,n},1] +
I^*_{\mathrm{Id}}[u^{\ptdMrond,n},1]\Bigr)\\[10pt]
\dsp \hspace*{20pt} \leq \sum_{n=1}^N
\Delt\;\Bleft\mathbb{P}^\ptTau(\mathbb{S}^\Delt[{\EuScript
S}])^n\;,\;u^{\ptTau,n}\Bright +\Bleft
\frac{1}{2}(u^{\ptTau,0})^2\;,\; 1^{\ptTau} \Bright.
\end{array}
\end{equation*}
Using the $L^\infty$ estimate (i) of the present proposition
together with Proposition~\ref{PropConsistency}(i), we
finally get (iii) with the constant
$$
C= C(\reg(\Tau))\times M\times \|{\EuScript S}\|_{L^1}+
\frac{1}{2}m_{\Om}\times(\|u_0\|_{L^\infty})^2.
$$
(iv) We adapt to the discrete framework the calculation that led
to estimate \eqref{Transl-base} in the proof of
Theorem~\ref{th:EPS-ESexists}(ii).
Denote by $ J(\Del), J^*(\Del)$, respectively, the integrals
\begin{equation*}
\begin{array}{l}
\dsp  \int_Q
|u^{\ptMrond,\Delt}(t\!+\!\Del,x)\!-\!u^{\ptMrond,\Delt}(t,x)|\;
|A(u^{\ptMrond,\Delt})(t\!+\!\Del,x)\!-\!A(u^{\ptMrond,\Delt})(t,x)|,\\[10pt]
\dsp \int_Q
|u^{\ptdMrond,\Delt}(t\!+\!\Del,x)\!-\!u^{\ptdMrond,\Delt}(t,x)|\;
|A(u^{\ptdMrond,\Delt})(t\!+\!\Del,x)\!-\!A(u^{\ptdMrond,\Delt})(t,x)|.
\end{array}
\end{equation*}
Let us first take $k\in\{1,\dots,N\}$ and estimate the quantity
\begin{equation*}
J_0(k):=\dsp\sum_{n=k+1}^{N} \Delt \Bleft
u^{\ptTau,n}-u^{\ptTau,(n-k)}\;,\;
A(u^{\ptTau,n})-A(u^{\ptTau,(n-k)}) \Bright.
\end{equation*}
To do this, for $n=(k+1),\dots,N$ we take the sum in $i$ from $(n-k+1)$ to
$n$ of equations \eqref{AbstractScheme} and make the scalar product
$\Bleft\cdot\;,\;\cdot\Bright$ with the discrete functions
$v^{\ptTau,n}$, where
$v^{\ptbarTau,n}:=A(u^{\ptbarTau,n})-A(u^{\ptbarTau,(n-k)})\in\R^\ptbarTau_0$
for $n=(k+1),\dots,N$. Summing in $n$ and assigning $v^{\ptbarTau,n}=0$
for $n=1,\dots,k$ and $n=(N+1),\dots,(N+k-1)$, we get
\begin{equation}\label{EstimOfTimeTransl}
\begin{array}{l}
\dsp \frac{J_0(k)}{\Delt}
=\sum_{n=k+1}^{N} \Bleft u^{\ptTau,n}-u^{\ptTau,(n-k)}\;,\; v^{\ptTau,n} \Bright \\[10pt]
\dsp \quad\qquad =\sum_{n=k+1}^{N} \Delt\Bleft \sum_{i=n-k+1}^{n}
\frac{u^{\ptTau,i}-u^{\ptTau,(i-1)}}{\Delt} \;,\; v^{\ptTau,n}
\Bright \\[10pt]
\dsp \quad\qquad
=\sum_{i=2}^N\;\sum_{j=\max\{1,k-i+2\}}^{\min\{k,N-i+1\}}\Delt\Bleft
 \frac{u^{\ptTau,i}-u^{\ptTau,(i-1)}}{\Delt}
\;,\; v^{\ptTau,(i+j-1)} \Bright
\\[10pt]
\dsp \quad\qquad =\sum_{j=1}^{k}\sum_{i=2}^N
\Delt \Bleft
-(\div_c {\mathfrak f})^\ptTau[u^{\ptbarTau,i}]
+\div^\ptTau[{\mathfrak a}(\grad^\ptTau w^{\ptbarTau,i})]\\[10pt]
\dsp\hspace*{110pt} -{\Prond}^\ptTau[w^{\ptbarTau,i}]
+\mathbb{P}^\ptTau(\mathbb{S}^\Delt[{\EuScript S}])^i\; , \;
v^{\ptTau,(i+j-1)}\Bright.
\end{array}
\end{equation}
We claim that the right-hand side of \eqref{EstimOfTimeTransl} is
bounded by a constant independent of $h$. Indeed, for each
$j=1,\dots,k$, define ${z_j}^{\ptbarTau,i}=v^{\ptbarTau,(i+j-1)}$, $i=1,\dots,N$.
First, from the property (ii) of the present
proposition and from formula \eqref{SummPenalWeak} we deduce
\begin{equation}\label{ControlOfA-A}
\|\grad^\ptTau z_j^{\ptbarTau,\Delt}\|_{L^p}\leq C,
\quad \sum_{i=1}^N \Delt\Bleft \Prond^{\ptbarTau}[z_j^{\ptTau,i}], z_j^{ \ptTau,i}
\Bright \leq C,
\quad \text{for all $j=1,\dots,k$}.
\end{equation}
In the sequel, we will omit the dependency of the entries of
$z_j^{\ptbarTau,\Delt}$ on $j$.

By definition of $(\div_c {\mathfrak f})^\ptTau[\cdot]$, taking
into account that $z_j^{\ptbarTau,n}\in\R^\ptbarTau_0$ and using
summation-by-parts, we deduce that for all $j=1,\dots,k$,
\begin{equation*}
\begin{split}
\dsp J_{1,j} & :=\Biggl|\sum_{i=1}^{N} \Delt \Bleft (\div_c
{\mathfrak f})^\ptTau[u^{\ptbarTau,i}] \;,\; z_j^{\ptTau,i}
\Bright\Biggr|\\
& =\Biggl| \sum_{i=1}^{N} \Delt
\Biggl(\frac{1}{d}\;\sum_{\ptKIL\in\ptErond} \mKIL
g_\ptKL(u_\ptK^i,u_\ptL^i)\;(z^i_\ptK\!-\!z^i_\ptL)\\
& \hspace*{30mm} +\frac{d\!-\!1}{d}\!\!\!\sum_{\ptdKIdL\in\ptdErond}
\!\! \mdKIdL
g_\ptdKdL(u_\ptdK^i,u_\ptdL^i)\;(z^i_\ptdK\!-\!z^i_\ptdL)\Biggr)
\Biggr|.
\end{split}
\end{equation*}
Since by (i), $u^{\ptMrond,\Delt},u^{\ptdMrond,\Delt}$ are bounded
by $M$, using property \eqref{Hypfluxes}(d) we bound all
values of $g_\ptKIL,g_\ptdKIdL$ above by $C\omega_M(M)$. It
follows by Remark \ref{ExplicDiscrGrad} that
$$
\frac{|z^i_\ptK\!-\!z^i_\ptL|}{\dKL} +
\frac{|z^i_\ptdK\!-\!z^i_\ptdL|}{\ddKdL}
\leq |\grad_{\!\ptS} z^{\ptbarTau,i}|,
$$
where $\SDM=\SDM^{\ptK,\ptL}_{\ptdK,\ptdL}$. Hence
\begin{equation*}
\dsp J_{1,j} \leq C(d-1)\omega_M(M) \sum_{i=1}^{N} \Delt
\sum_{\ptS\in\ptSrond} \mS |\grad_{\!\ptS} z^{\ptbarTau,i}|
\le
\mathrm{const}\, \|\grad^{\ptTau} z_j^{\ptbarTau,\Delt}\|_{L^1}.
\end{equation*}
Using \eqref{ControlOfA-A}, we can uniformly bound $J_{1,j}$.
Further, by Proposition~\ref{DualityDiff} and the H\"older inequality,
\begin{equation*}
\begin{split}
\dsp J_{2,j}& :=\Biggl| \sum_{i=1}^{N} \Delt
\Bleft\div^\ptTau[{\mathfrak a}(\grad^\ptTau
w^{\ptbarTau,i})]\;,\;z_j^{\ptTau,i} \Bright \Biggr|\\
& =\Biggl| \sum_{i=1}^{N} \Delt \Aleft {\mathfrak a}(\grad^\ptTau
w^{\ptbarTau,i})\;,\; \grad^\ptTau z_j^{\ptbarTau,i}\Aright
\Biggr|\leq \| {\mathfrak a}(\grad^\ptTau
w^{\ptbarTau,\Delt})\|_{L^{p'}} \| \grad^\ptTau
z_j^{\ptbarTau,\Delt}\|_{L^p}.
\end{split}
\end{equation*}
Using the growth assumption on $\mathfrak{a}$ together with
\eqref{ControlOfA-A} and (ii) of the present lemma, we can uniformly bound $J_{2,j}$.
Next, by \eqref{SummPenalWeak} and the Cauchy-Schwarz inequality,
\begin{equation*}
\begin{split}
    J_{3,j} & :=\Biggl| \sum_{i=1}^{N} \Delt \Bleft
    {\Prond}^\ptTau[w^{\ptbarTau,i}]  \;,\; z_j^{\ptTau,i} \Bright \Biggr| \\
    & \leq \frac{d\!-\!1}{d}\sum_{i=1}^{N} \Delt \!\!\!\!\!\!\!
    \sum_{\ptK\in\overline{\ptMrond}, \ptdK\in\overline{\ptdMrond}}
    \!\!\!\!\!\!\!  m_{\ptK\cap\ptdK}
    \frac{|w^i_\ptK-w^i_\ptdK|}{\sqrt{\size(\Tau)}}\frac{|z^i_\ptK-z^i_\ptdK|}{\sqrt{\size(\Tau)}}\\
    & \leq \frac{d\!-\!1}{d} \left(\sum_{i=1}^{N} \Delt
    \!\!\!\!\!\!\! \sum_{\ptK\in\overline{\ptMrond},
    \ptdK\in\overline{\ptdMrond}}\!\!\!\!\!\!\!\!\! m_{\ptK\cap\ptdK}
    \frac{|w^i_\ptK-w^i_\ptdK|^2}{\size(\Tau)} \right)^{1/2}\\
    & \hspace*{25mm} \times \left(\sum_{i=1}^{N} \Delt
    \!\!\!\!\!\!\!\sum_{\ptK\in\overline{\ptMrond},
    \ptdK\in\overline{\ptdMrond}} \!\!\!\!\!\!\!\!\! m_{\ptK\cap\ptdK}
    \frac{|z^i_\ptK-z^i_\ptdK|^2}{\size(\Tau)}\right)^{1/2}\\
    & = \left(\sum_{i=1}^{N} \Delt \Bleft
    \Prond^{\ptTau}[w^{\ptbarTau,i}], w^{\ptTau,i} \Bright
    \right)^{1/2} \left(\sum_{i=1}^{N} \Delt \Bleft
    \Prond^{\ptTau}[z_j^{\ptbarTau,i}], z_j^{ \ptTau,i} \Bright
    \right)^{1/2}.
\end{split}
\end{equation*}
Using again \eqref{ControlOfA-A} and (ii) of the present
proposition, we can uniformly bound $J_{3,j}$. Finally,  like in
\eqref{EstimWithW-2}, we have
$$
J_{4,j}:=\Biggl| \sum_{i=1}^{N} \Delt \Bleft
\mathbb{P}^\ptTau(\mathbb{S}^\Delt[{\EuScript S}])^i \;,\;
z_j^{\ptTau,i} \Bright\;\Biggr|\leq C(\reg(\Tau))\times
\|{\EuScript S}\|_{L^{p'}}\times\|\grad^\ptTau
z_j^{\ptbarTau,\Delt}\|_{L^p},
$$
which is also uniformly bounded, thanks to \eqref{ControlOfA-A}.
Gathering the estimates above, we conclude
$$
J_0(k)\leq \Delt\;\sum_{j=1}^k
\Biggl(J_{1,j}+J_{2,j}+J_{3,j}+J_{4,j} \Biggr)\leq C\;k\Delt.
$$
Using the definition of $\Bleft\cdot\,\cdot\Bright$ and the
$L^\infty$ estimate on $u^{\ptMrond,\Delt}$, cf.~(i), we get
\begin{equation}\label{EstJandJ*}
\begin{array}{l}
\dsp \frac{1}{d}J(k\Delt)+\frac{d\!-\!1}{d}J^*(k\Delt)\leq J_0(k)\!
+ \!\!\int^{N\Delt}_{(N-k)\Delt} \!\!\!\!\!\!\!\!\! m_\Om M
\max\{\pm A(\pm M)\} \leq C\;k\Delt.
\end{array}
\end{equation}
Now let $0<\Del<T$. We have $\Del/\Delt=(k-1)+\alpha$ for
some $k\in\{1,\dots,N\}$ and $\alpha\in[0,1)$.
Since $u^{\ptMrond,\Delt}$ is piecewise constant in $t$ with step $\Delt$, we have
\begin{equation}\label{EstJDelt}
\begin{split}
    J(\Del) & =J((k-1)\Delt+\alpha\Delt)\leq \alpha
    J(k\Delt)+(1-\alpha)J((k-1)\Delt)
    \\ & \leq \alpha
    C\,k\Delt+(1-\alpha)C\,(k-1)\Delt
    \leq C\,((k-1)+\alpha)\Delt=C\,\Del.
\end{split}
\end{equation}
From \eqref{EstJDelt}, together with the calculation used to pass from
\eqref{Transl-base} to \eqref{TranslEst-truc3} (cf.~the proof
of Theorem~\ref{th:EPS-ESexists}), we deduce the required estimate
$$
\int_Q
|A(u^{\ptMrond,\Delt})(t+\Del,x)-A(u^{\ptMrond,\Delt})(t,x)|\leq
\omega_A(\Del), \qquad \Del>0.
$$
Similarly, time translates of $A(u^{\ptdMrond,\Delt})$
are controlled with $J^*(k\Delt)$ in \eqref{EstJandJ*}.
\end{proof}

\subsection{Existence of discrete solutions}

\begin{prop}\label{PropDiscrExist}
Let $\Tau$ be a double mesh of $\Om$ and $\Delt>0$. There exists a
solution $u^{\ptbarTau,\barDelt}$ of the finite volume scheme
\eqref{AbstractScheme}, \eqref{DiscrBC}, \eqref{DiscrIC}.
\end{prop}

\begin{proof}
First note that it is sufficient to prove existence of solutions
$u_\rho^{\ptbarTau,\barDelt}$ to \eqref{AbstractScheme},
\eqref{DiscrBC}, \eqref{DiscrIC} with $A(\cdot)$ replaced by a
strictly increasing function $A_\rho(\cdot)$. Indeed, using the
$L^\infty$ estimate (i) of Proposition~\ref{APrioriEstimates},
which is independent of the choice of $A(\cdot)$, we get
compactness of $u_\rho^{\ptbarTau,\barDelt}$ in the
finite-dimensional space $\R^{(N+1)\times\ptbarTau}$. Choosing a
sequence of strictly increasing functions $A_\rho$ that converges
to $A$ uniformly on all compact of $\R$, we pass to the limit in
the scheme \eqref{AbstractScheme}, \eqref{DiscrBC},
\eqref{DiscrIC} written for (a subsequence of) $A_\rho(\cdot)$ and
$u_\rho^{\ptbarTau,\barDelt}$ and obtain existence for general
$A(\cdot)$.

Let us now assume that $A(\cdot)$ is invertible
and rewrite the scheme in terms of $w^{\ptbarTau,\barDelt}$ with
$u^{\ptbarTau,\barDelt}=A^{-1}(w^{\ptbarTau,\barDelt} )$. The
existence of $w^{\ptTau,n}$ is shown by induction on $n=0,\dots,N$. For
$n=0$, solution is given by \eqref{DiscrIC}. Assume that
$w^{\ptTau,(n-1)}$ exists. Choose $\Bleft\cdot\,,\,\cdot\Bright$
as  the scalar product on $\R^\ptTau$. We are looking for a
solution $w^{\ptTau,n}$ to $L[ w^{\ptTau,n}]=0$, where the
operator $L$ is given by
$$
\begin{array}{l}
\dsp L:z^{\ptTau}\in\R^\ptTau \mapsto \frac{A^{-1}(z^{\ptTau}) -
A^{-1}(w^{\ptTau,(n-1)})}{\Delt} + (\div_c {\mathfrak
f})^\ptTau[A^{-1}(z^{\ptbarTau})]\\[7pt]
\dsp \qquad\qquad\qquad\qquad\qquad -\div^\ptTau[{\mathfrak
a}(\grad^\ptTau z^{\ptbarTau})] +
{\Prond}^\ptTau[z^{\ptbarTau}]-\mathbb{P}^\ptTau(\mathbb{S}^\Delt[{\EuScript
S}])^n.
\end{array}
$$
By Proposition \ref{DualityConv} with $\theta=\mathrm{Id}$ and
$\psi\equiv 1$, by Proposition \ref{DualityDiff} and by
Lemma~\ref{SummationPenalize}, there exists a constant
$C=C\Bigl(\|w^{\ptTau,{n-1}}
\|_{\R^\ptTau},\|\mathbb{P}^\ptTau(\mathbb{S}^\Delt[{\EuScript
S}])^n\|_{\R^\ptTau},\Delt\Bigr)$ such that
$$
\Bleft L[z^{\ptTau}] \,,\,z^{\ptTau}\Bright\geq \Aleft {\mathfrak
a}(\grad^\ptTau z^{\ptbarTau})\,,\,\grad^\ptTau
z^{\ptbarTau}\Aright\; \;-C\;\|z^{\ptTau}\|_{\R^\ptTau}.
$$
By the coercivity assumption on $\mathfrak{a}$ and by
Proposition~\ref{PoincareAndTranslations}(i) we have
\begin{equation}\label{Coercivity-2}
\Aleft {\mathfrak a}(\grad^\ptTau z^{\ptbarTau})\,,\,\grad^\ptTau
z^{\ptbarTau}\Aright \geq \mathrm{const}\; \|\grad^\ptTau
z^{\ptbarTau}\|^p_{L^p}\geq
\mathrm{const} \;\Bigl(\|z^\ptMrond\|^p_{L^p}+
\|z^\ptdMrond\|^p_{L^p}\Bigr).
\end{equation}
Because the right-hand side of \eqref{Coercivity-2} is equivalent
to $\Bigl(\|z^\ptTau\|_{\R^\ptTau}\Bigr)^p$, we conclude that
$\Bleft L[z^{\ptTau}] \,,\,z^{\ptTau}\Bright\geq 0$ for
$\|z^{\ptTau}\|_{\R^\ptTau}$ sufficiently large. The existence of
$w^{\ptTau,n}$ follows by the standard Brouwer fixed point
argument (see \cite[Lemme 4.3]{Lions:Book69_Fr}).
\end{proof}

We point out that the uniqueness and, more generally, continuous
dependency of the discrete solutions
on the data can be established as well 
(see \cite{EyGaHe:book,EGHMichel,ABH``double''} for results of that sort).
However, in view of the convergence result of Theorem \ref{ThFVConverge}
and the well-posedness of the continuous problem, we view these
questions to be of less importance.

\subsection{Discrete entropy inequalities}

\begin{prop}\label{PropApproxEntrIneq}
Let $\Tau$ be a double mesh of $\Om$ and $\Delt>0$. Consider a
solution $u^{\ptbarTau,\barDelt}$ to the scheme
\eqref{AbstractScheme},\eqref{DiscrBC},\eqref{DiscrIC}; recall
that $w^{\ptbarTau,\barDelt}=A(u^{\ptbarTau,\barDelt})$.

Let $\psi\in {\mathcal D}(\overline{Q})$, $\psi\geq 0$; set
$\psi^{\ptbarTau,\Delt}=\mathbb{P}^\ptbarTau\circ
\mathbb{S}^\Delt[\psi]$. Let $\theta:\R\to\R$ be a nondecreasing
function; assume that $\psi$ and $\theta$ are chosen so that
\eqref{CarrilloChoice} holds; assume that  $\Delt$ is small
enough. Then
\begin{equation}\label{ApproxEntrIneq}
\begin{array}{l}
\dsp  -\Bleft \eta(u^{\ptTau,N})\;,\; \psi^{\ptTau,N}
\Bright\;+\;\sum_{n=1}^{N-1} \Delt \Bleft \eta(u^{\ptTau,n})\;,\;
\frac{\psi^{\ptTau,(n+1)} - \psi^{\ptTau,n}}{\Delt}\Bright\\[10pt]
\dsp + \!\sum_{n=1}^N \!\Delt\Bleft {\mathfrak
q}(u^{\ptTau,n})\,,\,(\grad \psi)^{\ptTau,n} \Bright
 - \! \sum_{n=1}^N \!\Delt \Aleft k\Bigl(\grad^\ptTau \!
w^{\ptbarTau,n}\Bigr)\grad^\ptTau \!\widetilde
A_\theta(w^{\ptbarTau,n}) \, ,\, \grad^\ptTau\psi^{\ptbarTau,n}\!
\Aright\!\!\!\!\\[10pt]
\dsp +\Bleft \eta(u^{\ptTau,0})\;,\; \psi^{\ptTau,1} \Bright
+\sum_{n=1}^N \Delt
\Bleft\mathbb{P}^\ptTau(\mathbb{S}^\Delt[f])^n\,,\,\theta(u^{\ptTau,n})
\psi^{\ptTau,n}\Bright\\[10pt]
\dsp \qquad \geq \frac{d-1}{d} \sum_{n=1}^N
\Delt\!\!\!\!\!\!\!\sum_{\ptK\in\overline{\ptMrond},
\ptdK\in\overline{\ptdMrond}}\!\!\!\!\!\!\!\!\!
m_{\ptK\cap\ptdK}\;\theta(u_\ptK)\,\frac{(w^n_\ptK-w^n_\ptdK)(\psi^n_\ptK
-\psi^n_\ptdK)}{\size(\Tau)}\\[10pt]
\hspace*{50mm}\dsp
+R_\theta[u^{\ptMrond,n},\psi^n]+R^{*}_\theta[u^{\ptdMrond,n},\psi^n],
\end{array}
\end{equation}
where $A_\theta(\cdot),\widetilde A_\theta(\cdot)$ and
$\eta(\cdot),\mathfrak{q}(\cdot)$, $R_\theta[\cdot,\cdot],R_\theta^*[\cdot,\cdot]$
are introduced in Definition \ref{AThetaDef}
and in Proposition~\ref{DualityConv}, respectively.

Moreover, with the specific choice $\theta\equiv 1$ and
$\psi\in{\mathcal D}([0,T)\times\Om)$, there holds
\begin{equation}\label{ApproxSolFaible}
\begin{array}{l}
\dsp  \sum_{n=1}^{N-1} \Delt \Bleft u^{\ptTau,n}\!,\,
\frac{\psi^{\ptTau,(n+1)} \!\!-\! \psi^{\ptTau,n}}{\Delt}\Bright
+\Bleft u^{\ptTau,0}\!,\, \psi^{\ptTau,1} \Bright
+\! \sum_{n=1}^N
\Delt\Bleft {\mathfrak f}(u^{\ptTau,n})\,,\,(\grad
\psi)^{\ptTau,n} \Bright\!\!\!\!\\[10pt]
\dsp  - \sum_{n=1}^N \Delt \Aleft k\Bigl(\grad^\ptTau
w^{\ptbarTau,n}\Bigr)\,\grad^\ptTau w^{\ptbarTau,n}\! ,\,
\grad^\ptTau\psi^{\ptbarTau,n} \Aright  +\sum_{n=1}^N \Delt
\Bleft\mathbb{P}^\ptTau(\mathbb{S}^\Delt[{\EuScript S}])^n\,,\,\psi^{\ptTau,n}\Bright\!\!\!\!\\[10pt]
\dsp \quad = \frac{d-1}{d} \sum_{n=1}^N
\Delt\!\!\!\!\!\!\!\sum_{\ptK\in\overline{\ptMrond},
\ptdK\in\overline{\ptdMrond}}\!\!\!\!\!\!\!\!\!
m_{\ptK\cap\ptdK}\;\frac{(w^n_\ptK-w^n_\ptdK)(\psi^n_\ptK-\psi^n_\ptdK)}{\size(\Tau)}\\[10pt]
\hspace*{50mm}\dsp +
R_1[u^{\ptMrond,n},\psi^n]+R^{*}_1[u^{\ptdMrond,n},\psi^n].
\end{array}
\end{equation}

Finally, with the specific choices $\theta\equiv A$ and
$\psi\equiv \zeta(t)$, where $\zeta\in \mathcal D([0,T))$ is a
nonnegative, nonincreasing function with $\zeta(t)\equiv 1$ for
small $t$, we have with $B(z)=\int_0^z A(s)\ds$
\begin{equation}\label{ApproxTestFunctionA}
    \begin{split}
        & \sum_{n=1}^{N-1} \Delt \Bleft B(u^{\ptTau,n})\;,\;
        \frac{\zeta^{\ptTau,(n+1)} - \zeta^{\ptTau,n}}{\Delt}\Bright
        +\Bleft B(u^{\ptTau,0})\;,\;   1^{\ptTau} \Bright
        \\ & \quad+\sum_{n=1}^N \Delt
        \Bleft\mathbb{P}^\ptTau(\mathbb{S}^\Delt[{\EuScript S}])^n\,,\,w^{\ptTau,n}\zeta^{\ptTau,n}\Bright
        \\ & \quad\quad
        \geq \;  \sum_{n=1}^N \Delt \Aleft k\Bigl(\grad^\ptTau
        w^{\ptbarTau,n}\Bigr)\,\grad^\ptTau w^{\ptbarTau,n} \, ,\,
        \grad^\ptTau w^{\ptbarTau,n}\, \zeta^{\ptbarTau,n} \Aright .
    \end{split}
\end{equation}
\end{prop}

\begin{proof}
Inequality \eqref{ApproxEntrIneq} follows by an application of
Lemma \ref{DualityEvol}, Proposition \ref{DualityConv},
Proposition \ref{ParabEntropyDissipation} and Lemma
\ref{SummationPenalize}. Note that in \eqref{ApproxEntrIneq}, we
have neglected the positive terms
$I_\theta[u^{\ptMrond,n},\psi^n]$,
$I^{*}_\theta[u^{\ptdMrond,n},\psi^n]$. In \eqref{ApproxSolFaible}
the corresponding terms are zero because $\theta\equiv 1$, and we
use the equality of Proposition~\ref{DualityDiff} instead of the
inequality of Proposition~\ref{ParabEntropyDissipation}. Also
notice that the term with $\psi^{\ptTau,N}$ in
Lemma~\ref{DualityEvol} disappears because $\Delt$ is small and
$\psi$ vanishes in a neighborhood of $t=T$. Finally, in
\eqref{ApproxTestFunctionA} we have treated
$A(u^{\ptbarTau,\Delt})\,\zeta^{\ptbarTau,\Delt}$ as a mere test
function by applying Proposition \ref{DualityDiff} on the
right-hand side, but we have used Lemma~\ref{DualityEvol},
Proposition \ref{DualityConv} and the choice of the constant in
$x$ function $\psi^{\ptbarTau,\Delt}$ to deal with the remaining
terms.
\end{proof}

\subsection{Control of the remainder terms in Proposition \ref{PropApproxEntrIneq}}\label{subsec:PenalVanish}
For all $\psi\in{\mathcal D}(\overline{Q})$, the terms on the
right-hand side of \eqref{ApproxEntrIneq},\eqref{ApproxSolFaible}
coming from the penalization operator vanish as $h\to 0$. 
Indeed, using the estimates of
Proposition \ref{APrioriEstimates}(i),(ii), the Cauchy-Schwarz
inequality, Proposition \ref{PropConsistency}(iv), and
the boundedness of $\theta$ on $[-M,M]$, we obtain
$$
\begin{array}{l}
\dsp \biggl|\sum_{n=1}^N \Delt\sum_{\ptK\in\overline{\ptMrond},
\ptdK\in\overline{\ptdMrond}}
m_{\ptK\cap\ptdK}\;\theta(u_\ptK)\,\frac{(w^n_\ptK
-w^n_\ptdK)(\psi^n_\ptK-\psi^n_\ptdK)}{\size(\Tau)}\biggr|\\[10pt]
\quad \dsp \leq C\, \left( \sum_{n=1}^{N} \Delt \Bleft
\Prond^{\ptTau}[w^{\ptbarTau,n}], w^{\ptTau,n} \Bright
\right)^{1/2} \left(\sum_{n=1}^{N} \Delt\!\!\!\!\!\!\!
\sum_{\ptK\in\overline{\ptMrond}, \ptdK\in\overline{\ptdMrond}}
\!\!\!\!\!\!\!\!\! m_{\ptK\cap\ptdK}
\frac{|\psi^n_\ptK-\psi^n_\ptdK|^2}{\size(\Tau)}
\right)^{1/2}\\[9pt]
\quad \dsp \leq C\|\grad\psi\|_{L^\infty}\times \size(\Tau).
\end{array}
$$

Let us show that the terms $R_\theta [u^\ptMrond,\psi],R^*_\theta [u^\ptdMrond,\psi]$ in
\eqref{ApproxEntrIneq},\eqref{ApproxSolFaible} (which are defined in
Proposition \ref{DualityConv}) vanish as $h\to 0$. This holds true
thanks to their upper bounds in terms the quantities
$I_{\mathrm{Id}} [u^\ptMrond,1]$, $I^*_{\mathrm{Id}}[u^\ptdMrond,1]$, quantities
which are controlled by means of
Proposition \ref{APrioriEstimates}(iii) (known as the ``weak BV
estimate'', cf.~\cite{EyGaHe:book}).

\begin{prop}\label{ControlByWeakBV}
Let $g_\ptKL\in C(\R^2)$ be a function with properties
\eqref{Hypfluxes}(a),(d). For $a,b\in \R$, consider
$$
I^\ptKIL_{\mathrm{Id}}(a,b)=\int_a^b
\bigl(g_\ptKL(s,s)-g_\ptKL(a,b)\bigr)\, ds,
$$
$$
R^\ptKIL_\ptK(a,b)=|g_\ptKL(a,a)-g_\ptKL(a,b)|, \quad
R^\ptKIL_\ptL(a,b)=|g_\ptKL(b,b)-g_\ptKL(a,b)|.
$$

There exists a continuous strictly increasing convex function
$\Pi_M:\R^+\to\R^+$ that only depends on $C$ and $\omega_M(\cdot)$
in \eqref{Hypfluxes}(d) such that $\Pi_M(0)=0$, $\Pi_M'(0)=0$ and
the following bounds hold:
\begin{equation}\label{EstimViaSignedTerm}
\left\{
\begin{array}{l}
 \dsp R^\ptKIL_\ptK(a,b) \leq \Pi_M^{-1}
(I^\ptKIL_{\mathrm{Id}}(a,b)),\\[4pt]  R^\ptKIL_\ptL(a,b)\leq
\Pi_M^{-1} (I^\ptKIL_{\mathrm{Id}}(a,b)), \end{array} \right. \qquad
\text{for all $a,b\in[-M,M]$}.
\end{equation}
\end{prop}

The proof is based upon the following generalization of
\cite[Lemma 4.5]{EyGaHe:book}.

\begin{lem}\label{omega-1Lemma}
Let $g\in C([a,b])$ be a nondecreasing function equipped with a modulus of
continuity $\omega$. Then
\begin{equation*}
\int_a^b \bigl(g(s)-g(a)\bigr)\,ds \geq \int_0^{g(b)-g(a)}
\omega^{-1}(r)\,dr.
\end{equation*}
\end{lem}
\begin{proof}
Set $\delta=\omega^{-1}(g(b)-g(a))$. Since $|g(b)-g(s)|\leq
\omega(b-s)$ and $g$ is nondecreasing, we have
$$
g(s) \geq \left\{ \begin{array}{lcl} g(b)-\omega(b-s), & &
b-\delta\leq s\leq b\\
g(a), && a\leq s\leq b-\delta.
\end{array} \right.
$$
Hence setting $z=b-s$, integrating by parts, and setting
$r=\omega(z)$, we deduce
\begin{equation*}
\begin{split}
 \dsp\int_a^b \bigl(g(s)-g(a)\bigr)\,ds  & \;\geq\; \int_{b-\delta}^b
\Bigl(g(b)-g(a)-\omega(b-s)\Bigr)\,ds\\
&=\;\delta\omega(\delta)-\int_0^\delta \omega(z)\,dz=\int_0^\delta
z\,d\omega(z)=\int_0^{\omega(\delta)}
\omega^{-1}(r)\,dr.\\[-15pt]
\end{split}
\end{equation*}
\end{proof}

\begin{proof}[Proof of Proposition \ref{ControlByWeakBV}]
Consider the case $a\leq b$. By \eqref{Hypfluxes}(a), we have
$$I^\ptKIL_{\mathrm{Id}}(a,b)=\int_a^b
\bigl(g_\ptKL(s,s)-g_\ptKL(a,b)\bigr)\, ds \geq \int_a^b
\bigl(g_\ptKL(s,b)-g_\ptKL(a,b)\bigr)\, ds;
$$
applying Lemma~\ref{omega-1Lemma} to $g(\cdot)=g_\ptKL(\cdot,b)$
and recalling \eqref{Hypfluxes}(d), we deduce
$$
I^\ptKIL_{\mathrm{Id}}(a,b)\geq \int_0^{g_\ptKL(b,b)-g_\ptKL(a,b) }
(C\omega_M)^{-1}(r)\,dr =\int_0^{R^\ptKIL_\ptL(a,b)}
(C\omega_M)^{-1}(r)\,dr.
$$
Thus in order to estimate $R^\ptKIL_\ptL(a,b)$ as in
\eqref{EstimViaSignedTerm}, it is sufficient to take the function
$\dsp\Pi_M:R\in\R^+\mapsto \int_0^R(C\omega_M)^{-1}(r)\,dr$.
Clearly, $\Pi_M$  is continuous, strictly increasing, convex,
$\Pi_M(0)=0$, and $\Pi_M'(0)=0$.

The other estimate in \eqref{EstimViaSignedTerm} is obtained in
the same way, and the case $a>b$ is obtained by symmetry.
\end{proof}

\begin{cor}\label{UseOfWeakBVestim}
(i) Consider $I_{\mathrm{Id}}[u^\ptMrond,1]$ defined as in
\eqref{SignedTerm},\eqref{SignedTermEntries} with $\theta=\mathrm{Id}$,
and $\psi\equiv 1$. For general nondecreasing $\theta(\cdot)$ and
general $\psi\in{\mathcal D}(\overline{\Om})$, consider
$R_\theta[u^\ptMrond,\psi]$ defined in \eqref{TermRtheta},\eqref{TermRthetaEntries}.
Assume $\|u^{\ptMrond}\|_\infty\leq M$. Let $\Pi_M$
be the function given in Proposition \ref{ControlByWeakBV}. Let
$\Pi^*_M$ be the conjugate convex function of $\Pi_M$. Then
\begin{equation}\label{HypErrorMajoration}
\begin{split}
& \Bigl| R_\theta[u^{\ptMrond},\psi] \Bigr|
\\ &  \leq
2\|\theta\|_{C([-M,M])}
\inf_{\alpha>0}
\Biggl(\frac{\size(\Tau)}{\alpha} I_{\mathrm{Id}}[u^{\ptMrond},1] +
\frac{C}{\alpha}\Pi^*_M\biggl(2\alpha
\max_{\ptK\in\ptMrond,\,\ptL\in\ptNrond(\ptK)}
\frac{|\psi_\ptK\!\!-\!\psi_\ptKIL|}{d_{\ptK,\ptKIL}} \biggr)
\Biggr),\!\!\!\!\!\!\!\!\!\!
\end{split}
\end{equation}
where $C$ depends on $\reg(\Tau)$, $d$ and $\Omega$.

(ii) Assume we are given a sequence of meshes $\Tau$ with $\size(\Tau)\to 0$
and time steps $\Delt\to 0$. Let $u^{\ptTau,\Delt}$ be the corresponding
discrete functions such that
$\|u^{\ptMrond,\Delt}\|_\infty \leq M$ and $ \sum_{n=1}^N \Delt\;
I_{\mathrm{Id}}[u^{\ptMrond,n},1]\leq C$
uniformly in $\Tau,\Delt$. Choose $\psi\in {\mathcal D}(\overline{Q})$ and
take $\psi^n=(\mathbb{S}^\Delt[\psi])^n$.
Then $\sum_{n=1}^N \Delt\;
R^*_{\theta}[u^{\ptbarTau,n},\psi^n]\to 0$ as $\size(\Tau)\to 0$.

Analogous statements that involve $ \sum_{n=1}^N \Delt\;
I^*_{\mathrm{Id}}[u^\ptdMrond,1]$ and $\psi_\ptdK$, $\psi_\ptdKIdL$ with
$\dK\in\dMrond$, $\dL\in\dNrond(\dK)$ hold for
$\sum_{n=1}^N \Delt\; R^*_{\theta}[u^\ptdMrond,\psi]$.
\end{cor}

\begin{proof}
(i) By \eqref{TermRtheta} and
Proposition \ref{ControlByWeakBV}, for all $\alpha>0$ we have
\begin{align*}
    & \Bigl|R_\theta [u^{\ptMrond},\psi]\Bigr|
    \\ & \quad
    \leq 2\|\theta\|_{C([-M,M])}
    \sum_{\ptK\in\ptMrond,\,\ptL\in\ptNrond(\ptK)}
    \left(\frac{1}{\alpha}\mKIL d_{\ptK,\ptKIL}\right)
    \times\Pi_M^{-1}(I^{\ptKIL}_{\mathrm{Id}})
    \times \left(\alpha\,\frac{|\psi_\ptK\!-\!\psi_\ptKIL|}{d_{\ptK,\ptKIL}}\right).
\end{align*}
Note that $\dKL\leq \size(\Tau)$. Further, even in the case the
diamonds are not necessarily convex, the definition of
$\reg(\Tau)$ permits to control the multiplicity of the covering
of $\Om$ by the convex envelopes of $\K$ and $\KIL$,
$\K\in\Mrond,\,\L\in\Nrond(\K)$. Thus one can upper bound
$\sum_{\ptK\in\ptMrond,\,\ptL\in\ptNrond(\ptK)} \mKIL
d_{\ptK,\ptKIL}$ by $C(\reg(\Tau),d)\,m_\Om$.  Applying the
inequality $r \, s\leq \Pi_M(r)+\Pi^*_M(s)$ on the right-hand side
above, we deduce \eqref{HypErrorMajoration}.

(ii) First notice that for all $\psi\in{\mathcal D}(\overline{Q})$, there exists $C>0$ such that
$$
\max_{n=1,\dots,N,\,\ptK\in\ptMrond,\ptL\in\ptNrond(\ptK)}
\frac{|\psi^n_\ptK\!-\!\psi^n_\ptKIL|}{d_{\ptK,\ptKIL}}\leq C,
\qquad \text{for all $h>0$}.
$$
Applying (i) for each $n$ and summing over $n=1,\dots,N$, we get
\begin{equation}\label{HypErrorMaj-2}
\begin{split}
\sum_{n=1}^N \Delt\; R^*_{\theta}[u^{\ptbarTau,n},\psi^n]
& \leq C
\inf_{\alpha>0} \left(\frac{\size(\Tau)}{\alpha}
\sum_{n=1}^N \Delt\; I_{\mathrm{Id}}[u^{\ptMrond},1] +T
\frac{1}{\alpha}\Pi^*_M\Bigl(C\alpha)\right)\\
& \leq C \inf_{\alpha>0} \left(\frac{\size(\Tau)}{\alpha} +
\frac{1}{\alpha}\Pi^*_M\Bigl(C\alpha)\right),
\end{split}
\end{equation}
where $C$  stands for a generic constant independent of $h$.

We have $(\Pi_M)'(0)=0$. Therefore
$$
(\Pi^*_M)'(0)=\lim_{b\to 0}\; \inf_{a}\left(a-\frac{\Pi_M(a)}{b}\right)
\leq \lim_{b\to 0}\left(b-\frac{\Pi_M(b)}{b}\right)=0.
$$
Hence for all $C>0$, $\lim_{\alpha\to 0} \frac{1}{\alpha}\Pi_M^*(C\alpha)=0$.
We deduce that the right-hand side
of \eqref{HypErrorMaj-2} tends to zero as $\size(\Tau)\to 0$.
\end{proof}

\begin{rem}\normalfont
Notice that if ${\mathfrak f}$ is locally Lipschitz continuous,
both $\Pi_M$ and $\Pi^*_M$ are quadratic; thus we can bound
$\Bigl| R_\theta[u^\ptbarTau,\psi] \Bigr|$ by
$\mathrm{Const} \, h^\beta$ for all $\beta<1/2$. Using the H\"older inequality
instead of the Young inequality, one recovers the result of
\cite{EyGaHe:book} with $\beta=1/2$.
Whenever ${\mathfrak f}$ is locally H\"older continuous of order
$\gamma\leq 1$, we find $\Pi^*_M(s)=\mathrm{Const}\, s^{1+\gamma}$. It
follows that $\Bigl| R_\theta[u^\ptbarTau,\psi] \Bigr|
\leq \mathrm{Const}\, h^\beta$ with $\beta=\frac{\gamma}{\gamma+1}$, under
the assumptions of Corollary \ref{UseOfWeakBVestim}(ii).
\end{rem}

\subsection{Approximate continuous entropy inequalities}
Relying on Proposition \ref{PropApproxEntrIneq}, we now deduce the
limiting (as $h\to 0$) entropy inequalities and the limiting
weak formulation; one should notice that they continue to hold if we replace
$(\eta_c^\pm,\mathfrak{q}_c^\pm)$ by regular ``boundary''
entropy-entropy flux pairs $(\eta_{c,\eps}^\pm,\mathfrak{q}_{c,\eps}^\pm)$.

\begin{prop}\label{PropContEntrIneq}
Consider a family of double meshes $\Tau$ of $\Om$ and
associated time steps $\Delt>0$, parametrized by $h=\max\{\size(\Tau),\Delt\}$, $h\to 0$.
Assume that $\reg(\Tau)$ is uniformly bounded. Denote the corresponding
discrete solutions of \eqref{AbstractScheme},\eqref{DiscrBC},\eqref{DiscrIC} by $u^{\ptbarTau,\barDelt}$.
Fix $\psi\in {\mathcal D}([0,T)\times\overline{\Om})$, $\psi\geq 0$, and
set $\psi^{\ptbarTau,\Delt}=\mathbb{P}^\ptbarTau\circ \mathbb{S}^\Delt[\psi]$.
Fix $\theta$ as one of the functions $\eta_c^\pm$, $c\in\R$.
Assume either $(c,\psi)\in \R^\pm\times\D([0,T)\times \overline{\Om})$,
or $(c,\psi)\in \R\times\D([0,T)\times \Omega)$. Then
\begin{equation}\label{ApproxContEntrIneq}
    \begin{split}
        \liminf_{h\to 0} &\;\;\Biggl(\; \int_{Q}
        \frac{1}{d}\biggl(\eta^\pm_c(u^{\ptMrond,\Delt})
        +(d-1)\,\eta^\pm_c(u^{\ptdMrond,\Delt})\biggr)\pt \psi \\
        & \quad +  \int_Q\frac{1}{d}\biggl({\mathfrak q}_c^\pm(u^{\ptMrond,\Delt})
        +(d-1)\,{\mathfrak q}_c^\pm(u^{\ptdMrond,\Delt})\biggr) \cdot\grad\psi \\
        & \quad\quad- \int_{Q} k(\grad^\ptTau w^{\ptbarTau,\Delt})
        \grad^\ptTau \;\widetilde A_{(\eta_c^\pm)'}(w^{\ptbarTau,\Delt}) \cdot\grad \psi\\
        & \quad\quad \quad +\int_\Omega \frac{1}{d}\biggl(\eta^\pm_c(u^{\ptMrond,0})
        +(d-1)\,\eta^\pm_c(u^{\ptdMrond,0})\biggr) \psi(0,\cdot)\\
        & \quad \quad\quad\quad
        + \int_Q \frac{1}{d}\biggl((\eta^\pm_c)'(u^{\ptMrond,\Delt})
        +(d-1)\,(\eta^\pm_c)'(u^{\ptdMrond,\Delt})\biggr) {\EuScript S}\psi\Biggr)
        \ge 0.
    \end{split}
\end{equation}
Furthermore, if $\psi\in{\mathcal D}([0,T)\times\Om)$, we have
\begin{equation}\label{ApproxContSolFaible}
    \begin{split}
        \lim_{h\to 0}\;
        \Biggl(\; &\int_{Q}
        \frac{1}{d}\biggl(u^{\ptMrond,\Delt}+(d-1)\,u^{\ptdMrond,\Delt}\biggr) \pt\psi \\
        & + \int_Q \Biggl(\frac{1}{d}\biggl({\mathfrak f}(u^{\ptMrond,\Delt})
        +(d-1)\,{\mathfrak f}(u^{\ptdMrond,\Delt})\biggr)
        -k(\grad^\ptTau w^{\ptbarTau,\Delt})
        \grad^\ptTau w^{\ptbarTau,\Delt}\Biggr) \cdot\grad \psi \!\\
        & \quad +\int_\Omega \frac{1}{d}\biggl(u^{\ptMrond,0}
        +(d-1)\,u^{\ptdMrond,0}\biggr) \psi(0,\cdot)
        +\int_Q\,{\EuScript S}\psi\; \Biggr)=0.
    \end{split}
\end{equation}
\end{prop}

\begin{proof}
By the choice of $(c,\psi)$, \eqref{CarrilloChoice} holds.
Thus, by Proposition \ref{PropApproxEntrIneq}, \eqref{ApproxEntrIneq}
and \eqref{ApproxSolFaible} hold; it suffices
to develop these formulas using the definitions of
$\Bleft\cdot\,,\,\cdot\Bright$, $\Aleft\cdot\,,\,\cdot\Aright$.

The second  term in \eqref{ApproxEntrIneq} rewrites exactly as the
corresponding term in \eqref{ApproxContEntrIneq}. Regarding the
other terms on the left-hand side, we also use the uniform  bound
on $u^{\ptTau,n}$ in $L^\infty$, the uniform  bound on
$k(\grad^\ptTau w^{\ptTau,\Delt})\grad^\ptTau w^{\ptTau,\Delt}$ in
$(L^{p'}(Q))^d$, and the convergences
$$
\begin{array}{l}
\dsp
 \sum_{n=1}^N
\sum_{\ptK\in\ptMrond} \!\frac{\psi_\ptK^{(n+1)}\!\! -\!
\psi_\ptK^{n} }{\Delt} \;\char_{Q_\ptK^n} \to \pt \psi, \;
\sum_{n=1}^N \sum_{\ptK\in\ptdMrond}
\!\frac{\psi_\ptdK^{(n+1)}\!\! - \! \psi_\ptdK^{n}
}{\Delt}\; \char_{Q_\ptdK^n} \to \pt \psi \;\text{in $L^1(Q)$}, \\[10pt]
 \dsp  \sum_{n=1}^N
\sum_{\ptK\in\ptMrond}
(\mathbb{P}^\ptTau(\mathbb{S}^\Delt[{\EuScript S}])^n)_\ptK
\char_{Q_\ptK^n} \!\!\to {\EuScript S}, \; \dsp  \sum_{n=1}^N
\sum_{\ptdK\in\ptdMrond}
(\mathbb{P}^\ptTau(\mathbb{S}^\Delt[{\EuScript S}])^n)_\ptdK
\char_{Q_\ptdK^n} \!\!\!\!\to {\EuScript S}\;\text{in $L^1(Q)$},
\\[10pt]
 \dsp  \sum_{n=1}^N
\sum_{\ptK\in\ptMrond} \psi^n_\ptK \;\char_{Q_\ptK^n} \to \psi, \;
\dsp  \sum_{n=1}^N \sum_{\ptdK\in\ptdMrond} \psi^n_\ptdK\;
\char_{Q_\ptdK^n} \to \psi\;\text{in $L^\infty(Q)$},
\\[15pt]
\dsp
\grad^\ptTau \psi^{\ptbarTau,\Delt}\to \grad\psi \;\text{in $L^p(Q)$}\; \;
\text{and}\; \; \psi^{\ptMrond,1}(\cdot)\to
\psi(0,\cdot),\; \psi^{\ptdMrond,1}(\cdot)\to \psi(0,\cdot)
\;\text{in $L^1(\Om)$},
\end{array}
$$
as $h\to 0$ (here we have put Proposition \ref{PropConsistency} to use).
Finally, the terms on the right-hand side of
\eqref{ApproxEntrIneq} vanish as $h\to 0$, thanks to the
initial remarks made in Subsection \ref{subsec:PenalVanish}
and Corollary \ref{UseOfWeakBVestim}(ii).
In the same way, \eqref{ApproxContSolFaible} follows from
\eqref{ApproxSolFaible}.
\end{proof}

\section{Convergence and statement of main result}\label{sec:convergence}

We are now in a position to state and prove the main result of this paper.

\begin{theo}\label{ThFVConverge}
Consider a family of double meshes $\Tau$ of $\Om$ and
associated time steps $\Delt>0$, parametrized
by $h=\max\{\size(\Tau),\Delt\}$, $h\to 0$.
Assume that $\reg(\Tau)$ is uniformly bounded.
Then the corresponding discrete solutions $u^{\ptbarTau,\barDelt}$
of \eqref{AbstractScheme},\eqref{DiscrBC},\eqref{DiscrIC} exist, are
uniformly bounded, and converge to the unique entropy solution
$u$ of \eqref{eq:prob1} in the following strong sense:
\begin{align*}
    & \text{$u^{\ptMrond,\Delt}\to u$, $u^{\ptdMrond,\Delt}\to u$
    in $L^s(Q)$ for any $s<\infty$},\\
    &\text{$\grad^\ptTau w^{\ptbarTau,\Delt}\to \grad w$
    in $L^p(Q)$}, \quad \text{where $w=A(u)$}.
\end{align*}
\end{theo}

\begin{proof}
We follow step by step the
proof of Theorem \ref{th:EPS-ESexists}.

{\bf (i)} Discrete solutions $u^{\ptbarTau,\barDelt}$ exist by
Proposition \ref{PropDiscrExist}. Besides, they verify the
asymptotic entropy inequalities \eqref{ApproxContEntrIneq}
(where we can replace $\eta_c^\pm$ by $\eta_{c,\eps}^\pm$) and the
asymptotic weak formulation \eqref{ApproxContSolFaible}, both of
Proposition \ref{PropContEntrIneq}.

{\bf (ii)} Proposition \ref{APrioriEstimates} yields uniform
estimates on  both $u^{\ptMrond,\Delt}$ and
$u^{\ptdMrond,\Delt}$ in $L^\infty(Q)$; on the time
translates of both $w^{\ptMrond,\Delt}$ and  $w^{\ptdMrond,\Delt}$
in $L^1(Q)$; on the penalization term
$\sum_{n=1}^N \Delt\Bleft \Prond^{\ptTau}[w^{\ptbarTau,n}], w^{ \ptTau,n} \Bright$;
and on $\grad^\ptTau w^{\ptbarTau,\Delt}$ in $L^p(Q)$.
The latter estimate implies further uniform estimates: namely, an
estimate of the space translates of both $w^{\ptMrond,\Delt}$ and
$w^{\ptdMrond,\Delt}$ in $L^1(Q)$, by Proposition \ref{PoincareAndTranslations} (ii);
an estimate of $\grad^\ptTau \widetilde A_{(\eta_{c,\eps}^\pm)'}(w^{\ptbarTau,\Delt})$
in $L^{p}(Q)$, because $\widetilde A_{(\eta_{c,\eps}^\pm)'}(\cdot)$ is
Lipschitz and by construction of $\grad^\ptTau[\cdot]$; and finally an
estimate of $\mathfrak{a}(\grad^\ptTau w^{\ptbarTau,\Delt})$ in
$L^{p'}(Q)$, because of the growth assumption on $\mathfrak{a}$.

{\bf (iii)} Thanks to the estimates of {\bf (ii)},
there exists a (not labelled) sequence of $\Tau$,$\Delt$ with $h\to 0$ such that

$\bullet$  by the Fr\'echet-Kolmogorov theorem, each of the
sequences $w^{\ptMrond,\Delt}$  and  $w^{\ptdMrond,\Delt}$
converges strongly in $L^1(Q)$ and pointwise a.e.~in $Q$;

$\bullet$  by Proposition \ref{DiscrRellich}, the limits of
$w^{\ptMrond,\Delt}$ and $w^{\ptdMrond,\Delt}$ coincide (we denote
the limit of $w^{\ptMrond,\Delt}$,$w^{\ptdMrond,\Delt}$  by $w$),
and $\grad^\ptTau w^{\ptbarTau,\Delt}$ converges weakly in
$L^p(Q)$ to $\grad w$;

$\bullet$ ${\mathfrak a}(\grad^\ptTau w^{\ptbarTau,\Delt})$
converges weakly in $L^{p'}(Q)$ to a limit field $\chi$;

$\bullet$ the sequences $u^{\ptMrond,\Delt}, u^{\ptdMrond,\Delt}$ converge
to $\mu,\mu^*:Q\times(0,1)\to \R$, respectively, in the sense of
nonlinear $L^\infty$ weak-$\star$ convergence \eqref{nonlinweakstar}.
Also by \eqref{nonlinweakstar}, the functions $w^{\ptMrond,\Delt}=A(u^{\ptMrond,\Delt})$
converge to $A(\mu)$ in the $L^\infty$ weak-$\star$ sense; since the functions
$w^{\ptMrond,\Delt}$ also converge strongly, $A(\mu)$ is
independent of $\alpha$ and coincides with $w$. In the same way,
we deduce that $A(\mu^*)$ is independent of $\alpha$ and coincides with $w$.
Also observe that $u^{\ptMrond,0}, u^{\ptdMrond,0}$ both
converge to $u_0$ a.e.~in $\Om$ and in $L^1(\Om)$.

{\bf (iv)} As in the proof of Theorem \ref{th:EPS-ESexists}, we can use the
chain rule and the Green-Gauss formula to deduce
\begin{equation}\label{ConvectDisappear}
    \begin{split}
        &\int_Q \int_0^1 \frac{1}{d}
        \left( {\mathfrak f}(\mu)+(d-1){\mathfrak f}(\mu^*)\right) \cdot\grad A(u)
        \\ & =\frac{1}{d}\int_0^1 \int_Q \left( {\mathfrak f}(\mu)\cdot \grad A(\mu)
        + (d-1) {\mathfrak f}(\mu^*)\cdot\! \grad A(\mu^*)\right)
        =\!\int_0^T \int_{\pOm} \widetilde A_{\mathfrak{f}}(w)\cdot n=0,
    \end{split}
\end{equation}
where $\widetilde A_{\mathfrak{f}}$ is defined in \eqref{eq:Af-KHK}.

{\bf (v)} Next we pass to the limit in \eqref{ApproxContSolFaible}. Indeed, by {\bf (iii)},
\begin{equation}\label{StartMinty}
    \begin{cases}
        \dsp \pt \widetilde u
        +\div \int_0^1 \frac{1}{d}\left( {\mathfrak f}(\mu)
        +(d-1) {\mathfrak f}(\mu^*)\right)\dal =\div \chi + {\EuScript S}
        \\ \text{\qquad in\; $L^{p'}(0,T;W^{-1,p'}(\Om))\!+\!L^1(Q)$},
        \quad \widetilde u|_{t=0}=u_0.
    \end{cases}
\end{equation}
where
$$
\widetilde u(t,x)=\int_0^1\widetilde \mu(t,x,\al) \dal, \qquad
\widetilde \mu = \frac{1}{d} \left( \mu + (d-1)\mu^*\right).
$$

Let us identify $\dsp \chi$ (the weak limit of ${\mathfrak a}(\grad^\ptTau w^{\ptbarTau,\Delt})$)
with ${\mathfrak a}(\grad w)$, and consequently obtain
that the weak convergence is in fact strong in $L^p(Q)$.  To this end, we will use {\bf (iv)}
and \eqref{ApproxTestFunctionA} to establish the inequality
\begin{equation}\label{IneqForMinty}
    \int_Q \chi\cdot \grad w \geq \liminf_{h\to 0}
    \sum_{n=1}^N\Delt\;\Aleft {\mathfrak a}(\grad^\ptTau w^{\ptbarTau,n}),\,
    \grad^\ptTau w^{\ptbarTau,n} \Aright.
\end{equation}

Indeed, using \eqref{StartMinty}, we can represent
the left-hand side of \eqref{IneqForMinty} as
\begin{equation}\label{StartMinty-KHK}
    \begin{split}
        \int_Q \chi\cdot \grad w\, \zeta
        & =-\int_{0}^T \langle \pt\widetilde u, w\, \zeta\rangle
        \\ & \qquad \
        + \int_0^1 \int_Q \frac{1}{d}\biggl( {\mathfrak f}(\mu)
        +(d-1){\mathfrak f}(\mu^*)\biggr) \cdot\grad w\, \zeta
        +\int_Q {\EuScript S}w\, \zeta,
    \end{split}
\end{equation}
where $\zeta\in \mathcal D([0,T))$ is nonincreasing with $\zeta(t)\equiv 1$ for $t$ small.

Note that since $A$ is nondecreasing, since $\widetilde u(t,x)$
is a convex combination of the values $\mu(t,x,\cdot)$ and $\mu^*(t,x,\cdot)$, and
because $A(\mu)=A(\mu^*)=w$, we conclude that
$$
w=A(\widetilde u).
$$

To control $\int_{0}^T \langle \pt\widetilde u, w\, \zeta\rangle$, we argue along the lines
of the proof of Theorem \ref{th:EPS-ESexists}. The duality
product $\langle \pt \widetilde u,A(\widetilde u)\rangle$ is treated via
the weak chain rule (cf.~\cite{AltLuckhaus}).  Hence, exploiting also the
convexity of $B(z)=\int_0^z A(s)\ds$,
\begin{equation}\label{eq:disc-time}
    \begin{split}
        &\int_0^T \langle \pt \widetilde u,A(\widetilde u)\,\zeta\rangle
        \\ &\qquad = - \int_Q B(\widetilde u) \zeta' - \int_\Om B(u_0)
        \\ &\qquad = \int_Q B\left(\int_0^1 \widetilde \mu(t,x,\al)\dal\right) (-\zeta') -\int_\Om B(u_0)
        \\ &\qquad \leq -\int_Q \, \zeta' \int_0^1 B( \widetilde \mu(t,x,\al)) \dal -\int_\Om  B(u_0).
    \end{split}
\end{equation}

Using \eqref{ConvectDisappear} and \eqref{eq:disc-time},  we
deduce from \eqref{StartMinty-KHK} that
\begin{equation}\label{IneqForMinty-2}
    \int_Q \chi\cdot \grad w\, \zeta
   \geq \int_Q \zeta' \int_0^1 B( \widetilde \mu(t,x,\al))\dal
    +\int_\Om  B(u_0)+ \int_Q \,{\EuScript S}w.
\end{equation}

On the other hand, Proposition~\ref{PropApproxEntrIneq} permits to
evaluate the right-hand side of \eqref{IneqForMinty} as follows:
\begin{equation}\label{IneqForMinty-3}
    \begin{split}
        & \liminf_{h\to 0} \left( \sum_{n=1}^N \Delt \Aleft \Bigl(
        {\mathfrak a}(\grad^\ptTau w^{\ptbarTau,n})  \, ,\,
        \grad^\ptTau w^{\ptbarTau,n}\, \zeta^{\ptbarTau,n} \Aright \right)
        \\ & \quad \le   \liminf_{h\to 0} \left (\sum_{n=1}^{N-1} \Delt \Bleft B(u^{\ptTau,n})\;,\;
        \frac{\zeta^{\ptTau,(n+1)} - \zeta^{\ptTau,n}}{\Delt}\Bright
        +\Bleft B(u^{\ptTau,0})\;,\;   1^{\ptTau} \Bright \right.
        \\ & \qquad\qquad \qquad\qquad
        +\sum_{n=1}^N \Delt
        \Bleft\mathbb{P}^\ptTau(\mathbb{S}^\Delt[{\EuScript S}])^n\,,\,w^{\ptTau,n}\zeta^{\ptTau,n}\Bright \Biggr)
    \end{split}
\end{equation}
By the previously established convergences (see also the proof of
Proposition~\ref{PropContEntrIneq}), the right-hand side of
\eqref{IneqForMinty-3} is equal to the right-hand side of
\eqref{IneqForMinty-2}. Once we let $\zeta$ tend to
$\char_{[0,T)}$, this establishes \eqref{IneqForMinty}.

Starting from \eqref{IneqForMinty}, we apply the Minty-Browder
argument that we employed for the continuous problem in the proof
of Theorem \ref{th:EPS-ESexists}.

Take $v\in L^p(0,T;W^{1,p}_0(\Om))\cap L^\infty(Q)$, and set
$v^{\ptbarTau,\Delt}=\mathbb{P}^\ptbarTau\circ
\mathbb{S}^\Delt[v]$. In view of \eqref{IneqForMinty}, taking into
account the strong convergence of $\grad^\ptTau
v^{\ptbarTau,\Delt}$ to $\grad v$ in $L^p(Q)$, cf.~Proposition
\ref{PropConsistency}, and the monotonicity of
$\mathfrak{a}(\cdot)$ we obtain
\begin{equation}\label{IneqForMinty-4}
    \begin{split}
        \int_Q \chi\cdot \grad (w-v)
        & \geq \liminf_{h\to 0} \sum_{n=1}^N \Delt\;
        \Aleft {\mathfrak a}(\grad^\ptTau w^{\ptbarTau,n}) \,,\, \grad^\ptTau w^{\ptbarTau,n}
        -\grad^\ptTau v^{\ptbarTau,n}\Aright
        \\ & \geq \liminf_{h\to 0} \sum_{n=1}^N \Delt\; \Aleft {\mathfrak a}(\grad^\ptTau v^{\ptbarTau,n})
        \,,\, \grad^\ptTau w^{\ptbarTau,n}- \grad^\ptTau v^{\ptbarTau,n}\Aright.
    \end{split}
\end{equation}
As is a well-known property of Leray-Lions operators, the
strong convergence of $\grad^\ptTau v^{\ptbarTau,\Delt}$ to $\grad v$ in $L^p(Q)$
implies the strong convergence of
$\mathfrak{a}(\grad^\ptTau v^{\ptbarTau,\Delt})$ to
$\mathfrak{a}(\grad v)$ in $L^{p'}(Q)$.
Therefore \eqref{IneqForMinty-4} yields
$$
\int_Q \chi\cdot \grad (w-v) \geq \int_Q \mathfrak{a}(\grad v)\cdot \grad (w-v).
$$
Choosing $v=w\pm\lambda\psi$ with $\lambda\downarrow 0$ and
$\psi\in L^p(0,T;W^{1,p}_0(\Om))$, we conclude
$$
\chi=\mathfrak{a}(\grad w).
$$
Moreover, as in the proof of Theorem \ref{th:EPS-ESexists} and
\cite[Theorem 5.1]{ABH``double''}, relying on the strict monotonicity
of ${\mathfrak a}$ and utilizing an argument of \cite{Brow,BocMuratPuel}, we
also deduce the strong convergence of $\grad^\ptTau
w^{\ptbarTau,\Delt}$ to $\grad w$ in
$L^p(Q)$.

{\bf (vi)} Now we can pass to the limit in the weak and entropy formulations
listed in Proposition \ref{PropContEntrIneq}. The passage from \eqref{ApproxContSolFaible}
to \ref{def:entropyWeakp} is straightforward. In \eqref{ApproxContEntrIneq}, we first
work with regularized boundary entropies.  Taking the limit, all the terms converge to the corresponding
terms in \ref{def:entropy2p} in a straightforward way, except
for the third one. Let us show that
$$
\text{$\grad^\ptTau \widetilde A_{(\eta_{c,\eps}^\pm)'}( w^{\ptbarTau,\Delt})$
converges weakly to $\grad \widetilde A_{(\eta_{c,\eps}^\pm)'}(w)$ in $L^p(Q)$.}
$$
Indeed, both $\widetilde A_{(\eta_{c,\eps}^\pm)'}(w^{\ptMrond,\Delt})$ and
$\widetilde A_{(\eta_{c,\eps}^\pm)'}(w^{\ptdMrond,\Delt})$ converge
to $\widetilde A_{(\eta_{c,\eps}^\pm)'}(w)$ by the a.e.~convergence
of $w^{\ptMrond,\Delt}$,$w^{\ptdMrond,\Delt}$  to $w$
and the continuity of $\widetilde A_{(\eta_{c,\eps}^\pm)'}$. Using
the boundedness in  $L^p(Q)$ of $\grad^\ptTau \widetilde
A_{(\eta_{c,\eps}^\pm)'}(w^{\ptbarTau,\Delt})$) and the
compactness property of Proposition~\ref{DiscrRellich}, we conclude
that our claim holds. The subsequent arguments are the same as in the proof of
Theorem~\ref{th:EPS-ESexists}.

{\bf (vii)}  We conclude that $(\mu,\mu^*,w)$ is an entropy double-process solution of \eqref{eq:prob1}.
In view of Theorem \ref{th:EPSunique}, this brings to an end the
proof of Theorem \ref{ThFVConverge}; indeed, we obtain
the convergence to $u$ for each sequence of discrete solutions with $h\to 0$. Also, the fact
that $\mu$ and $\mu^*$ turn out to be independent of $\alpha$ means that
the convergence of $u^{\ptMrond,\Delt},u^{\ptMrond,\Delt}$ to $u$
is strong in $L^s(Q)$ for all finite $s$.
\end{proof}

\section{On the choice of FV scheme and various generalizations}\label{sec:Extensions}
In this section we discuss other possible choices of finite volume schemes for \eqref{eq:prob1}.

$\bullet$ The use of DDFV schemes is motivated by their
convenience when it comes to the discretization of nonlinear diffusion operators.
Other possibilities exist; among them, let us mention the schemes
studied in \cite{HandMikulaSgallari} (see also
\cite{AfifAmaziane,AGW,ABK-FVCA5}), in \cite{ABH-M2AN}, in
\cite{Droniou} (see also \cite{EymardDroniou}), and
in \cite{EGH:Sushi-IMAJNA} (see also
\cite{EGH:IMAJNA,EGH:AnisoCRAS,EH-SUCCES}). All these schemes
possess some variant of the ``integration-by-parts'' property of
Proposition~\ref{DualityDiff}.

The 2D schemes of \cite{ABH-M2AN} are restricted to Cartesian
meshes, so they do not allow for domains much more general than
rectangles. Notice that their generalization to 3D appears to be
straightforward. The techniques used in the present paper and in
the references we cite, such as \cite{EyGaHe:book,EGHMichel},
combined with those of \cite{ABH-M2AN}, allow to design rather
simple FV schemes on Cartesian meshes for problem \eqref{eq:prob1}
and to prove their convergence. In this case, the notion of
entropy double-process solution is not needed, and the theoretical
results in \cite{EGHMichel} can be adapted directly.

 This is
also the case of the ``complementary volumes'' schemes as
described in \cite{HandMikulaSgallari}. In 2D,  ideas quite
similar to that of \cite{HandMikulaSgallari} were used to
construct the schemes of \cite{AfifAmaziane,AGW,ABK-FVCA5}. All
these schemes
 work on meshes dual to conformal triangular 2D meshes, and
 the discrete gradient is reconstituted by affine per triangle interpolation.
 ``Complementary volumes'' schemes
are
 simpler than our DDFV scheme from the practical point of view,
since one discretizes the problem on the same mesh $\dMrond$
using, roughly speaking, half of the unknowns. The discrete
duality properties for the 2D ``complementary volumes'' scheme are
shown in the same way as for our DDFV schemes; the proof is based
upon Lemma~\ref{lem:ReconstrPropTriangles} (see Appendix B
and also \cite{AndrBend,ABK-FVCA5}).
Unfortunately, the straightforward generalization of these
``complementary volumes'' schemes to 3D fails to satisfy the
discrete duality property, except for very constrained geometries
of the meshes (see Remark~\ref{rem:ReconstrGeneralizationMultiD}).

The key feature of the  2D schemes of
\cite{AfifAmaziane,HandMikulaSgallari,AGW,ABH-M2AN,ABK-FVCA5} (see
also \cite{CoudiereVilaVilledieu})
 lies in the fact that the fluxes
across interfaces are reconstructed ``manually''. The approaches
of Droniou and Eymard \cite{EymardDroniou,Droniou} and those of
the HVF, SUCCES and SUSHI schemes of Eymard, Gallou\"et and Herbin
\cite{EGH:IMAJNA,EGH:AnisoCRAS,EH-SUCCES,EGH:Sushi-IMAJNA} are
different; they rely on introducing additional unknowns (either
for the fluxes, or for the values on some of the edges) and on
careful penalization of the finite differences.

The schemes HVF, SUCCES and SUSHI (among many others) were
designed for handling linear anisotropic, heterogeneous diffusion
problems with possibly discontinuous coefficients; in this
framework, their convergence is justified.  These schemes avoid
usage of double meshes and thus may have less unknowns; they work
both in 2D and 3D. We refer to Eymard, Gallou\"et and Herbin
\cite{EGH:Sushi-IMAJNA} for the description and comparison of
these and related (e.g., mimetic finite difference) schemes. 
Finally, let us also mention the schemes of 
Aavatsmark et al.~(see, e.g., \cite{Aavatsmark,Aavatsmark-et-al}), that 
are in a sense intermediate. The gradient reconstruction used in
\cite{Aavatsmark,Aavatsmark-et-al} also involves additional edge
unknowns, which are eliminated by solving, locally, an algebraic
system of equations.

The scheme of \cite{Droniou} designed for nonlinear Leray-Lions
kind problems can be directly compared to the DDFV schemes of
\cite{ABH``double''} and of the present paper. The scheme of
\cite{Droniou} is very interesting because of the extreme
generality of the geometries allowed for the mesh (and it works in
any space dimension). For this same reason, theoretical
justification of its convergence in the hyperbolic-parabolic
framework \eqref{eq:prob1} seems problematic. Indeed, the
conformity (orthogonality) condition was used in an essential way
in the derivation of the discrete entropy inequalities (see Remark
\ref{eq:ConditionsForEntropyIneq}). The same difficulty arises for
the double schemes of \cite{ABH``double''} in the case of
non-conformal meshes, cf.~Remark \ref{rem:conformal}.
In passing, let us point out that the conformity (orthogonality)
assumption on the meshes is the only condition that is known to
ensure the discrete maximum principle for the DDFV schemes.

In conclusion, the 2D and 3D conformal DDFV schemes studied herein,
although constrained by the orthogonality condition, by the
Delaunay condition, and by condition \eqref{eq:AssumptMesh},
combine some degree of flexibility (e.g., any polygonal/polyhedral
domain can be partitioned into triangles/tetrahedra
satisfying these restrictions) with the rigid structure properties
underlying our convergence proof. But because of the
conformity constraint, the advantage of simple local refinement
procedures for 2D DDFV schemes, pointed out in \cite{ABH``double''}, is lost.

$\bullet$ Our assumption that $\Mrond$ consists of simplexes is a
practical one simplifying the presentation of the scheme. In 2D,
it can be replaced by the more general assumption that any element
of $\Mrond$ admits a circumscribed circle. In 3D, we can assume
that each $\K\in \overline{\Mrond}$ admits a circumscribed ball,
and each interface $\KIL$ is a triangle satisfying
\eqref{eq:AssumptMesh}.

Notice that Remark \ref{rem:ReconstrGeneralization2DPolygon} (see
also \cite{AndrBend}) makes it possible to define a consistent
discrete duality scheme even when the interfaces $\KIL$ are not
necessarily triangles. Unfortunately, the discrete Poincar\'e
inequality may fail in this generality; this undermines the
subsequent convergence analysis. Yet one interesting case is
that of a Cartesian mesh $\Mrond$; the corresponding
DDFV schemes are alternatives to the scheme of \cite{ABH-M2AN}
discussed above. More generally,
one can start with a mesh $\Mrond$ made of rectangles (e.g.,
inside $\Om$) and triangles (e.g., near the boundary $\ptl\Om$) in 2D.

$\bullet$ As pointed out in Remark \ref{rem:Gen4D}, a different
kind of reconstruction formula is needed for problems in 4D and
higher dimensions. It would be interesting to conceive discrete
gradients consistent with affine functions, following the
principle formulated in Remark \ref{ExplicDiscrGrad}. One natural
way is indicated in \cite{CoudiereHubertPierre}.

$\bullet$ The choice of penalization in our double scheme can be
changed (see Remark \ref{rem:penalChoice}). One could also
penalize the differences $(\uK-\udK)$ instead of the differences
$(\wK-\wdK)=(A(\uK)-A(\udK))$; this would permit to avoid the use
of double-process solutions. But this choice would introduce
additional coupling between the sets of variables
$(\uK)_{\ptK\in\ptMrond}$ and $(\udK)_{\ptdK\in\ptdMrond}$ in the
``hyperbolic'' regions. Indeed, if, e.g., $A(u)\equiv 0$, there is
no coupling at all between the variables sitting on $\Mrond$ and
those sitting on $\dMrond$. Therefore our choice seems more
convenient in terms of practical implementation.

$\bullet$  Convection-diffusion problems with anisotropic linear
and nonlinear diffusion were considered in
\cite{ChenPerthame,ChenKarlsen} and in
\cite{BenKar:Renorm,BenKar:RenormII}. General DDFV schemes
do not seem easy to adapt to the nonlinear anisotropic framework,
because of the presence of ``privileged" directions of diffusion.
In this case, the schemes of \cite{ABH-M2AN} on Cartesian
meshes constitute a natural choice, and the
geometry of $\ptl\Om$ should be rather
taken into account via the approximation of the domain $\Omega$ by
domains with piecewise axes-aligned boundaries. Notice that for
the anisotropic $p-$Laplace kind diffusions
$$
\partial_{x_1}\Bigl(|\ptl_{x_1} A_1(u)|^{p_1-2}\ptl_{x_1} A_1(u)\Bigr)
+ \partial_{x_2} \Bigl(|\ptl_{x_2} A_2(u)|^{p_2-2}\ptl_{x_2} A_2(u) \Bigr)
$$
considered by Bendahmane, Karlsen in
\cite{BenKar:Renorm,BenKar:RenormII}, the discrete entropy
inequalities on Cartesian meshes are as easy to obtain as for the
isotropic case $\mathfrak{a}(\xi)=k(\xi)\xi$ considered in the
present paper.

$\bullet$ Taking into account sufficiently smooth dependencies on
$(t,x)$ of the convection and diffusion operators is possible,
although quite technical; see \cite{EyGaHe:book,ABH``double''} for
some results in that direction, and also \cite{Chen:2005wf,KO:Unique} for
well-posedness results for degenerate equations
with $(t,x)$ dependent convection-diffusion operators.
Discontinuous coefficients are important for the modeling of fractured media.
DDFV schemes for Leray-Lions
operators $\div \mathfrak{a}(x,\grad w)$
with discontinuous (piecewise smooth) in $x$ nonlinearity
$\mathfrak{a}$ are studied in the recent work \cite{AngotetCo}.
The case of $x$-discontinuous flux functions $\mathfrak{f}(x,u)$ has received
much attention in the last fifteen years (see, e.g., \cite{BurgKarlMishTow}
and the references cited therein), both from a theoretical and
numerical perspective. Let us mention here that
the problem of the choice of the appropriate entropy conditions
strongly depends on the underlying physical interpretation;
different models lead to qualitatively different admissible solutions.

$\bullet$ Inhomogeneous Dirichlet boundary conditions can be taken
into account, combining the techniques of \cite{MichelVovelle}
with those of \cite{ABH``double''}; both are rather involved, which
explains our choice of the homogeneous boundary data for the
presentation of the scheme and the convergence arguments.

\section*{Appendix A: Proof of uniqueness}
\addtocounter{section}{1}

This appendix is devoted to a proof of Theorem \ref{th:EPSunique}.
The proof is an adaptation of the ones in Carrillo \cite{Carrillo}
(for entropy solutions) and that in Eymard, Gallou\"et, Herbin,
Michel \cite{EGHMichel} (for entropy process solutions, which can
be viewed as entropy double-process solutions with
$\mu\equiv\mu^*$). The proof  is mainly divided into several
lemmas (Lemmas \ref{lem:chain}, \ref{lem:unique1},
\ref{lem:unique2} below). For simplicity, let us only consider the
case where the source term ${\EuScript S}$ is zero (see also
Remark \ref{rem:L1CCPrinciple}).

We begin by introducing the set
$$
\E=\{r\in\R:\text{$A^{-1}(r)$ is neither empty nor a singleton}\},
$$
and proving
\begin{lem}
\label{lem:chain} Let $(\nu,\nu^\star,v)$ be an entropy
double-process solution of \eqref{eq:prob1} with initial data
$v_0$. Then for all $W\in \R^d$, for any $\phi \in \D([0,T)\times
\Omega)$, $c\in \R$ such that $A(c)\not \in \E$ and also for any
$\phi \in \D([0,T)\times \oOm)$, $c\in \R^\pm$ such that $A(c)\not
\in \E$, we have with the notation of Section~\ref{sec:def}
the following equality:
\begin{equation}\label{eq:gen-entr-proc}
    \begin{split}
        &\int_0^1\!\!\int_{Q} \Biggl[ \frac{1}{d}\Bigl(\eta_c^\pm(\nu)
        +(d-1)\eta_c^\pm(\nu) \Bigr)\;\pt \phi
        +\frac{1}{d}\Bigl(\,\mathfrak{q}_c^\pm(\nu)
        \,+\,(d-1)\mathfrak{q}_c^\pm(\nu^\star)\,\Bigr)\cdot\grad \phi
        \\ & \qquad \quad
        - \sign^\pm(v-A(c)) \Bigl(\mathfrak{a}(\Grad v) -W\Bigr)
        \Bigr)\cdot\Grad \phi \Biggr] \dx \dt \dalpha
        +\int_\Omega \eta_c^\pm(v_0) \phi|_{t=0}\dx
        \\ & =\lim_{\eps\dto 0} \int_{Q}(\sign^\pm_\eps)'(v-A(c))
        \Bigl(\mathfrak{a}(\Grad v)- W \Bigr)\cdot \Grad v\,  \phi\; \dx\dt.
    \end{split}
\end{equation}
\end{lem}

\begin{proof}
We refer to \cite[Lemma 1]{Carrillo} and to \cite{EGHMichel} for
details on the proof. The idea is to use
$\psi:=(\sign^\pm_\eps(v-A(c))\phi)$ as a test function in
\ref{def:entropyWeakp}. It is admissible; indeed, we can
approximate it by functions in ${\mathcal D}([0,T)\times\Om)$ and
pass to the limit in all terms of \ref{def:entropyWeakp}, because
$\psi\in L^\infty(Q)\cap L^p(0,T;W^{1,p}_0(\Om))$ for any of the
two possible choices of $(\phi,c)$ (in particular, notice that
$\sign^\pm_\eps(v-A(c))\in L^p(0,T;W^{1,p}_0(\Om))$ in case $c\in\R^\pm$).
We have
$$
\sign^\pm(\nu-c)=\sign^\pm(v-A(c))=\sign^\pm(\nu^*-c),
$$
thanks to the relation \ref{def:entropy1p} (which reads
$A(\nu)\equiv v\equiv A(\nu^*)$ in our notation) and to the choice
of $A(c)\notin \E$; then we use the weak chain rule
to deal with the time derivative.
We also insert into \ref{def:entropyWeakp} the term
$$
\int_0^1\!\!\int_Q \sign^\pm_\eps(v-A(c)) \grad\phi \cdot W-
\int_0^1\!\!\int_Q (\sign^\pm_\eps)'(v-A(c)) \grad w\cdot W \,\phi,
$$
which is equal to $0=\dsp \int_Q \div (W\,\sign^\pm_\eps(v-A(c))\phi)$
for any of the two possible choices of $(\phi,c)$, by the Gauss-Green formula.
As $\eps \down 0$, the term containing
$(\sign^\pm_\eps)'(v-A(c))\;\Bigl(\mathfrak{f}(\nu)-\mathfrak{f}(c)\Bigr)\cdot
\grad v $ vanishes, as shown in \cite[Lemma 1]{Carrillo}.
\end{proof}

We are now interested in comparing two entropy double-process
solutions of \eqref{eq:prob1}, denoted by $(\nu,\nu^*,\theta)$ and $(\mu,\mu^*,w)$), of
which the first one is chosen to satisfy $\nu\equiv \nu^*$.
Consider the distribution ${\mathcal I}$ on $\D(\overline{Q})$ defined by
\begin{equation}\label{eq:DefDistrI}
    \begin{split}
        & {\mathcal I}[\phi]:=\int_0^1\!\!\int_0^1\!\!\int_{Q}
        \Biggl[ \frac{1}{d}\Bigl((\nu\!-\!\mu)^++(d-1)(\nu\!-\!\mu^\star)^+ \Bigr)\;\pt \phi
        \\ &\qquad\qquad\qquad\qquad\quad
        +\frac{1}{d}\Bigl(\,\sgn{\nu\!-\!\mu}(\mathfrak{f}(\nu)\!-\!\mathfrak{f}(\mu))
        \\ &\qquad\qquad \qquad\qquad \quad
        +(d-1)\sgn{\nu\!-\!\mu^\star}(\mathfrak{f}(\nu)\!-\!\mathfrak{f}(\mu^\star))\,\Bigr)\cdot\grad \phi \\
        & \qquad\qquad\qquad\qquad\quad\quad\quad
        -\sgn{v\!-\!w} \Bigl(\mathfrak{a}(\Grad v)\!
        -\!\mathfrak{a}(\Grad w)\Bigr)\cdot\Grad \phi \Biggr] \dx \dt \dalpha \dbeta
        \\ &\qquad\qquad\qquad
        +\int_\Omega (v_0\!-\!u_0)^+ \phi(0,x)\dx.
    \end{split}
\end{equation}

Let us prove that we can write ${\mathcal I}$ as
\begin{equation}\label{eq:I-split}
    {\mathcal I}={\mathcal {IP}}+{\mathcal{IN}}
\end{equation}
where ${\mathcal {IP}}[\phi]$ is defined by the analogue of
\eqref{eq:DefDistrI} with each of $\nu,v,v_0,\mu,\mu^*,w,u_0$
replaced by its positive part; and ${\mathcal {IN}}[\phi]$ is
defined by the analogue of \eqref{eq:DefDistrI} with each of
$\nu,v,v_0,\mu,\mu^*,w,u_0$ replaced by
$-\nu^-,-v^-,-v_0^-,-\mu^-,-(\mu^*)^-,-w^-,-u_0^-$.
To emphasize, whenever necessary, the dependency of
$\mathcal{ I,IP, IN}$ on the involved solutions, we will write
$\mathcal{ I}^{\nu,\nu,v}_{\mu,\mu^*\!,w},,\mathcal{IP}^{\nu,\nu,v}_{\mu,\mu^*\!,w},
\mathcal{IN}^{\nu,\nu,v}_{\mu,\mu^*\!,w}$, respectively.

To justify \eqref{eq:I-split}, we use the identity
\eqref{eq:Split(a-b)+} from the following easy lemma.

\begin{lem}\label{lem:F-aplus-amoins-splitting}
 For all $F:\R\to\R$ such that $F(0)=0$, for all $a,b\in\R$ there holds
\begin{equation}\label{eq:Split(a-b)+}
    \begin{split}
        \sign^+(a\!-\!b)\;(F(a)\!-\!F(b)) & =\sign^+(a^+\!-\!b^+)\;(F(a^+)\!-\!F(b^+))
        \\ &\quad+\sign^+((-a^-)\!-\!(-b^-))\;(F(-a^-)\!-\!F(-b^-)),
    \end{split}
\end{equation}
and
\begin{equation}\label{eq:SplittingIntoPlusMinusParts}
\left|
\begin{array}{l}
(i)\quad
\sign^-(b-a^+)\;F(b)=-\sign^+(a^+-b^+)\;F(b^+)\;+\;\sign^-(b)\;F(b),\\
(ii)\quad \sign^-(b-a^+)\;F(a^+)=-\sign^+(a^+-b^+)\;F(a^+).
\end{array} \right.
\end{equation}
\end{lem}

We apply  \eqref{eq:Split(a-b)+} to $a=\nu$, $b=\mu$ (or
$b=\mu^*$) with $F=\mathrm{Id}$ and with $F=\mathfrak{f}_i$,
$i=1,\dots,d$. Futhermore, observe that the analogue of
\eqref{eq:Split(a-b)+} still holds for a.e. $(t,x)\in Q$ if we
take $a=v(t,x)$, $b=w(t,x)$  and replace $F(a),F(b)$ and 
$F(\pm a^\pm),F(\pm b^\pm)$ by $\mathfrak{a}(\Grad v)$,$\mathfrak{a}(\Grad w)$ 
and by $\mathfrak{a}(\pm\Grad v^\pm)$,$\mathfrak{a}(\pm\Grad  w^\pm)$, respectively. 
Indeed, we have, e.g., $\mathfrak{a}(\grad v) =\mathfrak{a}(\grad
v^+)+\mathfrak{a}(\grad v^-)$ a.e.~on $Q$, because $v\in
L^p(0,T;W^{1,p}_0(\Om))$. Using all aforementioned identities, we
split each term in the definition \eqref{eq:DefDistrI} of
${\mathcal I}$ into the sum of the corresponding terms in the
definitions of  $\mathcal {IP}$ and $\mathcal {IN}$.

Now we estimate $\mathcal I$ ``inside" the domain''.

\begin{lem}\label{lem:unique1}
Let $(\nu,\nu,v)$ and $(\mu,\mu^\star,w)$ be
entropy double-process solutions of \eqref{eq:prob1} with data
$v_0,u_0\in L^\infty(\Omega)$, respectively. Then
${\mathcal I}[\phi] \geq 0$, $\forall \phi \in \D([0,T)\times \Omega)$, $\phi\ge 0$.
\end{lem}

\begin{proof}
The proof is an application of the doubling of variables
method of Kruzhkov \cite{Kruzkov}; it follows
\cite{Carrillo,EGHMichel,BenKar:Renorm}. We let $\nu$ depend on
variables $(t,x,\alpha)\in Q\times (0,1)$ and $\mu$ depend on
another set of variables $(s,y,\beta)\in Q\times (0,1)$. In what
follows, $\grad v$ means $\grad_x v$ and $\grad w$ means $\grad_y w$.
As to the test function $\phi$, it will depend on the
variables $(t,s,x,y)$, thus we will use the notations $\pt,\ps$
and $\Grad_x,\Grad_y$ for the corresponding derivatives of $\phi$.
We will work with nonnegative test
functions $\phi\in {\mathcal D}\left(([0,T)\times\Om)^2\right)$.
Let us introduce the sets on which the diffusion term
for the first, respectively, for the second solution degenerates:
$$
\E_\nu=\{(t,x)\in Q\,|\,v(t,x)\in \E\},\quad \E_\mu=\{(s,y)\in Q\,|\,w(s,y)\in \E\}.
$$
Denote by $\E_\nu^c$, $\E_\mu^c$ the complementary sets in $Q$ of
$\E_\nu$, $\E_\mu$, respectively. Observe that $\Grad v=0$ a.e.~in
$\E_\nu$ and $\Grad w=0$ a.e.~in $\E_\mu$ (recall \ref{def:entropy1p}).

{\bf (i)} First we apply Lemma \ref{lem:chain} with the
solution $(\nu,\nu,v)$. For all $(s,y,\beta)\in \E_\mu\times (0,1)$,
choose $W=\mathfrak{a}( \grad w(s,y))$ and take
the entropy $\eta_c^+(\cdot)=(\cdot-c)^+$ with
$c=\mu(s,y,\beta)$, then with $c=\mu^\star(s,y,\beta)$ in
\eqref{eq:gen-entr-proc}. We multiply the two resulting equations
by $\frac{1}{d}$ and by $\frac{(d-1)}{d}$, respectively, and
add them together. Then  we integrate in $(s,y,\beta)\in \E_\mu\times
(0,1)$. Similarly, for $(s,y,\beta)\in \E_\mu^c\times (0,1)$, we
add together, with weights $\frac{1}{d}$ and $\frac{(d-1)}{d}$,
respectively, the entropy inequalities \ref{def:entropy2p} for
$(\nu,\nu,v)$ corresponding to $\eta_c^+(\cdot)$ with
$c=\mu(s,y,\beta)$ and with $c=\mu^\star(s,y,\beta)$.  We
integrate the resulting inequality in $(s,y,\beta)\in \E_\mu^c\times (0,1)$.

{\bf (ii)} Next, we exchange the roles of $(\nu,\nu,v)$
and $(\mu,\mu^*,w)$. This time we use the entropy
$\eta_c^-(\cdot)=(\cdot-c)^-$; we use $W=\mathfrak{a}(\grad
v(t,x))$; and we only use one value $c=\nu(t,x,\alpha)$ in the
analogue of \eqref{eq:gen-entr-proc} (for all $(t,x,\alpha)\in
\E_\nu\times (0,1)$) and in the analogue of \ref{def:entropy2p}
(for all $(t,x,\alpha)\in \E_\nu^c\times (0,1)$).

{\bf (iii)} Adding the inequalities obtained in {\bf (i),(ii)}, by the symmetry
of the expressions involved (such as $(\nu-\mu)^+=(\mu-\nu)^-$, etc.), we get, keeping in mind
Remark \ref{RemTwoFormsOfDiffusion}, the following inequality:
\begin{equation}\label{bis:eq:entropy_u_p_ilus_v_I}
\!\!\!\!\!\!\begin{split}
        \!\!\!\!\!\!\!\!\!\!\!\!&\int_0^1\!\!\int_0^1 \iint_{Q\times Q }
        \Biggr[\frac{1}{d}\Bigl(\,(\nu\!-\!\mu)^++(d-1)(\nu\!-\!\mu^\star)^+\Bigr)(\pt\!+\!\ps)\phi \\
        \!\!\!\!\!\!\!\!\!\!\!\!& \quad\qquad \qquad\qquad
        + \frac{1}{d}\Bigl(\;\sgn{\nu\!-\!\mu}(\mathfrak{f}(\nu)\!-\!\mathfrak{f}(\mu))
        \\
        \!\!\!\!\!\!\!\!\!\!\!\!& \quad \qquad\qquad\qquad\qquad
        +(d-1)\sgn{\nu\!-\!\mu^\star}(\mathfrak{f}(\nu)\!-\!\mathfrak{f}(\mu^\star))\;\Bigl)
        (\Grad_x\!+\!\Grad_y) \phi \\
        \!\!\!\!\!\!\!\!\!\!\!\! &
        \quad\qquad\qquad\qquad
        - \sgn{v\!-\!w}\Bigl(\mathfrak{a}(\Grad v)
        \!-\!\mathfrak{a}(\Grad  w)\Bigr) \cdot \Bigl( \Grad_x\!
        +\!\Grad_y  \Bigr)\phi \Biggl]\;dx \dt \dy\ds \dalpha\dbeta\\
        \!\!\!\!\!\!\!\!\!\!\!\! &
        \qquad +\int_0^1\iint_{\Omega\times Q}
        \frac{1}{d}\Bigl(\,(v_0\!-\!\mu)^++(d-1)(v_0\!-\!\mu^\star)^+\Bigr)\,\phi\dx(\dy\ds)\dbeta
        \\
        \!\!\!\!\!\!\!\!\!\!\!\!& \qquad
        +\int_0^1 \iint_{Q\times \Omega} (\nu\!-\!u_0)^+\,\phi\, (\dx\dt)\dy\dalpha
        \\
        \!\!\!\!\!\!\!\!\!\!\!\! & \qquad \geq \lim_{\eps\dto 0}
        \iint_{\E_\nu^c\times \E_\mu^c}
        \dsgne{v\!-\!w}\,
        \Bigl(\mathfrak{a}(\Grad v)\!-\! \mathfrak{a}(\Grad  w)\Bigr)
        \cdot\Bigl(\Grad v\!- \!\Grad w\Bigr)\;\phi\dx\dt\dy\ds.
    \end{split}\!\!\!\!\!\!\!\!\!\!\!\!\!\!\!\!\!\!
\end{equation}
The last term in \eqref{bis:eq:entropy_u_p_ilus_v_I} is
nonnegative, because $\mathfrak{a}$ is monotone and $\phi\geq 0$.

{\bf (iv)} Let us now specify the test function. For $l,n\in\N$,
let $ \omega_n:\R^d\to\R$, $\omega^l:\R\to\R$ be standard
symmetric mollifiers with supports in $\{x\in\R^d\,|\,\|x\|\leq
\frac{1}{n}\}$ and in $\{t\in\R\,|\,|t|\leq \frac{1}{l}\}$,
respectively. We take the test function in
\eqref{bis:eq:entropy_u_p_ilus_v_I} to be
\begin{equation*}
\phi_{n,l}(t,x,s,y)=\phi\left(x,t\right)
\omega_n\left({x-y}\right) \omega^l\left({t-s}\right) \equiv \phi
\omega_n\omega^l,
\end{equation*}
where $\phi\in \D([0,T)\times \Om)$, $\phi\geq 0$. With this
choice, we have
\begin{equation}\label{bis:eq:test_tmp1}
    (\pt + \ps)\phi_{n,l}
    = \Bigl(\pt\phi\Bigr)\,\omega_n\omega^l,
    \quad
    \left(\Grad_x+\Grad_y\right)\phi_{n,l}
    = \Bigl(\Grad_x
    \phi\Bigr)\, \omega_n \omega^l.
\end{equation}
 Then we let $n,l\to\infty$. The first term in
\eqref{bis:eq:entropy_u_p_ilus_v_I} converges to the first term in
the right-hand side of \eqref{eq:DefDistrI}. This argument is
standard; one can use, e.g., the properties of the Lebesgue points
of $L^1$ functions and the upper-semicontinuity of the $L^1$
``$+$-bracket''
$$
\Bigl[u,f  \Bigr]_+:=\int_{\Om} \sgn u\,f + \int_{\{u=0\}} f^+.
$$
The two latter terms in the left-hand side of
\eqref{bis:eq:entropy_u_p_ilus_v_I} are treated with the help of
the triangular inequality and of the strong initial trace property
\eqref{eq:IC-strong} proved in Lemma~\ref{lem:IC} below. The
limit, as $n,l\to\infty$, of each of these terms is majorated by
one half of the last term in \eqref{eq:DefDistrI} (this is because
$\int_Q \omega^l(t)\,dt=\frac 12 = \int_Q \omega^l(-s)\,ds$). This
concludes the proof of the lemma.
\end{proof}
\begin{lem}\label{lem:IC}
Let $(\mu,\mu^\star,w)$ be an entropy double-process solution of
\eqref{eq:prob1} with initial datum $u_0\in L^\infty(\Omega)$.
Then the initial datum is also taken in the following strong
sense:
\begin{equation}\label{eq:IC-strong}
\lim_{h\downarrow 0} \frac 1h \int_0^h\!\!\int_\Om\!\int_0^1
        \biggl(\,\frac{1}{d}|\mu\!-\!u_0|+\frac{d\!-\!1}{d}|\mu^\star\!-\!u_0|\biggr)\,\dt\dx\dalpha=0.
\end{equation}
\end{lem}

Notice that another way to formulate \eqref{eq:IC-strong} is to say that
$$
\text{\rm ess}\lim_{\!\!\!\!\!t\downarrow 0} \int_\Om\!\int_0^1
\biggl(\,\frac{1}{d}|\mu\!-\!u_0|+\frac{d\!-\!1}{d}|\mu^\star\!-\!u_0|\biggr)
\dx\dalpha=0,
$$
in the spirit of the original definition of Kruzhkov
\cite{Kruzkov}.

\begin{proof}
The proof follows the one of Panov in \cite[Proposition 1]{Panov-Izv}. 
For $c\in\R$ and $h>0$, consider the functions
$$
p_h(\cdot;c):\;x\in\Om \mapsto \frac 1h \int_0^h\!\!\int_0^1
\biggl(\,\frac{1}{d}|\mu(t,x;\alpha)\!-\!c|
+\frac{d\!-\!1}{d}|\mu^\star(t,x;\alpha)\!-\!c|\biggr)\,\dt\dalpha .
$$
Because $\mu,\mu^*$ are bounded, the set
$\Bigl(\,p_h(\cdot;c)\,\Bigr)_{h>0}$ is  bounded in $L^\infty(\Om)$.
Therefore for any sequence $h_n\to 0$, there exists a subsequence
(not relabeled) such that for all $c\in\mathbb{Q}$,
$\Bigl(\,p_h(\cdot;c)\,\Bigr)_{h>0}$ converges in $L^\infty(\Om)$ weak-$\star$ to
some limit denoted by $p(\cdot;c)$.

Fix $\xi\in \mathcal D(\Om)$, $\xi\geq 0$. From Definition~\ref{def:Entropy-Solutionp}, taking in
\ref{def:entropyWeakp} test functions approaching
$\psi(t,x):=\Bigl(1\!-\!\frac t{h_n}\Bigr)^+\,\xi(x)$ we readily
infer the inequalities
\begin{equation}\label{eq:weak-trace-ineq}
    \forall c\in \mathbb{Q} \quad
    \int_\Om p(x;c)\,\xi(x)\,\dx
    \leq \int_\Om |u_0(x)-c|\,\xi(x)\,\dx.
\end{equation}
By the density argument, we extend \eqref{eq:weak-trace-ineq} to
all $\xi\in L^1(\Om)$, $\xi\geq 0$.

Now for all $\delta>0$, there exists a number $N(\delta)\in
\mathbb{N}$, a collection $(c_i^\delta)_{i=1}^{N(\delta)}\subset
\mathbb{Q}$ and a partition of $\Om$ into disjoint union of
measurable sets $\Omega^\delta_1,\dots,\Omega^\delta_{N(\delta)}$
such that $\|u_0-u_0^\delta\|_{L^1}\leq \delta$, where
$$
u_0^\delta:=\sum\nolimits_{i=1}^{N(\delta)}
c_i^\delta\,\char_{\Om^\delta_i}.
$$
Because $\char_{\Om}=\sum\nolimits_{i=1}^{N(\delta)}
\char_{\Om^\delta_i}$, applying \eqref{eq:weak-trace-ineq} with
$c=c^\delta_i$ and $\xi=\char_{\Om^\delta_i}$ we deduce
\begin{align*}
    & \lim_{n\to\infty} \frac 1{h_n}
    \int_0^{h_n} \int_\Om \int_0^1
    \biggl(\,\frac{1}{d}|\mu-u^\delta_0|
    +\frac{d-1}{d}|\mu^\star-u^\delta_0|\biggr)
    \dt\dx\dalpha \\
    & \qquad
    = \lim_{n\to\infty} \int_\Om
    \sum\nolimits_{i=1}^{N(\delta)}
    p_{h_n}(x;c^\delta_i)\char_{\Om^\delta_i}(x)\dx
    = \int_\Om \sum\nolimits_{i=1}^{N(\delta)}
    p(x;c^\delta_i)\char_{\Om^\delta_i}(x)\dx \\
    & \qquad \leq
    \int_\Om \sum\nolimits_{i=1}^{N(\delta)}
    |u_0(x)-c^\delta_i|\,\char_{\Om^\delta_i}(x)\dx
    =\norm{u_0-u^\delta_0}_{L^1}\;\leq \delta.
\end{align*}
Using once more the bound $\|u_0-u^\delta_0\|_{L^1}\leq \delta$ (in
the first term of the previous calculation), we can
send $\delta$ to zero and infer
the analogue of \eqref{eq:IC-strong}, with a limit
taken on some subsequence of $(h_n)_{n>1}$.
Because $(h_n)_{n>1}$ was an arbitrary sequence
convergent to zero, \eqref{eq:IC-strong} is justified.
\end{proof}

Lemma \ref{lem:unique1} tells us that ${\mathcal I}[\cdot]$, which is
a distribution on $[0,T)\times \overline{\Omega}$, is
nonnegative when restricted to $\D([0,T)\times\Omega)$ and thus
it is a  locally finite measure on $[0,T)\times\Omega$.

Now we show that ${\mathcal I}[\phi]$ is nonnegative also for nonnegative
test functions $\phi$ that do not necessarily vanish on the boundary $[0,T)\times \ptl\Om$.

\begin{lem}\label{lem:unique2}
Let $(\nu,\nu,v),(\mu,\mu^\star,w)$ be entropy double-process
solutions of \eqref{eq:prob1} with initial data $v_0,u_0\in L^\infty(\Omega)$, respectively.
Then ${\mathcal I}[\phi]\geq 0$, $\forall \phi \in \D([0,T)\times \oOm)$, $\phi\ge 0$.
\end{lem}
\begin{proof}
We begin by modifying steps {\bf (i)-(iv)} of the proof of the previous lemma;
we refer to this proof for the notation and a part of the calculations.

{\bf (i)} We use \eqref{eq:gen-entr-proc} and
\ref{def:entropy2p} in the same way as in the proof
of Lemma \ref{lem:unique1}; but we choose the
values $c=\mu^+(s,y,\beta)$, $c=(\mu^\star)^+(s,y,\beta)$ and
$W=\mathfrak{a}(\Grad w^+(s,y))$ instead of the values
$c=\mu(s,y,\beta)$, $c=\mu^\star(s,y,\beta)$ and
$W=\mathfrak{a}(\Grad w(s,y))$, respectively.

Notice that for all $a,b\in\R$, $\eps\geq 0$, we have
$\sign^+_\eps(a-b^+)=\sign^+_\eps(a^+-b^+)$ and, moreover, this
expression is zero whenever $a\leq 0$. Thus we can replace
$\nu$,$v$,$\grad v$ by $\nu^+$,$v^+$,$\grad v^+$ everywhere in
this calculation and obtain
\begin{equation}\label{eq:Part-1-ofComparison}
     \begin{split}
        \!\!\!\!\!\!\!\!\!\!\!\! &\int_0^1\!\!\int_0^1 \iint_{Q\times Q}
        \Biggr[\frac{1}{d}\Bigl(\,(\nu^+\!-\!\mu^+)^++(d-1)(\nu^+\!-\!(\mu^\star)^+)^+\Bigr)\;\pt\phi
        \\
        \!\!\!\!\!\!\!\!\!\!\!\! &\qquad\qquad\qquad\qquad
        + \frac{1}{d}\Bigl(\;\sgn{\nu^+\!-\!\mu^+}(\mathfrak{f}(\nu^+)\!-\!\mathfrak{f}(\mu^+))
        \\
        \!\!\!\!\!\!\!\!\!\!\!\!&\qquad \qquad\qquad\qquad
        +\;(d-1)\sgn{\nu^+\!-\!(\mu^\star)^+}(\mathfrak{f}(\nu^+)\!-\!\mathfrak{f}((\mu^\star)^+))\;\Bigl)
        \cdot\Grad_x \phi
        \\
        \!\!\!\!\!\!\!\!\!\!\!\!& \qquad\qquad\qquad
        - \sgn{v^+\!-\!w^+}\Bigl(\mathfrak{a}(\Grad v^+)\!
        -\!\mathfrak{a}(\Grad w^+)\Bigr) \cdot \;\Grad_x\phi \Biggl]\dx\dt \dy\ds\dalpha\dbeta
        \\
        \!\!\!\!\!\!\!\!\!\!\!\!& \qquad
        +\int_0^1\iint_{\Omega\times Q}
        \frac{1}{d}\Bigl(\,(v_0^+\!-\!\mu^+)^++(d-1)
        (v_0^+\!-\!(\mu^\star)^+)^+\Bigr)\,\phi\dx (dyds)\dbeta
        \\
        \!\!\!\!\!\!\!\!\!\!\!\!& \qquad \geq
        \lim_{\eps\dto 0} \iint_{\E_\nu^c\times \E_\mu^c}
        \dsgne{v^+\!-\!w^+}\, \Bigl(\mathfrak{a}(\Grad v^+)\!-\!
        \mathfrak{a}(\Grad w^+)\Bigr)\cdot \Grad v^+\; \phi \dx\dt\dy\ds.
   \end{split}\!\!\!\!\!\!
\end{equation}

{\bf (ii)} We follow the proof of Lemma \ref{lem:unique1} but
choose $c=\nu^+(t,x,\alpha)$,$W=\mathfrak{a}(\Grad v^+(t,x))$
instead of $c=\nu(t,x,\alpha)$,$W=\mathfrak{a}(\Grad w(t,x))$.

Let us apply identities \eqref{eq:SplittingIntoPlusMinusParts} to
$a=\nu$, $b=\mu$ (or $b=\mu^*$) with $F=\mathrm{Id}$ and with
$F=\mathfrak{f}_i$, $i=1,\dots,d$. Moreover, as in the proof of \eqref{eq:I-split}, we
also have the analogue of \eqref{eq:SplittingIntoPlusMinusParts}$(i)$
with $a=v$, $b=w$ with $F(w)$ replaced
by $\mathfrak{a}(\grad w)$. In the same way, we also have the analogue of
\eqref{eq:SplittingIntoPlusMinusParts}$(ii)$ with $a=v$, $b=w$, and $F(v)$
replaced by $\mathfrak{a}(\grad v)$. Furthermore,
$\sign^\pm(\cdot)$ can be replaced by  $(\sign^\pm_\eps)'(\cdot)$
in the above properties. In conclusion, we obtain
\begin{equation}
\label{eq:Part-2-ofComparison}
  \!\!\!\!\!\!\!\!  \begin{split}
        \!\!\!\!\! &\int_0^1\!\!\int_0^1 \iint_{Q\times Q}
        \Biggr[\frac{1}{d}\Bigl(\,(\nu^+\!-\!\mu^+)^++(d-1)(\nu^+\!-\!(\mu^\star)^+)^+\Bigr)\;\ps\phi
        \\
         \!\!\!\!\! &\qquad\qquad\qquad
        + \frac{1}{d}\Bigl(\;\sgn{\nu^+\!-\!\mu^+}(\mathfrak{f}(\nu^+)\!-\!\mathfrak{f}(\mu^+))
        \\
         \!\!\!\!\!\!\!\! &\qquad \qquad\qquad\qquad
        +\;(d-1)\sgn{\nu^+\!-\!(\mu^\star)^+}(\mathfrak{f}(\nu^+)\!-\!\mathfrak{f}((\mu^\star)^+))\;\Bigl)
        \cdot\Grad_y \phi
        \\
         \!\!\!\!\! & \qquad\qquad\qquad
        - \sgn{v^+\!-\!w^+}\Bigl(\mathfrak{a}(\Grad v^+)\!
        -\!\mathfrak{a}(\Grad w^+)\Bigr) \cdot \;\Grad_y\phi \Biggl]\dx\dt\dy\ds\dalpha\dbeta
        \\
         \!\!\!\!\! &\qquad +\int_0^1\!\! \int_{Q\times \Omega} (\nu^+\!-\!u_0^+)^+\,\phi\;(\dx\dt)\dy\dalpha
        \\
         \!\!\!\!\!& \qquad \geq \quad \lim_{\eps\dto 0}
        \iint_{\E_\nu^c\times \E_\mu^c}
        \dsgne{v^+\!-\!w^+}\,
        \Bigl(\mathfrak{a}(\Grad w^+)\!-\! \mathfrak{a}(\Grad v^+)\Bigr)
        \cdot \Grad w^+\;\phi\dx\dt\dy\ds
        \\
         \!\!\!\!\! & \qquad+ \lim_{\eps\dto 0}
        \iint_{\E_\nu^c\times \E_\mu^c}
        \dsgne{-w}\,
        \mathfrak{a}(\Grad w)\cdot \Grad w\;
        \phi\dx\dt\dy\ds
        \\
         \!\!\!\!\! & \qquad -\int_0^1\!\!\int_0^1 \iint_{Q\times Q}
        \frac{1}{d}\Biggl[\sign^-(\mu)\Bigl\{\mu\,\ps\phi
        +\Bigl(\mathfrak{f}(\mu)\!-\!\mathfrak{a}(\grad w)\Bigr)\cdot\grad_y \phi\Bigr\}\;
        \\
         \!\!\!\!\! & \qquad
        \qquad +\;(d-1)
        \sign^-(\mu^*)\Bigl\{\mu^*\,\ps\phi
        +\Bigl(\mathfrak{f}(\mu^*)\!-\!\mathfrak{a}(\grad w)\Bigr)
        \cdot\grad_y \phi\Bigr\}\Biggr] \dx\dt\dy\ds\dalpha\dbeta
        \\
         \!\!\!\!\! &\qquad
        -\int_{Q\times\Omega} (\!u_0)^-\,\phi\;(\dx\dt)\dy.
    \end{split} \!\!\!\!\!\!\!\!\!\!\!
\end{equation}
Notice that the sum of the last two terms in
\eqref{eq:Part-2-ofComparison} can be rewritten under the form
$-{\mathcal L}_{\mu,\mu^*}(\chi)$, where
\begin{equation}\label{eq:DefChiFromPhi}
    \chi(s,y)=\int_Q \phi(t,s,x,y)\dt\dx
    \in {\mathcal D}([0,T)\times \overline{\Omega}),
\end{equation}
and the distribution ${\mathcal L}_{\mu,\mu^*}$ is defined on
${\mathcal D}([0,T)\times \overline{\Omega})$ by
\begin{equation}\label{eq:defLrond}
    \begin{split}
    {\mathcal L}_{\mu,\mu^*}(\chi)
    & :=\int_0^1\!\!\int_Q \Biggl[\frac{1}{d}
    \Bigl(\eta_0^-(\mu)\!+\!(d-1)\eta_0^-(\mu^*)\Bigr)\,\ps\chi
    \\ & \qquad
    +\frac{1}{d}\Bigl( {\mathfrak q}_0^-(\mu)\!+\!(d-1){\mathfrak q}_0^-(\mu^*)\Bigr)
    \cdot \grad_y \chi \\
    & \qquad
    - k(\Grad w) \Grad \widetilde A_{(\eta_0^-)'}(w)
    \cdot\grad_y \chi \Biggr]\dy\ds\dbeta
    +\int_{\Omega} \eta_0^-(u_0)\,\chi\dy.
    \end{split}
\end{equation}

{\bf (iii)} Adding \eqref{eq:Part-1-ofComparison} and \eqref{eq:Part-2-ofComparison}, we obtain,
for any $0\le \phi\in {\mathcal D}\left(([0,T)\times \overline{\Om})^2\right)$
with corresponding $\chi$ defined in \eqref{eq:DefChiFromPhi}, the following inequality:
\begin{equation}\label{eq:ComparPositiveParts}
   \!\!\!\!\!\!\!\!  \begin{split}
  \!\!\!\!\!\!\!\!\!\!\!\!\!       &\int_0^1\!\!\int_0^1 \iint_{Q\times Q}
        \Biggr[\frac{1}{d}\Bigl(\,(\nu^+\!-\!\mu^+)^+
        +(d-1)(\nu^+\!-\!(\mu^\star)^+)^+\Bigr)\;(\pt+\ps)\phi
        \\
         \!\!\!\!\!\!\!\!\!\!\!\!\! &\quad\qquad\qquad
        + \frac{1}{d}\Bigl(\;\sgn{\nu^+\!-\!\mu^+}(\mathfrak{f}(\nu^+)\!-\!\mathfrak{f}(\mu^+))
        \\
         \!\!\!\!\!\!\!\!\!\!\!\!\! & \quad \qquad \qquad\qquad +\;(d-1)\sgn{\nu^+\!-\!(\mu^\star)^+}
        (\mathfrak{f}(\nu^+)\!-\!\mathfrak{f}((\mu^\star)^+))\;\Bigl)
        \cdot(\Grad_x\!+\!\Grad_y) \phi
        \\
         \!\!\!\!\!\!\!\!\!\!\!\!\! & \quad\qquad\qquad
        - \sgn{v^+\!-\!w^+}\Bigl(\mathfrak{a}(\Grad v^+)\!
        -\!\mathfrak{a}(\Grad w^+)\Bigr) \cdot \;(\Grad_x\!+\!\Grad_y)\phi \Biggl]\dx\dt\dy\ds\dalpha\dbeta
        \\
         \!\!\!\!\!\!\!\!\!\!\!\!\! & \quad \quad
        +\int_0^1\iint_{\Omega\times Q}
        \frac{1}{d}\Bigl(\,(v_0^+-\mu^+)^++(d-1)(v_0^+-(\mu^\star)^+)^+\Bigr)\,\phi\dx\,(\dy\ds)\dbeta
        \\
         \!\!\!\!\!\!\!\!\!\!\!\!\!&\quad\quad
        +\int_0^1 \int_{Q\times \Omega}
        (\nu^+\!-\!u_0^+)^+\,\phi\;(\dx\dt)\dy\dalpha
        \\
         \!\!\!\!\!\!\!\!\!\!\!\!\! & \quad \quad\geq \;  \lim_{\eps\dto 0}
        \iint_{\E_\nu^c\times \E_\mu^c}\!\!\!\!\!\!
        \dsgne{v^+\!\!-\!w^+}\,
        \Bigl(\mathfrak{a}(\Grad v^+)\!-\! \mathfrak{a}(\Grad w^+)\Bigr)\cdot (\Grad v^+\!\!-\!\Grad w^+)\;
        \phi\dx\dt\dy\ds
        \\
         \!\!\!\!\!\!\!\!\!\!\!\!\! & \qquad +\;\lim_{\eps\dto 0}
        \iint_{\E_\nu^c\times \E_\mu^c}\!\!\!\!\!\!
        \dsgne{-w}\,
        \mathfrak{a}(\Grad w)\cdot \Grad w\; \phi\dx\dt\dy\ds
        -{\mathcal L}_{\mu,\mu^*}(\chi)
        \geq -{\mathcal L}_{\mu,\mu^*}(\chi),
   \end{split} \!\!\!\!\!\!\!\! \!\!\!\!\!\!\!\! \!\!\!\!\!
\end{equation}
where the last inequality is due to the monotonicity of $\mathfrak{a}(\cdot)$.

{\bf (iv)} Now fix $x_0\in\ptl\Om$. Since $\ptl\Om$ is supposed sufficiently regular,
there exists a vector $r_{x_0}$ and a positive number $R_{x_0}$ such that the segment
$(x,x+r_{x_0}]$ lies within $\Om$ for all  $x\in \ptl\Om\cap
B(x_0,R_{x_0})$, where $B(x,R)$ stands for the ball of $\R^d$ with
centre $x$ and radius $R$. Choose in \eqref{eq:ComparPositiveParts}
the sequence of test functions
\begin{equation*}
    \phi_{n,l}(t,x,s,y)=\phi\left(y,s\right)
    \omega_n\left({x-y+\frac{2}{n}\frac{r_{x_0}}{\|r_{x_0}\|}}\right)
    \omega_l\left({t-s}\right) \equiv \phi \omega_n\omega^l,
\end{equation*}
for which \eqref{bis:eq:test_tmp1} still holds. Notice that with
this choice, the associated function $\chi_{n,l}(y,s)$ in
\eqref{eq:DefChiFromPhi} writes as $\phi(y,s)
\theta_n(y)\theta^l(s)$, where
$$
\theta_n(y):=\int_\Om
\omega_n\left(x-y+\frac{2}{n}\frac{r_{x_0}}{\|r_{x_0}\|}\right)\;dx,
\qquad \theta^l(s):=\int_0^T \omega_l(t-s)\;dt;
$$
moreover, for all sufficiently large $n\in\N$ we have
\begin{equation}\label{eq:PropOfTheta-n}
\left|\begin{array}{l} \text{ $\phi\theta_n \in \D(\Om)$ for all
$\phi\in \D\biggl([0,T)\times\Bigl(\overline{\Om}\cap
B(x_0,R_{x_0}) \Bigr) \biggr)$};\\
\text{$\theta_n(y)=1$ for all $y\in B(x_0,R_{x_0})$ such that
$\dist (y, \ptl\Om)\geq \frac{3}{n}$}.
\end{array}\right.
\end{equation}

As in the proof of Lemma~\ref{lem:unique1}, passing to the limit
as $l,n\to\infty$ and taking into account the definition of the
distribution $\mathcal {IP}$, cf~\eqref{eq:I-split},
from \eqref{eq:ComparPositiveParts} we deduce
\begin{equation}\label{eq:PositiveHalf}
    \begin{split}
        &\mathcal {IP}[\phi]
        \geq\; -\liminf_{l,n\to\infty}
        {\mathcal L}_{\mu,\mu^*}[\phi\theta_n\theta_l].
    \end{split}
\end{equation}
Now we remark that according to \ref{def:entropy2p}, ${\mathcal L}_{\mu,\mu^*}$
defined by \eqref{eq:defLrond} is a nonnegative
distribution on $[0,T)\times\overline{\Om}$. Notice that the
values of $\theta_n$ are contained in the interval $[0,1]$. Therefore for
all $\chi\in\D([0,T)\times\overline{\Om})$, $\phi\geq 0$, one has
$$
{\mathcal L}_{\mu,\mu^*}[\phi\theta_n]
={\mathcal L}_{\mu,\mu^*}[\phi]
-{\mathcal L}_{\mu,\mu^*}
[\phi(1-\theta_n)]\leq {\mathcal L}_{\mu,\mu^*}[\phi].
$$
It follows that
\begin{equation}\label{eq:LimitMeasure}
\overline{\mathcal L}_{\mu,\mu^*}:\chi\in
\D([0,T)\times\overline{\Om})\mapsto \liminf_{n\to\infty}
{\mathcal L}_{\mu,\mu^*}(\chi\theta_n)
\end{equation}
is a nonnegative distribution on $[0,T)\times\overline{\Om}$;
thus, it is a measure on $[0,T)\times\overline{\Om}$. Since $\phi
\geq 0$ and $\theta^l \leq 1$, inequality \eqref{eq:PositiveHalf} yields
\begin{equation}\label{eq:MeasureIneq}
\mathcal {IP}[\phi]\geq -\liminf_{l,n\to\infty} {\mathcal L}_{\mu,\mu^*}
[\chi\theta_n\theta_l] \geq -\liminf_{n\to\infty}
{\mathcal L}_{\mu,\mu^*}[\phi\theta_n]=- \overline{\mathcal L}_{\mu,\mu^*}[\phi].
\end{equation}
It follows that $\mathcal {IP}$ is a measure on $[0,T)\times\overline{\Om}$.

The remaining steps of the proof are aimed at showing, in an indirect
way, that the positive part of the measure $\mathcal{IP}$ does not
charge the boundary $[0,T)\times \ptl\Om$ (in two particular
cases, a direct proof of this fact is given in
\cite{RouvreGagneux,AndrIgbida}). Notice that this property is
actually equivalent to the claim of the lemma; it accounts for the
dissipative nature of the boundary condition imposed for entropy solutions.

{\bf (v)} Take $\phi\in \D\biggl([0,T)\times\Bigl(\overline{\Om}
\cap B(x_0,R_{x_0}) \Bigr)\biggr)$.
Fix $m\in \N$. It is easily checked from \eqref{eq:PropOfTheta-n}
that for all sufficiently large $n\in\N$, for all $(t,x)\in Q$,
$$
\phi(s,y)(1-\theta_m(y))\theta_n(y)
=\phi(s,y)\theta_n(y)-\phi(s,y)\theta_m(y).
$$
Therefore by \eqref{eq:LimitMeasure},
\begin{align*}
    &\dsp\liminf_{m\to\infty} \overline{\mathcal L}_{\mu,\mu^*}[\phi(1-\theta_m)]
    \\ & \quad
    =\liminf_{m\to\infty}\liminf_{n\to\infty} {\mathcal L}_{\mu,\mu^*}[\phi\theta_n]
    -\liminf_{m\to\infty}\liminf_{n\to\infty}{\mathcal L}_{\mu,\mu^*}[\phi\theta_m]=0.
\end{align*}
Applying \eqref{eq:MeasureIneq} to the test function
$\phi(1-\theta_m)$, we deduce
\begin{equation}\label{eq:MeasureLbarVanishes}
    \begin{split}
        \mathcal {IP}^{\nu,\nu,v}_{\mu,\mu^*\!,w}[\phi] &
        =\mathcal{IP}^{\nu,\nu,v}_{\mu,\mu^*\!,w}[\phi\theta_m]
        +\mathcal{IP}^{\nu,\nu,v}_{\mu,\mu^*\!,w}[\phi(1-\theta_m)]
        \\ & \geq\limsup_{m\to\infty}
        \mathcal{IP}^{\nu,\nu,v}_{\mu,\mu^*\!,w}[\phi\theta_m].
    \end{split}
\end{equation}

{\bf (vi)} Definition \ref{def:Entropy-Solutionp} of entropy double-process solution
is invariant under the change of
$(\mu,\mu^*,w)$,$\mathfrak{f}$,$f$,$u_0$ into
$(-\mu,-\mu^*,-w)$,$-\mathfrak{f}$,$-f$,$-u_0$. Moreover, one
checks easily from the definition, cf.~\eqref{eq:I-split},
that $\mathcal{ IN}^{\nu,\nu,v}_{\mu,\mu^*\!,w}
=\mathcal{IP}^{-\mu,-\mu^*\!,-w}_{-\nu,-\nu,-v}$.
Therefore from \eqref{eq:MeasureLbarVanishes} we  deduce that for all
$\phi\in \D\biggl([0,T)\times\Bigl(\overline{\Om}\cap B(x_0,R_{x_0}) \Bigr)\biggr)$
$$
\begin{array}{l}
\dsp {\mathcal I}^{\nu,\nu,v}_{\mu,\mu^*\!,w}[\phi]=\mathcal
{IP}^{\nu,\nu,v}_{\mu,\mu^*\!,w}[\phi]+\mathcal
{IN}^{\nu,\nu,v}_{\mu,\mu^*\!,w}[\phi]= \mathcal
{IP}^{\nu,\nu,v}_{\mu,\mu^*\!,w}[\phi]+\mathcal
{IP}^{-\mu,-\mu^*\!,-w}_{-\nu,-\nu,-v}[\phi] \\[9pt]
\dsp\quad \geq \limsup_{m\to\infty}\biggl[\mathcal
{IP}^{\nu,\nu,v}_{\mu,\mu^*\!,w}[\phi\theta_m]+ \mathcal
{IP}^{-\mu,-\mu^*\!,-w}_{-\nu,-\nu,-v}[\phi\theta_m]\biggl]=
 \limsup_{m\to\infty}\mathcal
{I}^{\nu,\nu,v}_{\mu,\mu^*\!,w}[\phi\theta_m]\geq 0,
\end{array}
$$
where the last inequality is due to \eqref{eq:PropOfTheta-n}
and Lemma \ref{lem:unique1}.

{\bf (vii)} Not let $\phi$ be an arbitrary nonnegative
function in $\D([0,T)\times\overline\Om)$. Choose a  covering
$\bigcup_{i=1}^N B(x^i_0,R_{x^i_0})$, $N\in \N$, of the compact
set $\ptl\Om$. Introduce a partition of unity $(\xi_i)_{i=0}^N$ on
$\overline\Om$ associated with the covering  $ \Om \bigcup
\Bigl(\bigcup^N_{i=1} B(x^i_0,R_{x^i_0})\Bigr)$ of $\overline\Om$,
and apply Lemma \ref{lem:unique1} and the result of {\bf (vi)} to
the functions $\phi\xi_0\in \D([0,T)\times\Om)$ and to
$\phi\xi_i\in \D\biggl([0,T)\times\Bigl(\overline{\Om}\cap
B(x^i_0,R_{x^i_0}) \Bigr) \biggr)$, $i=1,\dots,N$, respectively.
The claim of the lemma follows.
\end{proof}

Now we conclude the proof of Theorem \ref{th:EPSunique}. We have
$u_0=v_0$. By a standard argument, choosing in
Lemma \ref{lem:unique2} $\phi=\phi(t)\in \D([0,T))$, we get for a.e.~$t\in (0,T)$,
\begin{equation}\label{L1-Contraction}
    \int_0^1\!\!\int_0^1\!\!\int_\Om
    \frac{1}{d}\Bigl(\,(\nu(t,x,\alpha)-\mu(t,x,\beta))^+
    +(d-1)(\nu(t,x,\alpha)-\mu^\star(t,x,\beta))^+\Bigr)\le 0.
\end{equation}
Now, \eqref{L1-Contraction} means that for a.e.~$(x,\alpha,\beta)\in\Omega\times(0,1)\times(0,1)$,
there holds
$$
\mu(t,x,\beta)=\nu(t,x,\alpha)=\mu^*(t,x,\beta),
$$
which means that $\mu\equiv \mu^* \equiv \nu$ and each of them is independent
of $\alpha,\beta$. This draws to a close the proof of Theorem \ref{th:EPSunique}.

\begin{rem}\normalfont \label{rem:L1CCPrinciple}
The proof of the $L^1$ contraction and comparison principle for
entropy solutions of \eqref{eq:prob1} (with ${\EuScript S}=0$) is
essentially contained in the above proof. For nonzero source terms
${\EuScript S}$, a more general version of inequalities
\eqref{eq:SplittingIntoPlusMinusParts} can be used; see
\cite{Carrillo} for the accurate treatment of this term.
\end{rem}

\section*{Appendix B: The reconstruction property}
Here we restate the result  of \cite[Lemma 8]{AGW} and discuss its
possible generalizations.

\begin{lem}\label{lem:ReconstrPropTriangles}
Consider a triangle $T\!\!\!T$ with vertices $t_1,t_2,t_3$ and let
$t_0$ be the centre of its circumscribed circle. Denote by
$|T\!\!\!T|$ its area. For $l\in \N/3\N$, denote by
$E_l$ the affine subspace $<\overrightarrow{t_{l-1} t_{l+1}}>$; denote by $T\!\!\!T_l$ the
triangle formed by $t_0,t_{l-1},t_{l+1}$ and by $|T\!\!\!T_l|$ its area, with the convention that
the area is negative if $t_0$ and $t_l$ lay on opposite sides
from the line passing by $t_{l-1},t_{l+1}$. Then
\begin{equation}\label{eq:ReconstrInTriangle}
    \frac{2}{|T\!\!\!T|} \sum_{l=1}^3 |T\!\!\!T_l|
    \mathrm{Proj}_{E_l}(\overrightarrow r)
    =\overrightarrow r, \qquad
    \text{for all $\overrightarrow r\in\R^2$}.
\end{equation}
\end{lem}

\begin{rem}\normalfont \label{rem:ReconstrGeneralizationMultiD}
For a multi-D generalization of the property
\eqref{eq:ReconstrInTriangle}, one could try to replace the
projections on lines $<E_l>$ by projections on hyperplanes that
contain the faces of the $d$-dimensional simplex $T\!\!\!T$. In
this case one should replace the factor $\frac{2}{|T\!\!\!T|}$ by
$\frac{d}{d-1}\frac{1}{|T\!\!\!T|}$, since
$|T\!\!\!T|=\sum_{l=1}^{d+1}|T\!\!\!T_l|$ and because the
dimension of $\mathrm{Proj}_{E_l}(\overrightarrow r)$ is
$(d-1)$, whereas the dimension of $\overrightarrow r$ is $d$.
The proof of Lemma~\ref{lem:ReconstrPropTriangles} given below
shows that this generalization fails, except for very particular
simplexes $T\!\!\!T$ (this is clear from the multi-dimensional
analogue of the identity \eqref{eq:ReconstrInTriangEquiv} below).
\end{rem}

\begin{rem}\normalfont \label{rem:ReconstrGeneralization2DPolygon}
Using the ``sine theorem'', another proof of
Lemma~\ref{lem:ReconstrPropTriangles} can be given,  which also
works for any 2D polygon that admits a circumscribed circle.
\end{rem}

\begin{proof}{Proof of Lemma \ref{lem:ReconstrPropTriangles}}
For $l\in \N/3\N$, denote by $d_l$ the
orthogonal projection if the point $t_l$ on the
affine subspace $E_l$; set $\overrightarrow
p_l=\overrightarrow{t_0 d_l}$ and $\overrightarrow
a_l=\overrightarrow{t_0 t_l}$. For $l,i\in \N/3\N$, set
$\overrightarrow b_{l,i}=\overrightarrow a_i-\overrightarrow a_l$.
Denote by $\overrightarrow n_l$ the exterior to
$T\!\!\!T$ unit normal vector to $E_l$. Notice that we have for
all $l\in \N/3\N$, $\overrightarrow d_l=(\overrightarrow
d_l\cdot \overrightarrow n_l)\,\overrightarrow n_l$,
and also, for all $i\in \N/3\N$ such that $i\neq l$,
$$
\frac{|T\!\!\!T_l|}{|T\!\!\!T|}=\frac{\overrightarrow p_l\cdot
\overrightarrow n_l}{\overrightarrow b_{l,i}\cdot \overrightarrow n_l},
$$
taking into account the sign of $|T\!\!\!T_l|$.
Since $\mathrm{Proj}_{E_l}+\mathrm{Proj}_{<\overrightarrow n_l>}$
is the identity operator, \eqref{eq:ReconstrInTriangle}
is equivalent to the statement that
$\frac{2}{|T\!\!\!T|} \sum_{l=1}^3 |T\!\!\!T_l|
\mathrm{Proj}_{<\overrightarrow n_l>}$
is the identity operator. All vector $\overrightarrow r\in \R^2$
can be uniquely represented under the form
$$
\overrightarrow r=\sum_{l=1}^{3} k_l\,\overrightarrow a_l \qquad
\text{with} \quad \sum_{l=1}^3 k_l=0,
$$
and thus, for all $l\in \N/3\N$,
$\overrightarrow r=\sum_{i\neq l,i=1}^{3}
k_i\,\overrightarrow  (a_i-\overrightarrow a_l)
=\sum_{i\neq l,i=1}^{3} k_i\,\overrightarrow b_{l,i}$.
Hence
$$
\begin{array}{l}
\dsp \frac{2}{|T\!\!\!T|} \sum_{l=1}^3 |T\!\!\!T_l|
\mathrm{Proj}_{<\overrightarrow n_l>}(\overrightarrow r)=2
\sum_{l=1}^3 \frac{|T\!\!\!T_l|}{|T\!\!\!T|}
\,\overrightarrow n_l\,
(\overrightarrow r \cdot \overrightarrow n_l) = 2
\!\!\!\!\sum_{i\neq l;\;i,l=1}^3\!\!\!\! \overrightarrow n_l
\frac{|T\!\!\!T_l|}{|T\!\!\!T|} \,k_i (\overrightarrow b_{l,i}\cdot \overrightarrow n_l)\\[10pt]
\dsp=  2 \sum_{l=1}^3  \overrightarrow n_l  \!\!\! \sum_{i\neq l,i=1}^3
\frac{\overrightarrow p_l\cdot \overrightarrow n_l}{\overrightarrow b_{l,i}\cdot \overrightarrow n_l}\, k_i
(\overrightarrow b_{l,i}\cdot \overrightarrow n_l) =
2\sum_{l=1}^3 \overrightarrow p_l  \!\!\!\sum_{i\neq l,i=1}^3 k_i
=- 2\sum_{l=1}^3 k_l \overrightarrow p_l.
\end{array}
$$
We conclude that \eqref{eq:ReconstrInTriangle} is equivalent to the identity
\begin{equation}\label{eq:ReconstrInTriangEquiv}
\sum_{l=1}^{3} k_l\,\overrightarrow a_l =- 2\sum_{l=1}^3 k_l
\overrightarrow p_l \qquad \text{for all $k_1,\dots,k_3\in\R$ such
that $\sum_{l=1}^3 k_l=0$}.
\end{equation}
Since $t_0$ is the centre of the circumscribed circle of
$T\!\!\!T$, the points $d_l$ are the centres of the corresponding
segments $[t_{i-1},t_{i+1}]$. Thus for all $i,j\in \N/3\N$, by the
Thales theorem we have $\overrightarrow p_i-\overrightarrow
p_j=-\frac{1}{2}(\overrightarrow a_i-\overrightarrow a_j)$. Hence
\eqref{eq:ReconstrInTriangEquiv} holds with $k_i\in\{0,1,-1\}$,
$i=1,\dots,3$. Hence it holds for all choice of $k_i$.
\end{proof}

We refer to \cite{AndrBend,ABK-FVCA5} for a different kind of
generalization of \cite[Lemma 8]{AGW} and a different proof of
Lemma~\ref{lem:ReconstrPropTriangles}.


\begin{thebibliography}{10}

\bibitem{Aavatsmark} I.~Aavatsmark.
\newblock{An introduction to multipoint flux approximations for
quadrilateral grids. Locally conservative numerical methods for
flow in porous media.},  \newblock{\em Comput. Geosci.} 6(3--4):
405--432, 2002.

\bibitem{Aavatsmark-et-al} I.~Aavatsmark, T.~Barkve, \O.~B\oe, and
T.~Mannseth.
\newblock{
Discretization on unstructured grids for inhomogeneous,
anisotropic media. Part I: derivation of the methods.}
\newblock {\em SIAM J.Sci.Comp.} 19(5):1700--1716, 1998.

\bibitem{AfifAmaziane} M.~Afif and B.~Amaziane.
\newblock{Convergence of finite volume schemes for a degenerate
convection-diffusion equation arising in flows in porous media},
\newblock{\em Comput. Methods Appl. Mech. Eng.} 191(46):5265--5286, 2002.

\bibitem{AltLuckhaus}
H.~W. Alt and S.~Luckhaus.
\newblock Quasilinear elliptic-parabolic differential equations.
\newblock {\em Math. Z.} 183(3):311--341, 1983.

\bibitem{AmmarCarWitt}
K.~Ammar, P.~Wittbold and J.~Carrillo.
\newblock Scalar
conservation laws with general boundary condition and continuous
flux function. \newblock {\em J. Diff. Eq.}  228(1):111--139,
2006.

\bibitem{AmmarWittbold}
K.~Ammar and P.~Wittbold. \newblock Existence of renormalized
solutions of degenerate elliptic-parabolic problems. \newblock
{\em Proc. Roy. Soc. Edinburgh Sect. A}  133(3):477--496, 2003.

\bibitem{AndrBend} B.~Andreianov and M.~Bendahmane.
\newblock{On Discrete Duality Finite Volume discretization\\ of
gradient and divergence operators in 3D}, preprint

\bibitem{ABK-FVCA5} B.~Andreianov,
M.~Bendahmane and K.H.~Karlsen.
\newblock{A gradient reconstruction formula for finite volume schemes and
discrete duality},  in  R. Eymard and J. M. Herard, eds.,
\newblock{Proceedings of Finite Volumes for Complex Applications V},
Herm\`es, 2008.

\bibitem{ABKO} B.~Andreianov, M.~Bendahmane,
K.H.~Karlsen and S.~Ouaro.
\newblock{Well-posed\-ness re\-sults for triply nonlinear degenerate parabolic
equations}, preprint,\\
http://www.math.ntnu.no/conservation/2008/031.html

\bibitem{ABH-M2AN}
B.~Andreianov, F.~Boyer and F.~Hubert.
\newblock Finite volume schemes for the $p-\!\!$Laplacian on Cartesian meshes.
\newblock {\em M2AN  Math. Model.
Numer. Anal.} 38(6):931-960, 2005.

\bibitem{ABH``double''}
B.~Andreianov, F.~Boyer and F.~Hubert.
\newblock Discrete duality finite volume schemes for Leray-Lions type elliptic problems on general 2D meshes.
\newblock {\em Num. Meth. PDE} 23(1):145--195, 2007.

\bibitem{AGW}
B.~Andreianov, M.~Gutnic and P.~Wittbold.
\newblock Convergence of finite volume
 approximations for a nonlinear elliptic-parabolic problem: a "continuous" approach.
 \newblock{\em  SIAM J. Numer. Anal.}  42(1):228--251, 2004

\bibitem{AndrIgbida}
B.~Andreianov and N.~Igbida.
\newblock Uniqueness for inhomogeneous
Dirichlet problem for
 elliptic-parabolic equations.
 \newblock {\em Proc. Royal Soc. Edinburgh A} 137(6):1119--1133, 2007.

\bibitem{A-D:Nat:Tang}
D.~Aregba-Driollet, R.~Natalini and S.~Tang.
\newblock Explicit diffusive kinetic schemes for nonlinear degenerate parabolic
  systems.
\newblock {\em Math. Comp.}  73(245):63--94, 2004.

\bibitem{Bardos:BC}
C.~Bardos, A.~Y. LeRoux and J.-C. N{\'e}d{\'e}lec.
\newblock First order quasilinear equations with boundary conditions.
\newblock {\em Comm. Partial Differential Equations} 4(9):1017--1034, 1979.

\bibitem{BarrettLiu}
J.W.~Barrett and W.B.~Liu.
\newblock
 Finite element approximation
of the parabolic $p$-Laplacian. \newblock{  SIAM J. Numer. Anal.}
31(2):413--428, 1994.

\bibitem{BenKar:Renorm}
M.~Bendahmane and K.~H. Karlsen.
\newblock Renormalized entropy solutions for quasilinear anisotropic degenerate
  parabolic equations.
\newblock {\em SIAM J. Math. Anal.} 36(2):405--422, 2004.

\bibitem{BenKar:RenormII}
M.~Bendahmane and K.~H. Karlsen.
\newblock Uniqueness of entropy solutions for
doubly nonlinear anisotropic degenerate parabolic equations.
\newblock
{\em Contemporary Mathematics} 371, Amer. Math. Soc., pp.1--27,
2005.

\bibitem{BenilanToureIII}
P.~B\'enilan and H.~Tour{\'e}.
\newblock Sur l'\'equation g\'en\'erale $u\sb t=a(\cdot,u,\phi(\cdot,u)\sb
  x)\sb x+v$ dans ${L}\sp 1$. {I}{I}. {L}e probl\`eme d'\'evolution.
\newblock {\em Ann. Inst. H. Poincar\'e Anal. Non Lin\'eaire} 12(6):727--761,
  1995.

\bibitem{BenilanWittbold}
P.~B\'enilan and P.~Wittbold.
\newblock On mild and weak solutions of elliptic-parabolic problems.
\newblock {\em Adv. Differential Equations} 1(6):1053--1073, 1996.

\bibitem{BocMuratPuel}
L. Boccardo, F. Murat and J.-P. Puel.
\newblock{ Existence of bounded solutions for nonlinear elliptic
unilateral problems,} {\newblock \em Ann. Math. Pura Appl.}
152(1):183--196, 1988.


\bibitem{BouGuaNat}
F.~Bouchut, F.~R. Guarguaglini and R.~Natalini.
\newblock Diffusive {B}{G}{K} approximations for nonlinear multidimensional parabolic equations.
\newblock {\em Indiana Univ. Math. J.} 49(2):723--749, 2000.

\bibitem{AngotetCo}
F.~Boyer and F.~Hubert.
\newblock \newblock Finite volume method for 2D linear and nonlinear elliptic
problems with discontinuities. \newblock {SIAM J. Num. Anal.}
46(6):3032--3070, 2008.

\bibitem{Brow} F. Browder.
{\newblock Existence theorems for nonlinear partial differential equations,}
{\newblock in \em Proc. Symp. Pure Math.}, S.S. Chern and S. Smale,
eds., AMS, Providence, Rhode Island, 1970,  pp.1--60.

\bibitem{BurgerEvjeKarlsen:1D_IBVP}
R.~B\"urger, S.~Evje and K.~H. Karlsen.
\newblock On strongly degenerate convection-diffusion problems modeling
  sedimentation-consolidation processes.
\newblock {\em J. Math. Anal. Appl.} 247(2):517--556, 2000.

\bibitem{BurgKarlMishTow}
R.~B\"urger, K.~H. Karlsen, S.~Mishra and J.~Towers.
\newblock On conservation laws with discontinuous flux.
\newblock In: {\em Y. Wang and K. Hutter (Eds.), Trends in Applications of
Mathematics to Mechanics}, Shaker Verlag, Aachen, pp.75--84, 2005.

\bibitem{Burger:bok}
M.~C. Bustos, F.~Concha, R.~B\"{u}rger and E.~M. Tory.
\newblock {\em Sedimentation and Thickening: Phenomenological Foundation and
  Mathematical Theory}.
\newblock Kluwer Academic Publishers, Dordrecht, The Netherlands, 1999.

\bibitem{Carrillo:94}
J.~Carrillo.
\newblock On the uniqueness of the solution of the evolution dam problem.
\newblock {\em Nonlinear Anal.} 22(5):573--607, 1994.

\bibitem{Carrillo}
J.~Carrillo.
\newblock Entropy solutions for nonlinear degenerate problems.
\newblock {\em Arch. Rational Mech. Anal.} 147(4):269--361, 1999.

\bibitem{CarrilloWittbold:99}
J.~Carrillo and P.~Wittbold.
\newblock Uniqueness of renormalized solutions of degenerate elliptic-parabolic
  problems.
\newblock {\em J. Differential Equations} 156(1):93--121, 1999.

\bibitem{ClaireChainais}
C.~Chainais-Hillairet. \newblock Finite volume schemes for a
nonlinear hyperbolic equation. Convergence towards the entropy
solution and error estimate. \newblock {\em M2AN Math. Model.
Numer. Anal.} 33(1):129--156, 1999.

\bibitem{ChenDiBen}
G.-Q. Chen and E.~DiBenedetto.
\newblock Stability of entropy solutions to the {C}auchy problem for a class of
  nonlinear hyperbolic-parabolic equations.
\newblock {\em SIAM J. Math. Anal.} 33(4):751--762, 2001.

\bibitem{Chen:2005wf}
G.-Q. Chen and K.~H. Karlsen.
\newblock Quasilinear anisotropic degenerate parabolic equations with
  time-space dependent diffusion coefficients.
\newblock {\em Commun. Pure Appl. Anal.} 4(2):241--266, 2005.

\bibitem{ChenKarlsen}
G.-Q. Chen and K.~H. Karlsen.
\newblock {$L^1$} framework for continuous dependence and error estimates for
  quasi-linear degenerate parabolic equations.
\newblock {\em Trans. Amer. Math. Soc.} 358(3):937--963, 2006.

\bibitem{ChenPerthame}
G.-Q. Chen and B.~Perthame.
\newblock Well-posedness for non-isotropic degenerate hyperbolic-parabolic equations.
\newblock {\em Ann. Inst. H. Poincar\'e Anal. Non Lin\'eaire} 20(4):645--668, 2003.

\bibitem{Chow}
S.-S.~Chow.
\newblock Finite element error estimates for non-linear elliptic equations of monotone
type.
\newblock {\em Numer. Math.} 54(4):373-393, 1989.

\bibitem{CockGripen}
B.~Cockburn and G.~Gripenberg.
\newblock Continuous dependence on the nonlinearities of solutions of
  degenerate parabolic equations.
\newblock {\em J. Differential Equations} 151(2):231--251, 1999.

\bibitem{CoudiereHubertPierre}
Y.~Coudiere and F.~Hubert. \newblock{General 3D discrete duality
FV meshes.}
\newblock {\em in preparation}

\bibitem{CoudierePierre}
Y.~Coudiere and Ch.~Pierre. \newblock Stability and convergence of
a finite volume method for two systems of reaction-diffusion
equations in electro-cardiology. \newblock {\em Nonlinear Anal.
Real World Appl.} 7(4):916--935, 2006

\bibitem{CoudiereVilaVilledieu}
Y.~Coudi\`ere, J.-P.~Vila and Ph.~Villedieu. \newblock
Convergence rate of a finite volume scheme for a two-dimensional
convection-diffusion problem.  \newblock {\em M2AN Math. Model.
Numer. Anal.} 33(3):493--516, 1999.

\bibitem{DomOmnes}
K.~Domelevo and P.~Omn\`es.
\newblock A finite volume method for the Laplace equation on almost arbitrary two-dimensional grids.
\newblock {\em M2AN Math. Model. Numer. Anal.} 39(6):1203--1249, 2005.

\bibitem{Droniou}
J.~Droniou.
\newblock Finite volume approximations for fully nonlinear
elliptic equations in divergence form.
\newblock {\em M2AN Math. Model. Numer. Anal.} 40(6):1069--1100, 2006.

\bibitem{EspKar}
M.~S. Espedal and K.~H. Karlsen.
\newblock Numerical solution of reservoir flow models based on large time step
  operator splitting algorithms.
\newblock In {\em Filtration in Porous Media and Industrial Applications
  (Cetraro, Italy, 1998)}, volume 1734 of {\em Lecture Notes in Mathematics},
  pp:9--77. Springer, Berlin, 2000.

\bibitem{EvjeKarlsen:DD}
S.~Evje and K.~H. Karlsen.
\newblock Discrete approximations of {$BV$} solutions to doubly nonlinear
  degenerate parabolic equations.
\newblock {\em Numer. Math.} 86(3):377--417, 2000.

\bibitem{EvjeKarlsen:SJNA}
S.~Evje and K.~H. Karlsen.
\newblock Monotone difference approximations of {$BV$} solutions to degenerate
  convection-diffusion equations.
\newblock {\em SIAM J. Numer. Anal.} 37(6):1838--1860, 2000.

\bibitem{EKR:Hyp2000}
S.~Evje, K.~H. Karlsen and N.~H. Risebro.
\newblock A continuous dependence result for nonlinear degenerate parabolic
  equations with spatially dependent flux function.
\newblock In {\em Hyperbolic problems: theory, numerics, applications, Vol. I
  (Magdeburg, 2000)}, pp. 337--346. Birkh\"auser, Basel, 2001.

\bibitem{EymardDroniou}
R.~Eymard and J.~Droniou.
\newblock A mixed finite volume scheme for anisotropic diffusion problems on any grid.
\newblock {\em  Numer. Math.} 105(1):35--71, 2006.

\bibitem{EyGaHe:book}
R.~Eymard, T.~Gallou{\"e}t and R.~Herbin.
\newblock Finite Volume Methods. {\em Handbook of Numerical
Analysis}, Vol. VII, P.~Ciarlet, J.-L.~Lions, eds.,
\newblock North-Holland, 2000.

\bibitem{EGHMichel}
R.~Eymard, T.~Gallou{\"e}t, R.~Herbin and A.~Michel.
\newblock Convergence of a finite volume scheme for nonlinear degenerate
  parabolic equations.
\newblock {\em Numer. Math.} 92(1):41--82, 2002.


\bibitem{EGH:IMAJNA}
R.~Eymard, T.~Gallou{\"e}t and R.~Herbin. \newblock A cell-centered
finite-volume approximation for anisotropic diffusion operators on
unstructured meshes in any space dimension.  \newblock {\em IMA J.
Numer. Anal.} 26(2):326--353, 2006.

\bibitem{EGH:AnisoCRAS}
R.~Eymard, T.~Gallou{\"e}t and R.~Herbin. \newblock A new finite
volume scheme for anisotropic diffusion problems on general grids:
convergence analysis.  \newblock{\em C. R. Math. Acad. Sci. Paris}
344(6):403--406, 2006.

\bibitem{EGH:Sushi-IMAJNA}
R.~Eymard, T.~Gallou{\"e}t and R.~Herbin. \newblock Discretisation of
heterogeneous and anisotropic diffusion problems on general
non-conforming meshes. SUSHI: a scheme using stabilisation and
hybrid interfaces, \newblock{} preprint HAL,
http://hal.archives-ouvertes.fr/hal-00203269

\bibitem{EH-SUCCES}
R.~Eymard and R.~Herbin.
\newblock A new collocated FV scheme for the incompressible
Navier-Stokes equations on general non-matching grids. \newblock
{\em C. R. Math. Acad. Sci. Paris} 344(10):659-­662, 2007.

\bibitem{GaHu} T.~Gallou\"et and F.~Hubert.
\newblock On the convergence of the parabolic approximation 
of a conservation law in several space dimensions. no. 1, 141. 
\newblock{ Chinese Ann. Math. Ser. B} 20(1):7-10, 1999.

\bibitem{GuarMilTer}
F.~Guarguaglini, V.~Milisic and A.~Terracina.
\newblock A discrete BGK approximation for strongly degenerate parabolic
  problems with boundary conditions.
\newblock {\em J. Diff. Eq.}  202(2):183--207, 2004.

\bibitem{HandMikulaSgallari}
A.~Handlovi\v{c}ov\'a, K.~Mikula and F.~Sgallari.
\newblock{Semi-implicit complementary volume scheme for solving level set
like equations in image processing and curve evolution},
\newblock{\em Numer. Math.} 93(4):675--695, 2003.

\bibitem{Herm}
F.~Hermeline.
\newblock A finite volume method for the approximation of diffusion operators on distorted meshes.
\newblock {\em  J. Comput. Phys.}  160(2):481--499, 2000.

\bibitem{Herm-3D}
F.~Hermeline.
\newblock Approximation of 2D and 3D diffusion operators with discontinuous full-tensor coefficients
on arbitrary meshes.
\newblock {\em Comput. Meth. Appl. Mech. Engrg.} 196:2497-2526, 2007.

\bibitem{Igbida-Urbano}
N.~Igbida and J.~M. Urbano.
\newblock Uniqueness for nonlinear degenerate problems.
\newblock {\em NoDEA Nonlinear Differential Equations Appl.} 10(3):287--307,
  2003.

\bibitem{KO:Unique}
K.~H. Karlsen and M.~Ohlberger.
\newblock A note on the uniqueness of entropy solutions of nonlinear degenerate
  parabolic equations.
\newblock {\em J. Math. Anal. Appl.} 275(1):439--458, 2002.

\bibitem{KR:Rough_Diff}
K.~H. Karlsen and N.~H. Risebro.
\newblock Convergence of finite difference schemes for viscous and inviscid
  conservation laws with rough coefficients.
\newblock {\em M2AN Math. Model. Numer. Anal.} 35(2):239--269, 2001.

\bibitem{KR:Rough_Unique}
K.~H. Karlsen and N.~H. Risebro.
\newblock On the uniqueness and stability of entropy solutions of nonlinear
  degenerate parabolic equations with rough coefficients.
\newblock {\em Discrete Contin. Dyn. Syst.} 9(5):1081--1104, 2003.

\bibitem{Kruzkov}
S.~N. Kru{\v{z}}kov.
\newblock First order quasi-linear equations in several independent variables.
\newblock {\em Math. USSR Sbornik} 10(2):217--243, 1970.

\bibitem{Lions:Book69_Fr}
J.-L. Lions.
\newblock {\em Quelques m\'ethodes de r\'esolution des probl\`emes aux limites
  non lin\'eaires}.
\newblock Dunod, 1969.

\bibitem{Malek}
J.~M{\'a}lek, J.~Ne{\v{c}}as, M.~Rokyta and M.~R{u}{\v{z}}i{\v{c}}ka.
\newblock {\em Weak and measure-valued solutions to evolutionary {P}{D}{E}s}.
\newblock Chapman \& Hall, London, 1996.

\bibitem{Mascia_etal:2000}
C.~Mascia, A.~Porretta and A.~Terracina.
\newblock Nonhomogeneous dirichlet problems for degenerate hyperbolic-parabolic
  equations.
\newblock {\em Arch. Ration. Mech. Anal.} 163(2):87--124, 2002.

\bibitem{MichelVovelle}
A.~Michel and J.~Vovelle.
\newblock Entropy formulation for parabolic degenerate equations with general
  {D}irichlet boundary conditions and application to the convergence of {FV}
  methods.
\newblock  {\em SIAM J. Numer. Anal.}  41(6):2262--2293, 2003

\bibitem{Ohlberger:FVM}
M.~Ohlberger.
\newblock A posteriori error estimates for vertex centreed finite volume
  approximations of convection-diffusion-reaction equations.
\newblock {\em M2AN Math. Model. Numer. Anal.} 35(2):355--387, 2001.

\bibitem{Otto:BC}
F.~Otto.
\newblock Initial-boundary value problem for a scalar conservation law.
\newblock {\em C. R. Acad. Sci. Paris S\'er. I Math.} 322(8):729--734, 1996.

\bibitem{Otto:L1_Contr}
F.~Otto.
\newblock ${L}\sp 1$-contraction and uniqueness for quasilinear
  elliptic-parabolic equations.
\newblock {\em J. Differential Equations} 131(1):20--38, 1996.

\bibitem{Panov-Izv} E.Yu.~Panov. \newblock{On the theory of generalized entropy solutions of the Cauchy problem
for a first-order quasilinear equation in the class of locally
integrable functions.} (Russian) \newblock{Izvestiya Math.}
66(6):1171--1218, 2002.

\bibitem{Pierre}
Ch.~Pierre.\newblock{\em Mod\'elisation et simulation de l'activit\'e \'electrique du coeur
dans le thorax, analyse num\'erique et m\'ethodes de volumes
finis.} \newblock Ph.D. Thesis, Universit\'e de Nantes, 2005.

\bibitem{RouvreGagneux}
{\'E}.~Rouvre and G.~Gagneux.
\newblock Solution forte entropique de lois scalaires
  hyperboliques-paraboliques d\'eg\'en\'er\'ees.
\newblock {\em C. R. Acad. Sci. Paris S\'er. I Math.} 329(7):599--602, 1999.

\bibitem{Simondon}
F.~Simondon.
\newblock \'Etude de l'\'equation $\ptl_t b(u)-{\rm
div}\,a(b(u),\,\nabla u)=0$ [\dots] par la m\'ethode des
semi-groupes dans $L\sp{1}(\Omega ).$ (French) \newblock {\em
 Publ. Math.
Fac. Sci. Besan\c{c}on} 7, Univ. Franche-Comt\'e, Besan\c{c}on, 1983.

\bibitem{SougPerth:03}
P.~E. Souganidis and B.~Perthame.
\newblock Dissipative and entropy solutions to non-isotropic degenerate
  parabolic balance laws.
\newblock  {\em Arch. Ration. Mech. Anal.},  170(3):359--370, 2003; Addendum in
{\em Arch. Ration. Mech. Anal.} 174,no.3:443--447, 2004.

\bibitem{Tassa}
T.~Tassa.
\newblock Uniqueness of piecewise smooth weak solutions of multidimensional
  degenerate parabolic equations.
\newblock {\em J. Math. Anal. Applic.} 210(2):598--608, 1997.

\bibitem{Volpert}
A.~I. Vol'pert.
\newblock The spaces {BV} and quasi-linear equations.
\newblock {\em Math. USSR Sbornik}, 2(2):225--267, 1967.

\bibitem{VolHud}
A.~I. Vol'pert and S.~I. Hudjaev.
\newblock Cauchy's problem for degenerate second order quasilinear parabolic
  equations.
\newblock {\em Math. USSR Sbornik} 7(3):365--387, 1969.

\bibitem{Vovelle:02}
J.~Vovelle.
\newblock Convergence of finite volume monotone schemes for scalar conservation
  laws on bounded domains.
\newblock {\em Numer. Math.} 90(3):563--596, 2002.

\bibitem{WuYin}
Z.~Wu and J.~Yin.
\newblock Some properties of functions in {$BV_x$} and their applications to
  the uniqueness of solutions for degenerate quasilinear parabolic equations.
\newblock {\em Northeastern Math. J.} 5(4):395--422, 1989.

\bibitem{WZYL:2001}
Z.~Wu, J.~Zhao, J.~Yin and H.~Li.
\newblock {\em Nonlinear diffusion equations}.
\newblock World Scientific Publishing Co. Inc., River Edge, NJ, 2001.

\bibitem{Yin:Double}
J.~Yin.
\newblock On a class of quasilinear parabolic equations of second order with
  double-degeneracy.
\newblock {\em J. Partial Differential Equations} 3(4):49--64, 1990.

\end{thebibliography}
\end{document}